\documentclass{article}
\usepackage[margin=1.206in]{geometry} 

\usepackage[latin1]{inputenc}
\usepackage{amsmath}
\usepackage{amssymb,url}
\usepackage{mathrsfs}
\usepackage{amscd}%
\usepackage{graphicx}
\usepackage{subfigure}
\usepackage{multirow}
\usepackage{enumerate}
\usepackage{scrpage2}
\usepackage{amscd}
\usepackage{theoremref}
\usepackage{subfigure}
\usepackage{float}

\DeclareMathOperator{\Ran}{Ran}

\DeclareMathOperator*{\argmin}{arg\,min}

\newtheorem{theorem}{Theorem}[section]
\newtheorem{lemma}[theorem]{Lemma}
\newtheorem{proposition}[theorem]{Proposition}
\newtheorem{remark}[theorem]{Remark}

\newtheorem{definition}[theorem]{Definition}

\usepackage[numbers]{natbib}
\usepackage[boxed]{algorithm2e}

\usepackage[usenames]{color}

\title{Nonlinear Residual Minimization by Iteratively Reweighted Least Squares}

\author{Juliane Sigl \vspace{2mm}\\
\small{\it Technische Universit\"at M\"unchen}\\
\small Fakult\"at f\"ur Mathematik \\
\small Boltzmannstr.~3\\
\small85748 Garching bei M\"unchen\\
\small {\it juliane.sigl@ma.tum.de}\\
\small{+49 89 289 17989}
}

\begin{document}

\maketitle

\begin{abstract}

In this paper we address the numerical solution of minimal norm residuals of {\it nonlinear} equations in finite dimensions.  We take particularly inspiration from the problem of finding a sparse vector solution of phase retrieval problems by using greedy algorithms based on iterative residual minimizations in  the $\ell_p$-norm, for $1 \leq p \leq 2$. Due to the mild smoothness of the problem, especially for $p \to 1$, we develop and analyze a generalized version of Iteratively Reweighted Least Squares (IRLS). This simple and efficient algorithm performs the solution of optimization problems involving non-quadratic possibly non-convex and non-smooth cost functions, which can be transformed into a sequence of common least squares problems, to be tackled eventually by more efficient numerical optimization methods. While its analysis has been by now developed in many different contexts (e.g., for sparse vector, low-rank matrix optimization, and for the solution of PDE involving $p$-Laplacians) when the model equation is {\it linear}, no results are up to now provided in case of {\it nonlinear} ones. We address here precisely the convergence and the rate of error decay of IRLS for such nonlinear problems. 
The analysis of the convergence of the algorithm is based on its reformulation as an alternating minimization of an energy functional, whose main variables are the competitors to solutions of the intermediate reweighted least squares problems and their weights. Under a specific condition of coercivity often verified in practice and assumptions of local convexity, we are able to show convergence of IRLS to minimizers of the nonlinear residual problem. For the case where we are lacking the local convexity, we propose an appropriate convexification by quadratic perturbations.  Eventually we are able to show convergence of this modified procedure to at least a very good approximation of stationary points of the original problem.
In order to illustrate the theoretical results we conclude the paper with several numerical experiments. We compare IRLS with standard Matlab optimization functions for a simple and easily presentable example and furthermore numerically validate our theoretical results in the more complicated framework of phase retrieval problems, which are our main motivation. Finally we examine the recovery capability of the algorithm in the context of data corrupted by impulsive noise where the sparsification of the residual is desired.
\end{abstract}


%
\pagestyle{myheadings}
\thispagestyle{plain}
\markboth{
}{
}
\section{Introduction}

\subsection{ Iteratively Reweighted Least Squares}

In this paper we are addressing the numerical solution of nonlinear equations $A(x)=y$ where $A$ is a nonlinear mapping from $\mathbb R^k$ to $\mathbb R^m$ for $m \geq k$, and $y \in \mathbb R^m$ is a given datum. If $y \notin \Ran(A)$ the equation  $A(x)=y$ has no solution, but nevertheless often an approximate one with minimal residual in a certain norm is desired. A common choice for such a norm, especially in statistical and data analytical applications, is an $\ell_p$-norm for $1\leq p \leq 2$, depending on the preferred type of residual error \cite{Boyd04} modeling different kinds of noise e.g. Gaussian, Poisson, impulsive etc. and mixed types.  This leads to the $\ell_p$-norm  minimization of the residual of the equation in the Euclidean space $\mathbb R^k$,
\begin{equation}\label{lp}
\min\limits_{x \in \mathbb{R}^k} \left\| A(x)-y\right\|_{\ell_p^m}^p.
\end{equation}
As a motivation of this paper, the model \eqref{lp} appears naturally as an intermediate step used in greedy algorithms for the solution of nonlinear equations with sparse solutions \cite{Quasi14}. One typical example is the quadratic map $A(x) = (|\langle x, a_i \rangle|^2)_{i=1,\dots,m}$ encoding the amplitudes of the scalar products of a vector $x$ with respect to a given family of vectors $\{a_1,\dots,a_m\}$. In this case, the solution to the equation $A(x)=y$ eventually boils down to the recovery of the unknown signs of the scalar products, as a prototype of the more complex phase retrieval problem occurring in X-ray crystallography \cite{Drenth07,Fienup82,Gerchberg72}. (We shall use this particular application as a nontrivial test case for numerical experiments in Section 5.2.2)
\\

For $A$ smooth and $1\ll p \leq 2$, the objective residual function in \eqref{lp} is also smooth enough for the employment of a standard Newton method. However, these algorithms are usually guaranteed to converge only locally and for nonsmooth maps $A$ or $p \approx 1$ (or even $p=1$) one may have to use less efficient versions, such as the semi-smooth Newton method \cite{Ulbrich02}. As the objective function involves an $\ell_p$-norm, one can consider also other methods, which might better and more directly exploit the structure of the problem. In particular  Iteratively Reweighted Least Squares (IRLS) is a popular minimization strategy for optimization problems involving
non-quadratic possibly non-convex and non-smooth cost functions, which can be transformed into a sequence of common least squares problems, which can eventually be tackled by more efficient numerical optimization methods. 
\\
In our setting, the method is realized by substituting \eqref{lp} with a sequence of weighted quadratic problems 
 \begin{equation}\label{l2w}
x^{n+1}= \arg \min_{x \in \mathbb{R}^k} \left\| A(x)-y\right\|_{\ell^2(w^n)}^2,
\end{equation} where $\|\xi \|_{\ell_{2}(w^n)} = \left (\sum_{i=1}^m |\xi_i|^2 w^n_i \right)^{1/2}$ is a weighted $\ell_2$-norm, with a weight sequence $w^n_i \approx |(A(x^*)-y)_i|^{p-2}$, $i=1,\dots,m$. Here $x^*$ is the expected minimal solution of \eqref{lp}, which is of course not at hand. Therefore, one uses the practical iterative update rule $w^n_i = |(A(x^n)-y)_i|^{p-2}$, $i=1,\dots,m$, hoping for the realization of a contraction principle, which may eventually allow for the convergence of the iterates $x^n \to x^*$ for $n \to \infty$.
Notice now that, for $A$ smooth enough, the sequence of problems \eqref{l2w} can be addressed by efficient and standard  Newton methods, despite the fact that we are targeting a nonsmooth problem, for instance for $p=1$.
\\

The simplicity, adaptability, and its straightforward implementation make IRLS a very popular choice for beginners and first numerical test experiments. Besides it turns out to be extremely efficient in several contexts (sparse vector \cite{Daubechies10}, low-rank matrix \cite{Fornasier11}, bounded variation function  \cite{Chambolle97} solutions of minimal problems) and it can exhibit superlinear convergence also for nonsmooth optimization problems, see \cite{DaubechiesDFG08}. However, most of the known results of convergence of this algorithm are limited to the case where $A$ is a {\it linear} map. 
Let us now make a short account of the known results concerning the analysis of IRLS for the case where linear maps $A$ are involved as in \eqref{lp} or as a linear constraint in \eqref{lpminCS}.  
The first studies on iteratively reweighted least squares can be documented already in the 1960s. 
One of the first appearances can be found  in the approximation practice in the doctoral thesis of Lawson in 1961 \cite{Lawson61}, in
the form of an algorithm for solving uniform approximation problems, in particular by Chebyshev polynomials, by means of limits of weighted $\ell_p$-norm solutions.
This iterative algorithm is now well-known in classical approximation theory as Lawson's algorithm.  In \cite{Cline72} it is proved that this algorithm has in principle a linear convergence rate. In the 1970s extensions of Lawson's algorithm
for $\ell_p$-minimization, and in particular $\ell_1$-minimization, were proposed. Perhaps the most comprehensive mathematical analysis of the
performance of IRLS for $\ell_p$-minimization for $1<p<3$ was given in the work of Osborne \cite{Osborne852}. For quite a while, no significant developments have been reported, until the growing interest in the minimization of total variation regularized functionals in image processing introduced by \cite{Rudin92} moved back IRLS into the focus of the community in the early 1990s. Beside its very simple application to total variation minimization and its intuitive implementation, as presented in  \cite{Chambolle97}, IRLS allowed also for very efficient preconditioning strategies \cite{Vogel98} and showed to be more efficient than more generic optimization algorithms such as interior point methods. In the late 90s IRLS appeared as a method for the reconstruction of sparse signals in the papers \cite{Gorodnitsky97} and \cite{Chen98}, long before {\it compressed sensing} started to grow in popularity, after the pioneering work by Candes, Romberg, Tao and Donoho \cite{Candes062,Donoho062}.
In several papers \cite{Chartrand07,Chartrand08,Chartrand082,Daubechies10} a rigorous analysis of the behaviour of IRLS was carried out, for the situation where it is applied on the $\ell_p$-norm minimization problem with linear constraints of the form 
\begin{equation}
\label{lpminCS}
\min\limits_{Ax=y}\|x\|_{\ell_p}
\end{equation}
where $0 < p \leq 1,A \in \mathbb{R}^{m\times N}$ is a given matrix, and $y\in \mathbb{R}^m$ a given measurement vector. 
In  \cite{Fornasier11} a further extension of IRLS to the problem of low-rank matrix recovery from a minimal number of linear measurements has been developed and analyzed. Inspired by the approach in \cite{Chambolle97} for the solution of the total variation minimization problem, IRLS appears again in \cite{garau11} under the name of Ka\v{c}anov iteration, for the solution of quasi-linear elliptic equations.
\\

Especially since the beginning of this decade we witness a booming research on variations on the theme around IRLS, with new applications in statistics and signal processing, and it becomes hard to give a complete overview beyond the aforementioned milestones of the development up to present\footnote{On Google Scholar the appearance of the phrase 'Iteratively Reweighted Least Squares' was counted 3360 times in papers published since 2010, and more than 112 in their title since 1970, half of them after 2003.}. We may want to refer the reader to the paper \cite{Ochs14} and the reference therein, for an overview of the most recent literature. \\ 

\subsection{The contribution of this paper}

Often models of physical measurements in the applied sciences and engineering, however, are not linear but in practice linear models are assumed for simplicity and nonlinearities are neglected. Unfortunately linearization is not appropriate in many applications and the assumed model does not represent reality in a satisfactory way. A typical example is the phase retrieval problem mentioned above, which we shall use later as a test case for numerical experiments.\\
Therefore it is of utmost interest to investigate to which extent the analysis of convergence of IRLS where linear models are involved can be generalized to nonlinear ones, in particular we shall deduce conditions for its applicability and we state its limitations as well. More precisely our aim is to study the numerical approximation via IRLS of a solution of the $\ell_p$-norm-minimization problem \eqref{lp}, where $A$ might be nonlinear and mildly smooth, and $1 \leq p < 2$ (the case $p=2$ is just a least squares and there is no need of further iterations). Notice that we are not afraid here to include the case $p=1$.\\

As already shown explicitly in \eqref{l2w}, the extension of the IRLS to this type of problems from an implementation point of view is simply straightforward. Hence, for practical issues, there are no additional difficulties beyond the application of standard recipes (including possible preconditioning etc.). Nevertheless, to our knowledge, a rigorous theoretical analysis of the convergence of IRLS for nonlinear residual minimizations as in \eqref{lp} has not been done yet, especially for the cases where $A$ is mildly smooth and $p \approx 1$. It is the subject of the present paper.\\


In \cite{Chambolle97,Daubechies10,Fornasier11} the analysis of convergence of IRLS has been based on its variational formulation, involving an energy functional where both the unknown and the weights appear as minimizing variables, and the algorithm is re-interpreted as an alternating minimization over them. For the problems in \cite{Daubechies10,Fornasier11} a certain coercivity condition was additionally required, namely the Restricted Isometry Property (RIP),  which is a near-identity spectral property of small submatrices of the linear model $A$. We shall extend  the analysis done in the aforementioned papers, especially we take inspiration from \cite{Daubechies10}, and analyze the generalized version of IRLS for $\ell_p$-norm minimization of the residual as in \eqref{lp}, by requiring a relaxed version of the RIP as already introduced in our previous paper  \cite{Quasi14}.
First of all we start by introducing a similar energy functional as the ones proposed in \cite{Chambolle97,Daubechies10,Fornasier11}, more precisely of the form 
\begin{equation}
\mathcal{J}(x,w,\epsilon):=\frac{p}{2}\left[ \sum_{i=1}^m w_i(A_i(x)-y_i)^2+\sum_{i=1}^m\left(\epsilon^2w_i+\frac{2-p}{p}w_i^{p/(p-2)}\right)\right], \text{   }x\in \mathbb{R}^k
\end{equation}
with $\epsilon > 0$, and  weight vector $w \in \mathbb{R}^m$, with positive entries $w_i > 0, i = 1, \dots,m$ and $1\leq p<2$.

Under the mentioned coercivity assumption we prove the convergence  and corresponding error decay rates of the IRLS type algorithm based on the alternating minimization of $\mathcal{J}(x,w,\epsilon)$, whenever $x \to  \mathcal{J}(x,w,\epsilon)$ is locally convex.
For the case where we are lacking the local convexity, we propose an appropriate convexification by quadratic perturbations. Let us remark that this strategy of convexification is rather standard and well-known in the nonlinear optimization literature, for instance in sequential quadratic programming \cite{Boyd04}. The innovation here is in addressing problems with severe nonsmoothness due to the possible choice of $p=1$, as already considered, e.g., in \cite{Artina13}. Eventually we are able to show convergence of this modified procedure to at least a very good approximation of stationary points of the original problem.\\

In order to explain and better illustrate our theoretical results we conclude the paper with several numerical experiments. 
Comparisons with standard Matlab methods applied to the original $\ell_p$-minimization problem for a simple easily presentable example reveals that IRLS possibly converges to different local minimizers than standard methods when starting computations from the same initial point. Furthermore we numerically validate our theoretical results for the more complex task of finding a sparse solution of phase retrieval problems, and we check success via the correct reconstruction of sparse vectors. As we will see in this more sophisticated application, IRLS significantly outperforms standard methods. Finally we examine the recovery capability of the algorithm in the context of data corrupted by impulsive noise where the sparsification of the residual is desired.

\subsection{Outline of the paper}

The paper is organized as follows: In Section 2, we introduce definitions and notations, show the reformulation of the $\ell_p$-minimization into a reweighted $\ell_2$-least squares problem and give a short review of popular numerical methods for its solution. Finally we present the IRLS method tailored to problems of the type \eqref{lp} and in Section 3 its detailed analysis of the convergence and rate of convergence follow. Thereafter in Section 4 the method is extended to cases where we can not guarantee local convexity and convergence is deduced for the modified algorithm.  This work is concluded by numerical experiments in Section 5.

\section{Nonlinear Residual Minimization via Iteratively Reweighted Least Squares}

In this section, we introduce the main terms and notations used in this paper and point out how the problem stated in \eqref{lp} can be recast as an approximating nonlinear weighted least squares problem. In addition, we
shortly review the basics around nonlinear weighted least squares minimization and give a very short overview over practical, efficient algorithms.

\begin{definition}(Weighted $\ell_p$-spaces)\\
Define the Banach space $\ell_p^m(w):=(\mathbb{R}^m,\|\cdot\|_{\ell_p^m(w)}) $ endowed with the weighted norm
$$\|x\|_{\ell_p^m(w)}:=\left(\sum\limits_{i=1}^m w_i|x_i|^p\right)^{\frac{1}{p}}.$$
When the dimension $m$ is understood, we omit it in the symbol of the space $\ell_p = \ell_p^m$.
Furthermore we define the (unweighted) $\ell_p$-spaces by setting $\ell_p(w):=\ell_p(1)$, where $1$ here is the weight with entries identically set to $1$.\\
For a nonlinear map $A:\ell_p^k\rightarrow \ell_q^m$, we define the norm
$$\|A\|_{\ell_p^k\rightarrow \ell_q^m}:=\sup\limits_{\|x\|_{\ell_p^k}\leq 1}\|A(x)\|_{\ell_q^m}$$
and for the particular case of $p=q=2$, $\|A\|:=\|A\|_{\ell_p^k\rightarrow \ell_q^m}$ is the standard operator norm.
\end{definition}

The unit balls are indicated by $\mathcal{B}^{\ell_p}:=\left\{x \in \mathbb{R}^m:\left\|x\right\|_{\ell_p} \leq 1 \right\}$. More generally, the shifted balls centered around $\bar{x} \in \mathbb{R}^m$ with radius R are denoted by $\mathcal{B}^{\ell_p}(\bar{x},R):=\left\{x \in \mathbb{R}^m:\left\|x-\bar{x}\right\|_{\ell_p} \leq R \right\}$. In cases where the space is clear we simplify the notation and use $\mathcal{B}(\bar{x},R).$ \\

Furthermore we will consider the range of a map $A:\ell_p^k\rightarrow \ell_q^m$ and denote it by $\Ran(A)$. \\

The space of $n$-times continuously differentiable functions from a certain Euclidean space $\mathcal{E}_1$ to another Euclidean space $\mathcal{E}_2$ will be denoted by $\mathcal{C}^n(\mathcal{E}_1,\mathcal{E}_2)$.
By $I:\mathcal{E}_1\rightarrow \mathcal{E}_1$ we express the identity operator on $ \mathcal{E}_1$. (The spaces $\mathcal{E}_1,\mathcal{E}_2$ will be clear from the context.)\\

Positive, independent constants are in general named $ C,\tilde{C} ,\hat{C},\bar{C},C^*,C_1, C_2\ldots$

\subsection{Problem Setting and Characterization of $\ell_p$- and reweighted $\ell_2$-minimizers }
Now we are well equipped to formulate the problem setting.

Let $A:\mathbb{R}^k\rightarrow\mathbb{R}^m$ be a {\it nonlinear continuous map} with $k\geq m$ and $y \in \mathbb{R}^m$. In the following we fix $y$ and consider the nonlinear equation system 
$$A(x)=y,$$
which is an overdetermined system that has no solution if $y \notin \Ran(A)$. Nevertheless often an approximate solution with minimal residual in a certain norm, commonly an $\ell_p$-norm, is desired.\\
Our aim is to find a best approximating $x^*$, given the measured data $y$ as the solution to the $\ell_p$-norm-minimization problem
\begin{equation}\label{ellp}
\left \| A(x^*)-y\right\|_p^p=\min\limits_{x \in \mathbb{R}^k}\left\| A(x)-y\right\|_p^p,
\end{equation}
where $1\leq p\leq 2$.\\
We next consider the minimization in a weighted $\ell_2(w)$-norm. We suppose that the weight vector $w=(w_1,\cdots,w_m)$ is strictly positive, meaning that $w_i > 0$ for all $i\in \left\{1,\dots,m\right\}$. \\

We formally observe now that  
\begin{equation*}
 \left\| A(x)-y\right\|_{\ell_p}^p=\sum\limits_{i=1}^m(A_i(x)-y_i)^p= \sum\limits_{i=1}^m\frac{(A_i(x)-y_i)^2}{(A_i(x)-y_i)^{2-p}}.
\end{equation*}
and setting $w_i=|A_i(x)-y_i|^{p-2}$ we can formally write
$$ \left\| A(x)-y\right\|_{\ell_p}^p=\left\| A(x)-y\right\|_{\ell_2(w)}^2.$$
The minimizer with respect to this weighted norm will be denoted as
$$x^w=\argmin\limits_{x \in \mathbb{R}^k}\|A(x)-y\|_{\ell_2(w)}^2.$$

In the next subsection we will shortly summarize common techniques for solving problems of this type concentrating in particular on Newton type methods with approximations to the Hessian. 
\subsection{Weighted Nonlinear Least Squares Fitting} 
Nonlinear unconstrained optimization 
$$\min\limits_{x \in \mathbb{R}^k} f(x)$$
is a very common problem, which we specify here for $f(x)= \sum\limits_{i=1}^m w_i f_i(x)_i$ and $f_i(x)=(A_i(x)-y_i)^2$.\\
For solving such a problem a broad variety of standard tools is available. On the one hand  we can resort to algorithms based only on function evaluation values such as the Nelder-Mead algorithm \cite{Nelder65} or pattern search \cite{Torczon97} that do not require smoothness or methods involving information on the gradient or Hessian of the objective function such as gradient descent methods or Newton's method (depending on the smoothness of $A$). \\
On the other hand we are dealing with a very specific type of unconstrained problem, whose specific structure, in particular the 
special structure of the Hessian matrix for the least-squares objective function, can be exploited in more specialized algorithms tailored to this type of regression problem assuming $A$ is smooth enough. \\
In this case the Hessian can be split in two terms where one represents the linearized part of the problem and the other bringing in the nonlinearity, which is the critical one:
$$ \nabla^2 f(x)=\sum_{i=1}^m w_i \left[ \nabla f_i(x) \nabla f_i(x)^T+f_i(x)\nabla^2 f_i(x) \right]=\sum_{i=1}^m w_i \left[ \nabla A_i(x) \nabla A_i(x)^T+(A_i(x)-y_i)\nabla^2 A_i(x) \right].$$
 The first term involves the gradients of the functions only and therefore computations are easier to perform, while the second contains computationally expensive second order derivatives, but is very small if the data fit is already good and the residual is close to zero. 
Many methods specialized on nonlinear least squares are based on the idea of computing the first term of the Hessian exactly while the second is only approximated using first derivatives. \\

The simplest realization of this idea is the Gauss-Newton method, which iteratively computes the search direction $d$ via the formula for Newton's method 
$$\nabla^ 2f(x) d=-\nabla f(x) $$ and approximates the second order derivative $\nabla^ 2f(x)$ by only the first term, resulting in the system  
$$\nabla A(x) W \nabla A (x)^Td =-\nabla A (x) W A(x),$$
where $W$ is the diagonal matrix with the weights $w_i$ on its diagonal.\\ 
Once the search direction is computed one produces an updated approximation by setting 
$$x\rightarrow x+d$$
and so on. \\
If  at the solution vector $x^*$ obtained thereby it holds $f(x^*)=0$ and $\nabla A(x^*)$ is of full rank, the Gauss-Newton method behaves like Newton's method close to the solution without investing the effort of computing second order derivatives. A more detailed analysis can be found in \cite{Dennis77}.\\
A great number of related methods for nonlinear least-squares can be seen as improvements of the Gauss-Newton algorithm because they involve a kind of approximation to the second term in the formula for the Hessian matrix
$$\sum\limits_{i=1}^mf_i(x) \nabla^2 f_i(x)=\sum\limits_{i=1}^m(A_i(x)-y_i) \nabla^2 A_i(x).$$
The most well-known is also one of the simplest, referred to as the Levenberg-Marquardt algorithm, appearing first in the papers of Levenberg \cite{Levenberg44} and Marquardt \cite{Marquardt63}. This strategy uses the approximation 
$$\sum\limits_{i=1}^mw_if_i(x) \nabla^2 f_i(x)=\sum\limits_{i=1}^mw_i(A_i(x)-y_i) \nabla^2 A_i(x) \approx \lambda W,$$
where $\lambda\geq0$  is some scalar. Then the search direction is obtained by solving the linear system
$$\left[\nabla A(x) W \nabla A (x)^T+\lambda W\right]d=-\nabla A (x) W A(x).$$
For the implementation of the Levenberg-Marquardt it is possible to use a trust-region strategy as described in \cite{More77}. 
Other approximations to the Hessian of $f(x)$ are also possible, for example, a quasi-Newton approximation to the second term of the Hessian but we will close the excursion on this type of methods here and refer to related, standard literature \cite{Ruhe79,Gill78,Dennis81,Bates88}.
\begin{remark}
 If this approximation of the Hessian in the methods mentioned above is not accurate enough these strategies have rate of convergence that can not compete with Newton's method. More precisely the convergence rate will be at most linear. Hence, the life becomes difficult with these methods if, for any reason, the problem is not smooth enough!
\end{remark}
Besides these methods approximating the Hessian other classes of methods exist e.g. based on orthogonal decomposition of the Jacobian that shall not be covered here and for further information we refer to the literature \cite{Fraley88} .
\begin{remark}\thlabel{nonc2} At this point we note that the methods mentioned above only converge to stationary points or at best local minimizers in the context of nonconvex least squares problems! Global optimization in the nonconvex case is in general a challenging task. Nevertheless  under certain smoothness conditions, \cite{Zhou13} suggests a method involving a line-search strategy that guarantees to find the global minimizer.
\end{remark}

\subsection{Auxiliary Functional and Nonlinear Residual Iteratively Reweighted Least Squares Algorithm}
At this point we would like to introduce a helpful tool for the formulation and theoretical analysis of an iteratively reweighted least squares algorithm for problems of type \eqref{lp} in the form of the following functional:
\begin{definition}\thlabel{J}
 Given a real number $\epsilon > 0$, $ x \in \mathbb{R}^m$, and a weight vector $w \in \mathbb{R}^m$, with positive entries
$w_i > 0, i = 1, \dots,m$ and $1\leq p <2$ we define
\begin{equation}
\mathcal{J}(x,w,\epsilon):=\left[ \sum_{i=1}^m w_i(A_i(x)-y_i)^2+\sum_{i=1}^m\left(\epsilon^2w_i+\frac{2-p}{p}w_i^{p/(p-2)}\right)\right], \text{   }x\in \mathbb{R}^m.
\end{equation}
\end{definition}
A very similar auxiliary functional is appearing in \cite{Chambolle97,Daubechies10,Fornasier11} as well, but the one defined above is a version adapted to the setting of $\ell_p$-minimization of residuals. \\

IRLS actually performs an alternating minimization of the functional $\mathcal{J}$.

\begin{algorithm}
\KwSty{Algorithm 1: Nonlinear residual iteratively reweighted least squares (NR-IRLS)}:\\
\KwIn{$A:\mathbb{R}^k\rightarrow \mathbb{R}^{m}$, $y\in\mathbb{R}^m$ }
We initialize by taking $w^0:=(1,\cdots,1)$. We also set $\epsilon_0=1$. \\
We then recursively define for $n=0,1, \dots$
$$x^{n+1}=\argmin\limits_{x \in \mathbb{R}^k} \mathcal{J}(x,w^n,\epsilon_n)=\argmin\limits_{x \in \mathbb{R}^k} \left\|A(x)-y\right\|^2_{\ell_2(w^n)}$$
and
the auxiliary variables
$$\mathcal{N}^{n+1}=\min_i(|A_i(x^{n+1})-y_i)|) \text{  and  } \mathcal{M}^{n+1}=\max_i(|A_i(x^{n+1})-y_i|),$$
to define
$$\epsilon_{n+1}=\min\left( \max( \mathcal{N}^{n+1},\tilde{\epsilon}),\epsilon_n, \mathcal{M}^{n+1}\right)$$
where $\tilde{\epsilon}>0$ is a fixed small constant.
Then we define
$$w^{n+1}=\argmin\limits_{w \in \mathbb{R}^k_+}\mathcal{J}(x^{n+1},w,\epsilon_{n+1}) =  \left ( \left ( (A(x^{n+1})_i-y_i)^2 + \epsilon_{n+1}^2) \right)^\frac{p-2}{2} \right)_{i=1}^m.$$
\KwOut{$x^{(1)}$, $x^{(2)},\ldots$}
\label{algo:1}
\end{algorithm}

We stop the algorithm if $\epsilon_n = 0$. In this case we define $x^j := x^n$ for $j > n$. However, in general, the algorithm will generate an infinite sequence $(x^n)_{n \in \mathbb{N}}$ of distinct vectors and it is convenient to stop as soon as $\epsilon$ falls below an appropriately chosen threshold.\\
\begin{remark}\thlabel{nonc}\begin{enumerate}
\item{For $\epsilon, w$ fixed the functional $\mathcal J(\cdot, w,\epsilon)$ is not necessarily convex in the variable $x$ if $A$ is a nonlinear map.\\
In the case that $A$ is indeed linear and has full rank the functional $\mathcal J(\cdot, w,\epsilon)$ is strictly convex in the variable $x$ which implies that every stationary point is a global minimizer and this optimization problem can be solved efficiently. \\
In contrary nonlinearity of the map $A$ turns minimizing the functional $\mathcal J(\cdot, w,\epsilon)$ into a nonconvex $k$-dimensional minimization problem where several local minimizers and stationary points can occur and solving the problem becomes a harder task see \thref{nonc2}.\\
The reader is referred to the Appendix for further considerations on convexity for the linear case as well as nonlinear perturbations of linear random maps as in the numerical experiments in Section 5.}
\item{Each step of the algorithm requires the solution of a $k$-dimensional nonlinear weighted least squares problem, which can be solved numerically by the popular strategies mentioned above but we emphasize once more that in the general case those methods only find critical points, and this an intrinsic limitation of NR-IRLS!}
\end{enumerate}
\end{remark}

\section{Theoretical Analysis and Convergence Results for the NR-IRLS }
In the following section we will have a closer look at Algorithm 1 and point out some of its properties, in particular the boundedness of the iterates $(x^n)_{n\in \mathbb{N}}$ and the fact that these are getting arbitrarily close as $n\rightarrow \infty$. These results will be useful to finally develop the proof of convergence and to estimate the rate of convergence of Algorithm 1 under conditions determined along the way. \\

First of all our relevant domain for the search for the minimizer of  \eqref{ellp} shall be a suitably chosen, sufficiently large closed set $\mathcal{D}\in \mathbb{R}^k$ containing $x^*$ and $0$. We require that the first iterate vector $x^1$ is contained in $\mathcal{D}$. 

\begin{remark}\thlabel{x0}We stress again that without further assumptions on the nonlinear map $A$,  nonconvexity of the functional $\mathcal{J}(\cdot,w^0,\epsilon_0)$ can not be excluded and therefore several critical points are expected. Hence a different choice of the starting vector $x^0$ for your iterative method for solving the nonlinear least squares problem in the very first step will influence the behavior of the algorithm in the following iterations and hence the final solution of Algorithm 1 as well!
\end{remark}

At this point we want to state some essential assumptions on the measurement map $A$. We require $A \in \mathcal{C}^0(\mathbb{R}^k,\mathbb{R}^m)$ and that $A$ is bounded on $\mathcal{D}$.
Additionally we introduce  the following property, which we call  {\it boundedness and coercivity condition (BCC)}
\begin{definition}\thlabel{BCC}
Let $A:\mathbb{R}^k\rightarrow\mathbb{R}^m$ be a nonlinear continuous map. We say that $A$ fulfills the boundedness and coercivity condition (BCC) at $x^* \in \mathcal D$ if there exist $\alpha,\beta>0$ such that
$$\alpha \|x^*-z\|_{\ell_2}\leq \|A(x^*)-A(z)\|_{\ell_p}\leq \beta \|x^*-z\|_{\ell_2}$$
for all $z\in\mathcal{D}$.
\end{definition}
\begin{remark}
The lower bound in the BCC imposes that the level set of $x^*$ is a singleton only containing $x^*$ itself. This is necessary to conclude the identifiability from the nonlinear measurements $A(x^*)$ without any further assumptions on $x^*$ itself. The upper bound on the other hand is the request of Lipschitz continuity at $x^*$.
\end{remark}

Now we want to return back to the functional in \thref{J} and our first quite straightforward observation is that after the $n$-th step we obtain
$$\mathcal{J}(x^{n+1},w^{n+1},\epsilon_{n+1})=\sum\limits_{i=1}^m[(A(x^{n+1})_i-y_i)^2+\epsilon_{n+1}^2]^{p/2}.$$

Moreover due to the minimization properties resulting from Algorithm 1, the following monotonicity property holds.
\begin{lemma}\thlabel{monotonicity}
The inequalities 
$$\mathcal{J}(x^{n+1},w^{n+1},\epsilon_{n+1})\leq \mathcal{J}(x^{n+1},w^{n},\epsilon_{n+1})\leq \mathcal{J}(x^{n+1},w^{n},\epsilon_{n})\leq\mathcal{J}(x^{n},w^{n},\epsilon_{n}) $$
hold for all $n\geq0$.
\end{lemma}
Proof: 
Here the first inequality follows from the minimization property that defines $w^{n+1}$, the second
inequality from $\epsilon_{n+1}\leq \epsilon_n$, and the last inequality from the minimization property that defines
$x^{n+1}$.$\blacksquare$\\

Due to \thref{monotonicity} we can state that $\mathcal{J}(x^n,w^n,\epsilon_n)\leq \mathcal{J}(x^1,w^0,\epsilon_0)$, where the right hand side is a constant and this will help to obtain the boundedness of the iterates $(x^n)_{n\in \mathbb{N}}$:

\begin{lemma} \thlabel{xbound}
Let $A:\mathbb{R}^k\rightarrow\mathbb{R}^m$ be a nonlinear continuous map fulfilling the boundedness and coercivity condition (BCC) at $x^* \in \mathcal D$  and $y\in \mathbb{R}^m$ be given. Then 
the sequence of iterates $\left(x^n\right)_{n}$ defined by Algorithm 1 is bounded and hence lies in the ball $\mathcal{B}(0, R^*)$, where $R^*=\frac{1}{\alpha } \mathcal{J}(x^{0},w^{0},\epsilon_{0})^{1/p}+\frac{1}{\alpha}\|A(x^*)-y\|_{\ell_p}+\|x^*\|_{\ell_2}$.
\end{lemma}
Proof: For all $n \in \mathbb{N}$
\begin{align*}
\|x^n\|_{\ell_2}&\leq\|x^n-x^*\|_{\ell_2}+\|x^*\|_{\ell_2}\leq\frac{1}{\alpha} \|A(x^n)-A(x^*)\|_{\ell_p}+\|x^*\|_{\ell_2}\\
&\leq  \frac{1}{\alpha }\left(\sum\limits_{i=1}^m[(A(x^n)_i-y_i)^2+\epsilon_n^2]^{p/2}\right)^{1/p}+\frac{1}{\alpha}\|A(x^*)-y\|_{\ell_p}+\|x^*\|_{\ell_2}\\
&\leq \frac{1}{\alpha } \mathcal{J}(x^{n},w^{n},\epsilon_{n})^{1/p}+\frac{1}{\alpha}\|A(x^*)-y\|_{\ell_p}+\|x^*\|_{\ell_2}.
\end{align*}
By the monotonicity property in \thref{monotonicity} we obtain
$$\|x^n\|_{\ell_2}\leq \frac{1}{\alpha } \mathcal{J}(x^{1},w^{0},\epsilon_{0})^{1/p}+\frac{1}{\alpha}\|A(x^*)-y\| _{\ell_p}+\|x^*\|_{\ell_2}=R^*,$$
where all terms composing on the right hand side are bounded.$\blacksquare$\\

\begin{remark} Hence the iterates are contained in a ball of radius $R^*$ and
although the minimization in the first step of the algorithm is a global minimization it actually turns out to be local on this specific ball $\mathcal{B}(0,R^*)$. Nevertheless this ball can be quite large especially when we have a very small BCC constant $\alpha$.\\
\end{remark}

As $R^*>0$ is depending on the unknown solution $x^*$ we would also like to have an estimate only depending on values given to us already or that are fixed in advance. On the way towards such a sovereign upper bound on $R^*$, we observe that due to the assumed minimality of $x^*$
$$\|A(x^*)-y\|_{\ell_p} \leq \|A(0)-y\|_{\ell_p}.$$
Moreover 
\begin{align*}
\|x^*\|_{\ell_2}&=\|x^*-0\|_{\ell_2} \leq \frac{1}{\alpha} \|A(x^*)-A(0)\|_{\ell_p}\\
&\leq \frac{1}{\alpha} \left(\|A(x^*)-y\|_{\ell_p}+\|A(0)-y\|_{\ell_p}\right) \leq \frac{2}{\alpha} \|A(0)-y\|_{\ell_p}
\end{align*}
Hence we obtain the upper bound $\hat{R}$
$$R^*\leq\hat{R}:= \frac{1}{\alpha} \left(\mathcal{J}(x^0,w^0,\epsilon_0)^{1/p}+3\|A(0)-y\|_{\ell_p}\right).$$



As already indicated in \thref{nonc} the possible nonconvexity of functional $\mathcal J(\cdot, w,\epsilon)$ can severely impede the optimization process in the first step of Algorithm 1 and is more problematic for a theoretical study of the overall behavior of the algorithm as well.

For a complete analysis of the convergence we will consider the case where the functional $\mathcal J(\cdot, w,\epsilon)$ is locally convex case first and then continue with the introduction of appropriate adjustments of the algorithm in the subsequent section, for the case where the local convexity fails.

At this point we recall the concept of strong convexity of a function as in \cite{Ansari13}
\begin{definition}\thlabel{sconvex}
A continuous function $f: \mathbb{R}^k \rightarrow \mathbb{R}$ is called {\it C-strongly convex} at a point $\tilde{x}$ in a set $S\in\mathbb{R}^k$ if
\begin{equation}
f(tx+(1-t)\tilde{x}) \le t f(x)+(1-t)f(\tilde{x}) -  C t(1-t) \|x-\tilde{x}\|_{\ell_2}^2 \text{  for  }0\le t\le1, C>0\text{  for all  }x\in S.
\end{equation}
\end{definition}
\begin{remark}\thlabel{r1}
\begin{itemize}
\item[(i)]{If the function $f$ is continuously differentiable, the condition in \thref{sconvex} above can be formulated as follows
$$ f(x) - f(\tilde{x}) \ge \nabla f(\tilde{x})^T (x-\tilde{x}) + {C} \|x-\tilde{x}\|_{\ell_2}^2 \text{  for  }0\le t\le1, C>0\text{  for all  }x\in S.$$}
\item[(ii)]{If the function $f$ is twice continuously differentiable, the condition in \thref{sconvex} above can be deduced via Taylor's formula and the mean value theorem
$$ f(x) = f(\tilde{x}) +\nabla f(\tilde{x})^T (x-\tilde{x}) + \frac{1}{2} (x-\tilde{x})^T \nabla^2f(z) (x-\tilde{x})$$ 
for a certain  $z \in \{ t \tilde{x} + (1-t) x : t \in [0,1] \}.$ The matrix appearing in the last term can be characterized by an integral as follows
$$ \nabla^2f(z)=\int_0^1(1-t) \nabla^2 f(\hat{x}(t)) dt\text{ with }\hat{x}(t)=tx+(1-t)\tilde{x}.$$
Then the strong convexity constant $C$ equals the smallest eigenvalue of $ \nabla^2f(z)$ { as soon as it is strictly positive definite}.
}
\end{itemize}
\end{remark}
From strong convexity of a function at a minimizer we obtain the following result:
\begin{theorem}\thlabel{fc}
Let $f:\mathbb{R}^k\rightarrow\mathbb{R}$ be a $\mathcal {C}^{0}$ nonlinear map  and $x^*$ a minimizer. If $f$ is ${C}$-strongly convex at $x^*$ in the set $S$, then the following lower bound holds
\begin{align}\label{subs1}
f(x) - f(x^*) \ge {C} \|x-{x^*}\|^2_{\ell_2}   
\end{align}
where $x\in S$, ${C}>0$.
\end{theorem}

Proof:
By \thref{sconvex} we have that 
\begin{align*}
f(tx+(1-t)x^*) \le t f(x)+(1-t)f(x^*) -  {C} t(1-t) \|x-x^*\|_{\ell_2}^2,
\end{align*}
holds if and only if
$$
\frac{f(tx+(1-t)x^*)-f(x^*)}{t} +{C} (1-t) \|x-x^*\|_{\ell_2}^2 \le  f(x)-f(x^*). 
$$
As this condition has to hold for all $t \in [0,1]$, letting $t\rightarrow 0^+$ yields
$$f(x)- f(x^*) \ge \tilde{C}+{C} \|x-{x^*}\|_{\ell_2}^2 ,$$
where $ \tilde{C}=\liminf_{t\rightarrow {0}^+}\frac{f(tx+(1-t)x^*)- f(x^*)}{t}$. Due to the optimality of $x^*$ the term $\tilde{C}$ is always nonnegative, and
$$ f(x) - f(x^*) \ge  {C} \|x-{x^*}\|_{\ell_2}^2.$$

From now on we shall assume that the functional $\mathcal{J}(\cdot,w^n,\epsilon_n)$ is strongly convex as in \thref{sconvex} locally at $x^{n+1}$ for all $n\geq 0$ and we describe a uniform property as follows.
\begin{definition}\thlabel{USCC1}
Let $A:\mathbb{R}^k\rightarrow\mathbb{R}^m$ be a nonlinear map in ${\mathcal{C}^{0}}$ and $\mathcal{J}(x,w^n,\epsilon_n)$ for $w^n,\epsilon_n$ as generated by Algorithm \ref{algo:1}  for all $n\geq 0$, with the corresponding minimizer $x^{n+1}$. 
We say that the \emph{first uniform strong convexity condition (USCC-1)} is fulfilled if there exists a uniform constant $C>0$ such that for all $n\geq 0$ the following condition holds
\begin{align}\label{subs1}
&\mathcal{J}(x^n,w^n,\epsilon_n)- \mathcal{J}(x^{n+1},w^n,\epsilon_n) =\|A(x^n)-y\|^2_{\ell_2(w_n)}-\|A(x^{n+1})-y\|^2_{\ell_2(w_{n})}\ge {C} \|x^n-{x^{n+1}}\|_{\ell_2}^2
\end{align}
\end{definition}
\begin{remark} The USCC-1 is of course fulfilled if the functional $\mathcal{J}(\cdot,w^n,\epsilon_n)$ for $w^n,\epsilon_n$ fixed in each step $n$ is strongly convex at $x^{n+1}$ in the set $\mathcal{B}(0,R^*)$ with a constant $C>0$ independent of $n$. 
\end{remark}

Moreover we want to point out that a further desirable property of the map $A$ is the strong convexity of the function $x \to \|A(x)-y\|^2_{\ell_p}$ at its minimizer $x^*$. 
\begin{definition}\thlabel{USCC2}
Let $A:\mathbb{R}^k\rightarrow\mathbb{R}^m$ be a nonlinear map in ${\mathcal{C}^{0}}$ and $x^*$ a minimizer of $\|A(x)-y\|^2_{\ell_p}$
We say that the \emph{second uniform strong convexity condition (USCC-2)} is fulfilled if there exists a uniform constant $\hat{C}>0$ such that for all $n\geq0$ the following condition holds
\begin{align}\label{subs2}
\|A(x^n)-y\|^2_{\ell_p}- \|A(x^*)-y\|^2_{\ell_p}\ge \hat{C} \|x^n-{x^*}\|_{\ell_2}^2  .
\end{align}

\end{definition}

\begin{remark}
\begin{itemize}
\item[(i)]{In the case $A(x^*)=y$ the condition \eqref{subs2} is equivalent to the lower bound of the BCC with $\hat{C}=\alpha^2.$}

\item[(ii)]{Notice that condition \eqref{subs2} is fulfilled as soon as the function $\|A(\cdot)-y\|^2_{\ell_p}$ is totally convex at $x^*$ according to the definition of total convexity in \cite{Resmerita06}. Furthermore it follows from results in \cite{Resmerita06} that in this case the function $\|A(\cdot)-y\|^2_{\ell_p}$ is also strictly convex in $\mathcal{B}(0,2R^*)$.  }
\end{itemize}
\end{remark}

Now we are prepared with all necessary tools to formulate preliminary results on the convergence or NR-IRLS.
\subsection{Preliminary Results}
In this section we formulate several Lemmata that will be fundamental ingredients for the proof of convergence of Algorithm 1. \\
As a first result we would like to state that from the sequence $\mathcal{J}(x^n,w^n,\epsilon_n)$  being convergent it follows that the iterates $x^0, \cdots,x^n,x^{n+1}, \cdots \in \mathbb{R}^k$ of Algorithm 1 are getting arbitrarily close for $n\rightarrow \infty$ assuming that the USCC-1 is fulfilled.

\begin{lemma}\thlabel{cauchyx}
Let $A:\mathbb{R}^k\rightarrow\mathbb{R}^m$ be a nonlinear map in $\mathcal{C}^0$ and given $y\in \mathbb{R}^m$, if the USCC-1 property as in \thref{USCC1} is fulfilled with constant $C$, for the iterates of Algorithms 1 it holds
$$\lim_{n \rightarrow \infty} \|x^{n}-x^{n+1}\|_{\ell_2}^2 = 0. $$
\end{lemma}

Proof: For each $n=1,2,...$ we have 
\begin{align*}
\left[ \mathcal{J}(x^{n},w^n,\epsilon_{n})-\mathcal{J}(x^{n+1},w^{n+1},\epsilon_{n+1})\right]&\geq  \left[ \mathcal{J}(x^{n},w^n,\epsilon_{n})-\mathcal{J}(x^{n+1},w^{n},\epsilon_{n})\right]\\
&\geq C \left\| x^{n+1}-x^n\right\|^2_{\ell_2}
\end{align*}

From the monotonicity as in \thref{monotonicity} and the boundedness of the sequence $\left(\mathcal{J}(x^{n},w^n,\epsilon_{n})\right)_{n\in \mathbb{N}}$  we know that 
$$\lim_{n \rightarrow \infty} (\mathcal{J}(x^{n},w^n,\epsilon_{n})-\mathcal{J}(x^{n+1},w^{n+1},\epsilon_{n+1}))= 0,$$
hence also 
$$\lim_{n \rightarrow \infty} \|x^{n}-x^{n+1}\|_{\ell_2}^2 =0 .\text{   }\blacksquare$$\\

From the monotonicity of $\epsilon_n$, we know that $\epsilon:= \lim_{n\rightarrow \infty} \epsilon_n$ exists and is non-negative. The following functional will play a role in our proof of convergence, especially for $\epsilon>0$.
\begin{definition} ($\epsilon$-perturbed $\ell_p$-norm residual)
Let $A:\mathbb{R}^k\rightarrow\mathbb{R}^m$ be a nonlinear map and given $y\in \mathbb{R}^m$ we define the $\epsilon$-perturbed $\ell_p$-norm residual to be the following functional 
$$f_\epsilon(x):=\sum^m_{i=1} ( (A(x)_i-y_i)^2+\epsilon^2)^{p/2}.$$
\end{definition}
Notice that, if we knew that $x^{n}$ converged to a point $\bar{x}$ then, having in mind that 
\begin{equation}\label{JJ}
\mathcal{J}(x^{n},w^{n},\epsilon_{n})=\sum^m_{i=1} ( (A(x^{n})-y)_i^2+\epsilon_n^2)^{p/2},
\end{equation} 
the  $\epsilon$-perturbed $\ell_p$-norm residual in $\bar{x}$, $f_\epsilon(\bar{x} )$ would be the limit of $\mathcal{J}(x^{n},w^{n},\epsilon_{n})$ for $n\to \infty$. 
We denote a minimizer in dependence of $\epsilon$ with
\begin{equation}\label{xeps}
x^\epsilon \in \argmin\limits_{x} f_\epsilon(x),
\end{equation}
where we consider a global minimizer as we did for the minimization in the first step of Algorithm 1 as well.
Such minimizers are characterized by the following criterion. \\
\begin{lemma}\thlabel{fepsilon}
Let $\epsilon>0$ and define $w(z,\epsilon)=((A(z)-y)_i^2+\epsilon^2)^{(p-2)/2})_{i=1}^m$. If 
$$  \left\|A(z)-y\right\|^2_{\ell_2(w(z,\epsilon))} \leq  \left\|A(\tilde{z})-y\right\|^2_{\ell_2(w(z,\epsilon))} \text{  for all  }\tilde{z},$$
then it follows that $z=x^\epsilon\in \argmin\limits_{x} f_\epsilon(x).$
\end{lemma}
Proof: We want to prove that if 
$$\left\|A(z)-y\right\|^2_{\ell_2(w(z,\epsilon))} \leq  \left\|A(\tilde{z})-y\right\|^2_{\ell_2(w(z,\epsilon))} \text{  for all  }\tilde{z}.$$
then $f_\epsilon(z)\leq f_\epsilon(\tilde{z})$ for all  $\tilde{z}.$\\
We start with the given inequality
$$\left\|A(z)-y\right\|^2_{\ell_2(w(z,\epsilon))}= \sum_i \frac{(A_i(z)-y_i)^2}{[(A_i(z)-y_i)^2+\epsilon^2]^{(2-p)/2}} \leq  \sum_i \frac{(A_i(\tilde{z})-y_i)^2}{[(A_i(z)-y_i)^2+\epsilon^2]^{(2-p)/2}} =\left\|A(\tilde{z})-y\right\|^2_{\ell_2(w(z,\epsilon))}$$

As a first step we add $\epsilon^2$ to the numerator of each term of the sum and take the square root (monotone!) of the expressions. We obtain the inequality 

$$\left( \sum_i \frac{[(A_i(z)-y_i)^2+\epsilon^2]}{[(A_i(z)-y_i)^2+\epsilon^2]^{(2-p)/2}} \right)^{1/2}\leq  \left(\sum_i \frac{[(A_i(\tilde{z})-y_i)^2+\epsilon^2]}{[(A_i(z)-y_i)^2+\epsilon^2]^{(2-p)/2}}  \right)^{1/2}.$$
The left-hand-side can already be expressed as $f_\epsilon(z)$ and we use the $\frac{1}{2}$- triangle inequality for the square root. 
$$(f_\epsilon(z))^{1/2}\leq  \left(\sum_i \frac{[(A_i(\tilde{z})-y_i)^2+\epsilon^2]}{[(A_i(z)-y_i)^2+\epsilon^2]^{(2-p)/2}}  \right)^{1/2}\leq  \sum_i \frac{[(A_i(\tilde{z})-y_i)^2+\epsilon^2]^{1/2}}{[(A_i(z)-y_i)^2+\epsilon^2]^{(2-p)/4}} .$$
Now we apply H\"older's inequality and obtain
\begin{align*}
(f_\epsilon(z))^{1/2}&\leq\left(\sum_i((A_i(\tilde{z})-y_i)^2+\epsilon^2)^{p/2}\right)^{1/p}\cdot \left(\sum_i((A_i({z})-y_i)^2+\epsilon^2)^{(p-2)/4\cdot p/(p-1)}\right)^{2(p-1)/2p}\\
&=(f_\epsilon(\tilde{z}))^{1/p} \cdot \left[\left(\sum_i((A_i({z})-y_i)^2+\epsilon^2)^{(p-2)/4\cdot p/(p-1)}\right)^{2(p-1)/(p-2)}\right]^{(p-2)/2p}.
\end{align*}
We note now that $\frac{1}{(a+b)^\tau}\leq \frac{1}{a^\tau}+\frac{1}{b^\tau}$ for $a,b,\tau>0$. 
Hence as $2(p-1)/(p-2)$ is negative, using this factor on each summand gives the estimate
 $$(f_\epsilon(z))^{1/2}\leq(f_\epsilon(\tilde{z}))^{1/p} \cdot \left[\left(\sum_i((A_i({z})-y_i)^2+\epsilon^2)^{p/2}\right)\right]^{(p-2)/2p}=(f_\epsilon(\tilde{z}))^{1/p} \cdot(f_\epsilon(z))^{(p-2)/2p} .$$
 We rearrange and find 
$$(f_\epsilon(z))^{1/2-(p-2)/2p} =(f_\epsilon(z))^{1/p}\leq (f_\epsilon(\tilde{z}))^{1/p}.$$
As the p-th square root is monotone, we have 
$$f_\epsilon(z)\leq f_\epsilon(\tilde{z}).$$
Thus we have established the result.$\blacksquare$\\


\subsection{Convergence and Error Decay Rates for Algorithm 1}

Finally we can state convergence results for Algorithm 1 with the help of the tools and conditions already mentioned above, more precisely \thref{BCC} and \thref{USCC1}. Furthermore if $\|A(\cdot)-y\|^2_{\ell_p}$ is $\hat{C}$-strongly convex at $x^*$ certain requirements on the constant $\hat{C}$ in \eqref{subs1} will lead us to an error decay rate, which is linear in the case that $y \in \Ran(A)$, but one will have an additional error term scaling in the ineliminable discrepancy $\|A(x^*)-y\|_{\ell_p}$ otherwise.

\begin{theorem}\thlabel{conv}
Fix $y\in \mathbb{R}^m$and let $A:\mathbb{R}^k\rightarrow\mathbb{R}^m$ be a nonlinear map in $\mathcal{C}^0$ and the functionals $ \mathcal{J}(x,w^n,\epsilon_n)$  for $w^n,\epsilon_n$ as generated by Algorithm \ref{algo:1} for all $n\geq 0$, for which we consider the conditions
\begin{enumerate}
\item[(a)] the boundedness and coercivity condition (BCC), i.e., there exist $\alpha,\beta>0$ such that, for all $z\in\mathcal{B}(0,R^*)$:
\begin{equation*}
\alpha \|x^*-z\|_{\ell_2}\leq \|A(x^*)-A(z)\|_{\ell_p}\leq \beta \|x^*-z\|_{\ell_2};
\end{equation*}
\item[(b)] and the first uniform strong convexity condition (USCC-1), i.e., there exists a uniform constant $C>0$ such that for all $n\geq 0$  the following conditions holds
\begin{align}\label{subs1}
\mathcal{J}(x^n,w^n,\epsilon_n)- \mathcal{J}(x^{n+1},w^n,\epsilon_n) =\|A(x^n)-y\|^2_{\ell_2(w_n)}-\|A(x^{n+1})-y\|^2_{\ell_2(w_{n})} \ge {C} \|x^n-{x^{n+1}}\|_{\ell_2}^2.  
\end{align}
\end{enumerate}
Then the sequence $(x^n)_{n \in \mathbb{N}}$ generated by Algorithm 1 converges to a vector $\bar{x}$.
\begin{enumerate}
\item[(i)]if $\epsilon=\lim\limits_{n\rightarrow\infty} \epsilon_n=0$, and condition (a) holds, then $\bar{x}=x^*$ is the solution to the $\ell_p$-minimization problem \eqref{ellp}. Moreover  $y\in Ran(A)$ and $y=A(x^*)$.
\item[(ii)]if $\epsilon=\lim\limits_{n\rightarrow\infty} \epsilon_n>0$, and both conditions (a) and (b) hold, then $\bar{x}=x^\epsilon$ as defined in \eqref{xeps} and $x^\epsilon \in \mathcal{B}(0,R^*)$. For simplicity here we assume that $x^\epsilon$ is actually the unique global minimizer of $f_\epsilon$. 
\end{enumerate}
\begin{enumerate}
\item[(c)] Let the error at the $n$-th step be denoted as $E^n$ and the unavoidable error  $E^*=\|A(x^*)-y\|_{\ell_p}^2$. If condition (a) is fulfilled as well as the the second uniform convexity condition(USCC-2) i.e there exists a uniform constant $\hat{C}>0$ such that for all $n\geq 0$, the following conditions hold
$$\|A(x^n)-y\|^2_{\ell_p} - \|A(x^*)-y\|^2_{\ell_p} \ge \hat{C} \|x^n-{x^*}\|_{\ell_2}^2  $$ 
for all $n\geq 0$, where $\hat{C}>0$ is such that  $\mu:=\frac{2^{1+2/p}(m^2+1)\beta^2}{\hat{C}}<1$ and $\nu=\frac{2^{1+2/p}(m^2+1-2^{-2/p})}{\hat{C}} $, we can additionally infer the following property:
\end{enumerate}
\begin{enumerate}
\item[(iii)] the error decay rate can be characterized in terms of the just defined errors $E^n$ and $E^*$ as follows:
\begin{equation}\label{E1}E^{n+1}\leq\mu E^n +\nu E^*\end{equation}
or 
\begin{equation}\label{E2}E^{n+1}\leq\mu^n E^0+\sum\limits_{r=1}^n \mu^r \nu E^*.\end{equation}
Taking the limits for $n\rightarrow \infty$ gives an asymptotic error of the order of $E^*$
\begin{equation}\label{E3}\bar{E} :=  \|A(\bar{x})-y\|_{\ell_p}^2 \leq \frac{\nu}{1-\mu} E^*.\end{equation}
\end{enumerate}
\end{theorem}
Proof:
(i) In this case we want to prove that $x^n$ converges, and that its limit is the solution to the $\ell_p$-minimization problem \eqref{ellp}. Suppose that $\epsilon_{n_0}=0$ for some $n_0$. Then by the definition of the algorithm, we know that the iteration is stopped at $n=n_0$, and $x^n=x^{n_0}, n\geq n_0$. Therefore $\bar{x}=x^{n_0}$. \\
From the definition of $\epsilon_n$, it then also follows that $ \max_i((A_i(x^{n+1})-y_i))^2=0$, hence $\left\| A(\bar x)-y\right\|^p_{\ell_p}=0$ and, in view of (a) we have $\bar x= x^*$.
\vspace{2mm}
\\Suppose now that $\epsilon_n>0$ for all $n$. Since $\epsilon_n\rightarrow 0$, there is an increasing sequence of indices $(n_l)$ such that $\epsilon_{n_l}< \epsilon_{n_l-1}$ for all $l$. \\
Note that $x^n$ is a bounded sequence according to \thref{xbound} and there exists a subsequence $(t_s)_{s \in \mathbb{N}}$ of $(n_l)_{l \in \mathbb{N}}$ such that $(x^{t_s})_{s \in \mathbb{N}}$ converges to a limit point $\tilde{x}$. 
We observe that by the definition of $\epsilon_{t_s}$
$$\sum_i((A(x^{t_s})-y_i)_i^2+\epsilon_{t_s}^2)^{p/2}< \sum_i 2^{p/2}\max_j|A(x^{_{t_{s}}})_j-y_j|^p.$$
In fact, by the definition of $\left( \epsilon_n\right)_{n \in \mathbb{N}}$ we know that if $\epsilon_{t_s}$ falls below $\tilde{\epsilon}$, then $\epsilon_{t_s}=\max_j|A(x^{_{t_{s}}})_j-y_j| <\tilde{\epsilon}$ and hence $\epsilon_{t_s} \rightarrow 0$ implies that $\max_j|A(x^{t_{s}})_j-y_j| \rightarrow 0$.
Using the continuity of $A$ it follows that 
$$0 \leq \sum_i|A(\tilde{x})_i-y_i|^p=\lim_{s\rightarrow \infty}\sum_i((A(x^{t_s})_i-y_i)^2+\epsilon_{t_{s}}^2)^{p/2} \leq \lim\limits_{s\rightarrow \infty}2^{p/2}m\epsilon_{t_s}^p=0.$$
We must still show that $x^n\rightarrow x^*$. Since $x^{t_s}\rightarrow x^*$ and $\epsilon_{t_s}\rightarrow0$, it holds that $\mathcal{J}(x^{t_s}, w^{t_s},\epsilon_{t_s}) \rightarrow 0=\sum_i|A(x^*)_i-y_i|^p.$
 By the monotonicity property of $\mathcal J$, we get $\mathcal{J}(x^{n}, w^{n},\epsilon_{n}) \rightarrow 0 =\sum_i|A(x^*)_i-y_i|^p$.
 \vspace{2mm}
 \\Since \eqref{JJ} implies $\mathcal{J}(x^{n}, w^{n},\epsilon_{n})- m \epsilon_n^p\leq \sum_i|A_i(x^n)-y_i|^p\leq\mathcal{J}(x^{n}, w^{n},\epsilon_{n}),$
we obtain $$\lim_{n\to\infty} \sum_i|A_i(x^n)-y_i|^p =\sum_i|A_i(x^*)-y_i|^p =0.$$
By using the BCC we obtain 
\begin{align*}
0&\leq \lim \sup_{n \rightarrow \infty}\left\| x^n-x^*\right\|_2< \lim\sup_{n \rightarrow \infty} \left ( \frac{1}{\alpha} \left(\sum_i|A_i(x^n)-y_i|^p\right)^{1/p}+\frac{1}{\alpha}\left( \sum_i|A_i(x^*)-y_i|^p\right)^{1/p} \right)\\
&= \frac{2}{\alpha}\lim\limits_{n \rightarrow \infty}  \left(\sum_i|A_i(x^n)-y_i|^p\right)^{1/p}=0.
\end{align*}
This statement completes the proof that $x^n \rightarrow x^*$ in this case.\\

(ii) We shall first show that $x^n\rightarrow x^\epsilon$, $n\rightarrow \infty$ with $x^\epsilon \in \argmin\limits_{x} f_\epsilon(x)$. We already observed that $(x^n)_{n\in \mathbb{N}_0}$ is a bounded sequence in $\mathcal{B}(0,R^*)$ and hence this sequence has accumulation points. Let $(x^{n_l})_{l\in \mathbb{N}_0}$  be any convergent subsequence of $(x^n)_{n\in \mathbb{N}_0}$  and let $\bar{x}$ be its limit. We want to show that $\bar{x}=x^\epsilon$. \\
By continuity of $A$ it follows that $\lim\limits_{l\rightarrow \infty} w^{n_l}_i= [(A(\bar{x})_i-y_i)^2+\epsilon^2]^{(p-2)/2}={w}(\bar{x},\epsilon)_i:=\bar{w}_i, i=1,\cdots, m$. \\
On the other hand, by invoking \thref{cauchyx} , we obtain also that $x^{n_l+1}\rightarrow \bar{x}, i\rightarrow\infty$. \\
From the minimality properties we know that for all $n_l$
\begin{equation}\label{comp}
\left\|A(x^{n_l+1})-y \right\|_{\ell_2(w^{n_l})}\leq\left\|A({z})-y \right\|_{\ell_2(w^{n_l})}, \text{for all }z \in \mathbb{R}^k.
\end{equation}
By passing to the limit for $n_l\rightarrow \infty$ for $z$ fixed, we obtain also
$$\left\|A(\bar{x})-y \right\|_{\ell_2(\bar{w})}\leq\left\|A({{z}})-y \right\|_{\ell_2(\bar{w})}.$$
Now \thref{fepsilon} implies that $\bar{x}=x^\epsilon$. As we assumed that $x^\epsilon$ is the unique minimizer of $f_\epsilon$, it is  the unique accumulation point of $(x_n)_{n\in \mathbb{N}}$ and is also its limit. This establishes (ii).\\

(iii) We would like to get an error bound at the $(n+1)$-th iteration from the chain of estimates as developed in the following. We start with the error 
\begin{equation}\label{potato}
\|x^n-x^*\|_{\ell_2}^2 \geq \frac{1}{\beta^2} \|A(x^n)-A(x^*)\|_{\ell_p}^2 \geq \frac{1}{\beta^2} \left(\frac{1}{2} \|A(x^n)-y\|_{\ell_p}^2 -\|A(x^*)-y\|_{\ell_p}^2\right),
\end{equation}
where we used the BCC. 

We now need to pass to our functional $\mathcal{J}$, because we only have monotonicity of it along the iterations thanks to \thref{monotonicity}. Let $\|\epsilon_n\|_{\ell_2(w_n)}:=\| \epsilon_n\cdot (1,\dots,1)^T\|_{\ell_2(w_n)}$. We observe that
\begin{align*}
\|A(x^n)-y\|_{\ell_p}^2&=\left(\sum\limits^m_{i=1} |A_i(x^n)-y_i|^p\right)^{\frac{2}{p}}\geq\left(\sum\limits^m_{i=1}\frac{(A_i(x^n)-y_i)^2+\epsilon_n^2 -\epsilon_n^2}{((A_i(x^n)-y_i)^2+\epsilon_n^2)^{(2-p)/2}}\right)^{\frac{2}{p}}\\
&\geq2^{1-2/p} \mathcal{J}(x^n,w^n,\epsilon_n)^{\frac{2}{p}}-\|\epsilon_n\|_{\ell_2(w_n)}^{\frac{4}{p}} \geq 2^{1-2/p}\mathcal{J}(x^{n+1},w^{n+1},\epsilon_{n+1})^{\frac{2}{p}}-\|\epsilon_n\|_{\ell_2(w_n)}^{\frac{4}{p}}\\
&\geq 2^{1-2/p} \|A(x^{n+1})-y\|_{\ell_p}^2 -\|\epsilon_n\|_{\ell_2(w_n)}^{\frac{4}{p}}.
\end{align*}
From \eqref{potato} and the latter estimate we obtain
\begin{align*}
\|x^n-x^*\|_{\ell_2}^2&\geq \frac{1}{2\beta^2}\left[ 2^{1-2/p}\|A(x^{n+1})-y\|_{\ell_p}^2 -\|\epsilon_n\|_{\ell_2(w_n)}^{\frac{4}{p}}-2\|A(x^*)-y\|_{\ell_p}^2\right].
\end{align*}
We add and subtract the term $\frac{1-2^{-2/p}}{\beta^2}\|A(x^*)-y\|_{\ell_p}^2$ and rearrange 
$$ \|x^n-x^*\|_{\ell_2}^2+ \frac{1-2^{-2/p}}{\beta^2} \|A(x^*)-y\|_{\ell_p}^2+\frac{1}{2\beta^2}\|\epsilon_n\|_{\ell_2(w_n)}^{\frac{4}{p}}\geq \frac{1}{2^{2/p}\beta^2} \left(\|A(x^{n+1})-y\|_{\ell_p}^2 -\|A(x^*)-y\|_{\ell_p}^2\right)$$
We rearrange again and use \thref{USCC2} to obtain 
\begin{align*}
&\|x^n-x^*\|_{\ell_2}^2+\frac{1-2^{-2/p}}{\beta^2} \|A(x^*)-y\|_{\ell_p}^2+\frac{1}{2\beta^2}\|\epsilon_n\|_{\ell_2(w_n)}^{\frac{4}{p}}\\
&\geq \frac{\hat{C}}{2^{1+2/p}\beta^2} \|x^{n+1}-x^*\|_{\ell_2}^2.
\end{align*}
We examine the expression $\|\epsilon_n\|_{\ell_2(w_n)}^{\frac{4}{p}}$ a bit closer and estimate it from above 

\begin{align*}
\|\epsilon_n\|_{\ell_2(w_n)}^{\frac{4}{p}}&=\left(\sum\limits_{i=1}^m\frac{\epsilon^2_n}{[(A_i(x^n)-y_i)^2+\epsilon_n^2]^{(2-p)/2}}\right)^{\frac{2}{p}}\leq \left( \sum\limits_{i=1}^m \epsilon_n^p\right)^{\frac{2}{p}}\\
&\leq m^2  \|A(x^n)-y\|_{\ell_p}^2\leq 2m^2\|A(x^n)-A(x^*)\|_{\ell_p}^2 +2m^2\|A(x^*)-y\|_{\ell_p}^2 \\
&\leq2{\beta^2}{m^2} \|x^n-x^*\|_{\ell_2}^2 +2{m^2} \|A(x^*)-y\|_{\ell_p}^2,
\end{align*}
where we used the definition of $\epsilon_n$ and the fact that the maximum absolute value of a single entry of a vector is smaller than the norm of the vector. 
Hence we combine the results above to obtain \eqref{E1}
\begin{align*}
E^n&=\|x^{n+1}-x^*\|_{\ell_2}^2\leq \frac{2^{1+2/p}(m^2+1)\beta^2}{\hat{C}}\|x^n-x^*\|_{\ell_2}^2+\frac{2^{1+2/p}(m^2+1-2^{-2/p})}{\hat{C}}  \|A(x^*)-y\|_{\ell_p}^2=\mu E^{n+1}+\nu E^*.
\end{align*}

By recurrently substituting $E^n$ by its predecessors we obtain  \eqref{E2}
\begin{align*}
E^{n+1}\leq\mu^n E^0+\sum\limits_{r=1}^n \mu^r \nu E^*.
\end{align*}
Taking the limit $n\rightarrow \infty$ gives \eqref{E3}.
$\blacksquare$\\

\begin{remark}
\begin{enumerate}
\item{Notice that for $\bar{x}=x^*$ \eqref{E3} is trivial while for $\bar{x}=x^\epsilon$ \eqref{E3} gives information on $x^\epsilon$ as a quasi-minimizer.} 
\item{Requiring the inequality in \thref{fepsilon} to be valid for all $\tilde{z}$ and not only for $z \in \mathcal{B}(0,2R^*)$ is due to the global minimization of $\mathcal{J}(x,w^{n_l},\epsilon_{n_l})$ w.r.t. $x$. This minimization gives us $x^{n_l+1}$, that is the global minimizer that is compared to all other $z$ in \eqref{comp} in step (ii). }
\item{The values of $\mu$ and $\nu$ represent the worst upper bounds up to the point where $\epsilon_n=\mathcal{M}^n$ in Algorithm 1.  In cases where $\epsilon_n=\mathcal{N}^n$, we could choose $\tilde{\mu}= \frac{2^{2+2/p}\beta^2}{\hat{C}}, \tilde{\nu}=\frac{2^{2+2/p}-2}{\hat{C}}$ instead of $\mu,\nu$ for these steps, hence better constants.}
\end{enumerate}
\end{remark}

\section{Local Convexification of the Auxiliary Functional}

In the case where we can not fulfill the uniform strong convexity condition of \thref{USCC1} one is in the unpleasant situation of facing a not even locally convex optimization problem. In this case convergence of the NR-IRLS as defined above can not be guaranteed as we did. For this harder case we suggest a strategy of adaptive adjustment of the algorithm towards local convexification around the current iterate, which can ensure again convergence of the generalized IRLS to at least a critical point. 

\subsection{Starting points for the convexified algorithm}
In the algorithm as stated above, for the first iteration of the generalized IRLS, where $w_0=(1,\dots,1)^T$, just a usual nonlinear $\ell_2$-least squares step is executed to obtain the minimizer of the functional $\mathcal{J}(x,w^0,\epsilon_0)$. Already at this point in the locally nonconvex case one might possibly encounter several local minimizers and which one is set to be the next iterate $x^1$ depends on the particular starting point $x_0$ of the optimization strategy used in as already pointed out in \thref{x0}. 

As the convexification that we have in mind is locally around the current iterate it is important to think about the influence of the first iterate and therefore also the choice of $x^0$ for the optimization process in the first step. The overall result of the algorithm might vary strongly for different choices!

At this point one can take a little excursion to complex analysis to infer that if $A(x)$ is of analytic meaning that $A(x)=\left( \pi_1(x),\dots,\pi_i(x),\dots,\pi_m(x)\right)$, where the $\pi_i(x)$ are analytic, then we only have isolated zeros of $\nabla \|A(x)-y\|_{\ell_p}^p$ for $p>1$ and only a finite number of them on any compact set.\\
Our suggestion is now to invest the computational effort of finding several or even all critical points occurring in the first iteration. Our hope motivating this is that the location of critical points of the functional $\|A(x)-y\|_{\ell_p}^p$ does not change strongly with varying $p$ and that a global minimizer for $p$ might be not too far from a local $\ell_2$-minimizer (or critical point).
For finding these local minimizers one can use the methods described in the beginning of the paper for instance a Levenberg-Marquardt type algorithm with different (maybe randomly chosen) starting points. 
The stationary points found will be denoted as $x^{*_1}_{\ell_2}, x^{*_2}_{\ell_2},\dots,x^{*_L}_{\ell_2}$. For each of them we would suggest to use them as starting points for an adjusted version of NR-IRLS that we will introduce in the following. 
After having executed this adjusted version of the generalized IRLS for all $x^{*_1}_2, x^{*_2}_2,\dots,x^{*_L}_2$ and having obtained $L$ possible solutions $x^{*_1}, x^{*_2},\dots,x^{*_L}$ for the $\ell_p$-minimization problem, we chose the $x^{*_s}$ giving the lowest value of $\|A(x^{*_s})-y\|_{\ell_p}$ as our best candidate for the desired $\ell_p$-minimizer.\\

\subsection{Convexification of the Auxiliary Functional and Convexified Algorithm}

To overcome the disadvantage of not fulfilling the condition in \thref{USCC1}, we want to introduce another way to ensure uniform local convexity and want to gain convergence by a convexifying adaption of Algorithm 1.\\
Let us now consider the case where $w,\epsilon$ are fixed. 
Our idea is now to construct a convexification of our original functional $\mathcal{J}$ as follows 
\begin{equation}\label{Jconv}
 \mathcal{J}_{\omega,u}(x,w,\epsilon)=\mathcal{J}(x,w,\epsilon)+\omega \|x-u\|^2_{\ell_2},
 \end{equation}
for a parameter $\omega>0$ and $u \in \mathbb{R}^k$.  This kind of regularization of a functional is also called Moreau envelope in nonsmooth convex optimization \cite{Rockafellar09}.

We want to embed this type of straightforward convexification \eqref{Jconv} into the first step of the Nonlinear Residual IRLS algorithm and perform a minimization of a regularized problem instead of a minimization as performed in the original formulation of NR-IRLS to obtain a corresponding sequence of minimizers $x^n$.\\
To incorporate \eqref{Jconv} into the original formulation of the algorithm we have to  make an appropriate choice for the newly introduced parameters $u$ and $\omega$.
We decide to fix $\omega>0$ generously large enough and constant for all iterations and set $u=x^n$ at the $n$-th step.  Thereby we obtain our next iterate as follows
\begin{equation}\label{inner}
x^{n+1} = \argmin\limits_{x} \, \mathcal{J}_{\omega, x^{n}}(x,w^n,\epsilon_n), \\
\end{equation}

\begin{remark}

\begin{itemize}
\item[(a)]  For the solution of the convex minimization problem appearing here a great variety of well-understood convex optimization methods is available see Section 2.2. 
\item[(b)]  Explanations why this choice for the adaption of the algorithm will resolve our lack of local convexity and recommendations for the concrete choice of $\omega$ will follow later in the theoretical analysis section.
\end{itemize}
\end{remark}

Now we can turn towards the formulation of an adapted NR-IRLS algorithm:\\

For initialization we set $x^1= x^{*_l}_{\ell_2}$ and define as usual 
$$\mathcal{N}^1=\min_i(|A_i(x^1)-y_i)|) \text{  and  } \mathcal{M}^1=\max_i(|A_i(x^1)-y_i|),$$
to obtain
$$\epsilon_1=\min\left( \max( \mathcal{N}^1,\tilde{\epsilon}),\epsilon_0, \mathcal{M}^1\right)$$
where $\tilde{\epsilon}$ is a small number and
\begin{align*}
w^1&=\left[( A(x^1)-y).^2+\epsilon^2_1\right]^{(p-2)/2}.
\end{align*}
Now we will adapt the functional $\mathcal{J}$ as in \eqref{inner} 
\begin{equation}\label{adJ}
\mathcal{J}_{\omega_1,x^1}(x,w^1,\epsilon_1)=\frac{p}{2}\left[\sum\limits^m_{i=1}\left[w^1_i( A_i(x)-y_i)^2\right]+ \sum_{i=1}^m(\epsilon_1^2w^1_i+\frac{2-p}{p}(w^1_i)^{p/(p-2)})\right]+\omega\|x-x^1\|^2_{\ell_2}.
\end{equation}

We obtain the next candidate for a minimizer  by performing the step
\begin{equation}
x^{2} = \argmin\limits_{x} \, \mathcal{J}_{\omega, x^{1}}(x,w^1,\epsilon_1).\\
\end{equation}


Continuing in this fashion results in the following adapted algorithm 

\begin{algorithm}[h]
\KwSty{Algorithm 2: Convexified nonlinear residual iteratively reweighted least squares }:\\

\KwIn{$A:\mathbb{R}^k\rightarrow \mathbb{R}^{m}$, $y\in\mathbb{R}^m$ }
We initialize by taking $\epsilon_1= \min\limits_i ( A_i(x^{*_l}_2)-y_i)^2$ and $w^1=\left[( A(x^{*_l}_2)-y).^2+\epsilon^2_1\right]^{(p-2)/2}$. \\
We fix $x^1=x^{*_l}_{\ell_2}$ for the first step. We then recursively define for $n= 1,2,3,\dots$ 
\begin{align*}\label{JJ} 
\mathcal{J}_{\omega_n,x^{n}}(x,w^n,\epsilon_n)=\frac{p}{2}\left[\sum\limits^m_{i=1}\left[w^n_i( A_i(x)-y_i)^2\right]+ \sum_{i=1}^m(\epsilon_n^2w^n_i+\frac{2-p}{p}(w^n_i)^{p/(p-2)})\right]+\omega\|x-x^{n}\|^2_{\ell_2}.
\end{align*}
 We obtain the next iterate by performing a step 
\begin{equation}
x^{n+1} = \argmin\limits_{x} \, \mathcal{J}_{\omega,  x^{n}}(x,w^n,\epsilon_n).
\end{equation}
 We redefine 
 $$\mathcal{N}^{n+1}=\min_i(|A_i(x^{n+1})-y_i)|) \text{  and  } \mathcal{M}^{n+1}=\max_i(|A_i(x^{n+1})-y_i|),$$
giving us
$$\epsilon_{n+1}=\min\left( \max( \mathcal{N}^{n+1},\tilde{\epsilon}),\epsilon_{n}, \mathcal{M}^{n+1}\right),$$
where $\tilde{\epsilon}$ is a small constant and
\begin{align*}
w^{n+1}&=\left[( A(x^{n+1})-y).^2+\epsilon^2_{n+1}\right]^{(p-2)/2}.
\end{align*}
%
\KwOut{$x^{2}$, $x^{3},\ldots$}
\end{algorithm}

\subsection{Convergence Analysis for the Convexified Algorithm}

In the following section we examine Algorithm 2 and its behavior mostly analogous to the original version of NR-IRLS. We again discuss some preliminary properties and facts already familiar from Section 3 for the adapted algorithm now that will be useful to show its convergence later on. Having obtained these results we shall develop the proof of convergence  of Algorithm 2 in several lemma.\\

A first important note is that we also have a type of monotonicity property for the adjusted functionals: 

\begin{lemma}\thlabel{mono2} The inequalities 
\begin{align}
\mathcal{J}(x^{n},w^n,\epsilon_n)& =\mathcal{J}_{\omega,x^n}(x^n,w^n,\epsilon_n)\geq \mathcal{J}_{\omega,x^{n}}(x^{n+1},w^n,\epsilon_n)\\
&\geq \mathcal{J}_{\omega,x^{n}}(x^{n+1},w^n,\epsilon_{n+1})\geq \mathcal{J}_{\omega,x^{n}}(x^{n+1},w^{n+1},\epsilon_{n+1})\\
&\geq \mathcal{J}_{\omega,x^{n+1}}(x^{n+1},w^{n+1},\epsilon_{n+1})= \mathcal{J}(x^{n+1},w^{n+1},\epsilon_{n+1})
\end{align}
hold for all $n\geq 0$.
\end{lemma}
Proof: 
Here the first and inequality follows from the minimization property defined via the first step of Algorithm 2,  the second
inequality from $\epsilon_{n+1}\leq \epsilon_n$, and the second last inequality from the minimization property that defines
$w_{n+1}$. The last inequality results from the fact that the norm of a difference of different vectors is greater than 0.$\blacksquare$\\

From the monotonicity property of the sequence $(\mathcal{J}(x^n,w^n,\epsilon_n))_{n\in \mathbb{N}}$ and its boundedness from below we know that this sequence is convergent.

Moreover note that also in the convexified case the sequence  ${(x}^{n})_{n \in \mathbb{N}}$ resulting as the output of Algorithm 2 is bounded, as it can be shown analogously to \thref{xbound} and hence ${(x}^{n})_{n \in \mathbb{N}}\in \mathcal{B}(0,R^*).$\\


At this point we want to explain further why we chose the formulation in \eqref{inner} to adapt the functional $\mathcal{J}$. The introduction of the regularization term is the key to regain the USCC-1 for the adapted functional \eqref{Jconv} as we can guarantee the existence of a positive USCC-1 constant $\tilde{C}$ with an appropriate choice of $\omega$. 

\begin{lemma}\thlabel{USCC12}
Let $A:\mathbb{R}^k\rightarrow\mathbb{R}^m$ be a { nonlinear map in ${\mathcal{C}^{0}}$} and $\mathcal{J}(x,w,\epsilon)$ as defined in \thref{J} and $\mathcal{J}_{\omega,u }(x,w,\epsilon)$ as defined in \eqref{Jconv}. 
{ We additionally assume that
\begin{align} \label{Lipcond}
& \left |  t   [\|A(tx^n+(1-t)x^{n+1})-y \|_{\ell_2(w^n)}^2  - \|A(x^n)-y \|_{\ell_2(w^n)}^2 ] \right .  \\  \nonumber 
& \left .  +(1-t)[  \|A(tx^n+(1-t)x^{n+1})-y \|_{\ell_2(w^n)}^2 - \|A(x^{n+1})-y \|_{\ell_2(w^n)}^2]  \right | \\ \nonumber
&\leq L t (t-1) \|x^{n} -x^{n+1}\|_{\ell_2}^2.
\end{align}
for some $L>0$ indepedent of $n \in \mathbb N$ and for all $t \in [0,1]$.}
Let $(x^n)_{n\in \mathbb{N}}$ be the sequence of minimizers output by Algorithm 2. Then for $\omega>0$ large enough the USCC-1 is fulfilled for the adapted functional  in \eqref{Jconv}, i.e., there exists a uniform constant $\tilde{C}>0$ such that for all $n\ge0$ holds
 \begin{align*}
\mathcal{J}_{\omega,x^n}(x^n,w^n,\epsilon_n)- \mathcal{J}_{\omega,x^n}(x^{n+1},w^n,\epsilon_n) \ge {\tilde{C}} \|x^{n+1}-x^n\|_{\ell_2}^2 
\end{align*}
\end{lemma}
\begin{remark}\label{remrem}
{ Before starting the proof of this lemma, let us discuss for a moment the validity of \eqref{Lipcond}. Were $A \in \mathcal C^2$ and $\epsilon^n \geq \epsilon$ for all $n \in \mathbb N$, then the Hessian 
$$
\nabla^2 F_{w^n}(x)= \sum_{i=1}^m w_i^n \left[ \nabla A(x)_i \nabla A(x)_i^T+(A(x)_i-y_i)\nabla^2 A(x)_i\right],
$$
of the map $$x \to F_{w^{n}}(x)=\|A(x)-y \|_{\ell_2(w^n)}^2$$ would be actually uniformly bounded on $B(0,R^*)$, say by a constant $L'>0$. By considering two Taylor expansions around the point $x=  tx^n+(1-t)x^{n+1}$ in the following expressions,
we actually obtain easily a uniform estimate of the type \eqref{Lipcond}:
\begin{eqnarray*}
&& \left |  t   [\|A(tx^n+(1-t)x^{n+1})-y \|_{\ell_2(w^n)}^2  - \|A(x^n)-y \|_{\ell_2(w^n)}^2 ] \right .  \\ 
&& \left . +(1-t)[  \|A(tx^n+(1-t)x^{n+1})-y \|_{\ell_2(w^n)}^2 - \|A(x^{n+1})-y \|_{\ell_2(w^n)}^2]  \right | \\
&=& \left |  -t \nabla F_{w^n}(tx^n+(1-t)x^{n+1})^T (x^n - tx^n+(1-t)x^{n+1}) \right . \\ 
&& \phantom{XXXXX} +  -t (x^n - tx^n+(1-t)x^{n+1})^T \nabla^2 F_{w^n}(\xi_{t}^n)  (x^n - tx^n+(1-t)x^{n+1})   \\ 
&& -(1-t)  \nabla F_{w^n}(tx^n+(1-t)x^{n+1})^T(x^{n+1} - tx^n+(1-t)x^{n+1}) +\\
&& \left . \phantom{XXXXX} + -(1-t)(x^{n+1} - tx^n+(1-t)x^{n+1})^T \nabla^2 F_{w^n}(\eta_{t}^n)  (x^{n+1} - tx^n+(1-t)x^{n+1}) \right | \\
\end{eqnarray*}
Now, we have that
$$
 -t \nabla F_{w^n}(tx^n+(1-t)x^{n+1})^T (x^n - tx^n+(1-t)x^{n+1}) = -t (1-t)\nabla F_{w^n}(tx^n+(1-t)x^{n+1})^T (x^n -x^{n-1}) 
$$
as well as 
$$
-(1-t)  \nabla F_{w^n}(tx^n+(1-t)x^{n+1})^T(x^{n+1} - tx^n+(1-t)x^{n+1}) = t (1-t)\nabla F_{w^n}(tx^n+(1-t)x^{n+1})^T (x^n -x^{n-1}).
$$
Hence the first order terms in the sum above are one the opposite of the other and they delete each other. One is therefore left only with the second order terms. By using the boundedness of the Hessians 
and by observing that
$$
\|(x^n - tx^n+(1-t)x^{n+1}\|_{\ell_2}^2 = (1-t)^2 \|x^n - x^{n+1}\|_{\ell_2}^2, \mbox{ and } \|(x^{n+1} - tx^n+(1-t)x^{n+1}\|_{\ell_2}^2 = t^2 \|x^n - x^{n+1}\|_{\ell_2}^2
$$
one obtains
\begin{eqnarray*}
&& \left |  t   [\|A(tx^n+(1-t)x^{n+1})-y \|_{\ell_2(w^n)}^2  - \|A(x^n)-y \|_{\ell_2(w^n)}^2 ] \right .  \\ 
&& \left . +(1-t)[  \|A(tx^n+(1-t)x^{n+1})-y \|_{\ell_2(w^n)}^2 - \|A(x^{n+1})-y \|_{\ell_2(w^n)}^2]  \right | \\
&\leq& L' t(1-t)^2 \|x^n - x^{n+1}\|_{\ell_2}^2 + L' t^2(1-t)   \|x^n - x^{n+1}\|_{\ell_2}^2 \leq L t(t-1)  \|x^n - x^{n+1}\|_{\ell_2}^2,
\end{eqnarray*}
where in the last inequality we used that $t \in [0,1]$ and $L = 2 L'$. We claim therefore that  \eqref{Lipcond} is a reasonable assumption also in case 
where $A$ is not as smooth. But  it is here crucial that $\epsilon^n \geq \epsilon$ for all $n \in \mathbb N$ as it is actually used later in the
proof of Theorem \ref{mainth}.

}
\end{remark}

Proof: 
{ In view of the condition \eqref{Lipcond}, 
we obtain the following estimates for $t \in [0,1]$
\begin{eqnarray*}
&& \left | \mathcal{J}(tx^n+(1-t)x^{n+1},w^n,\epsilon_n)  - [ t\mathcal{J}(x^n,w^n,\epsilon_n) +(1-t)\mathcal{J}(x^{n+1},w^n,\epsilon_n)] \right |\\
&\leq & \left |  t   [\|A(tx^n+(1-t)x^{n+1})-y \|_{\ell_2(w^n)}^2  - \|A(x^n)-y \|_{\ell_2(w^n)}^2 ] \right .  \\ 
&& \left . +(1-t)[  \|A(tx^n+(1-t)x^{n+1})-y \|_{\ell_2(w^n)}^2 - \|A(x^{n+1})-y \|_{\ell_2(w^n)}^2]  \right | \leq L t (t-1) \|x^{n} -x^{n+1}\|_{\ell_2}^2.
\end{eqnarray*}

Consequently we can write
\begin{align*}
\mathcal{J}(tx^n+(1-t)x^{n+1},w^n,\epsilon_n) \le t\mathcal{J}(x^n,w^n,\epsilon_n) +(1-t)\mathcal{J}(x^{n+1},w^n,\epsilon_n)  -  {C} t(1-t) \|x^n-x^{n-1}\|_{\ell_2}^2,
\end{align*}
where $C=-L$ is a uniform constant that is not necessarily positive as we do not assume yet that strong convexity holds for $ \mathcal{J}(\cdot,w^n,\epsilon_n) $ in this case, but certainly $C>-\infty$.} \\
We then add $\omega\|tx^n+(1-t)x^{n+1}-x^n\|^2_{\ell_2}$ on both sides of the inequality
\begin{align*}
&\mathcal{J}(tx^n+(1-t)x^{n+1},w^n,\epsilon_n) +\omega\|tx^n+(1-t)x^{n+1}-x^n\|^2_{\ell_2} \\
&\le t\mathcal{J}(x^n,w^n,\epsilon_n) +(1-t)\mathcal{J}(x^{n+1},w^n,\epsilon_n)  -  {C} t(1-t) \|x^n-x^{n-1}\|_{\ell_2}^2+ \omega\|tx^n+(1-t)x^{n+1}-x^n\|^2_{\ell_2}
\end{align*}
and reformulate
\begin{align*}
&\mathcal{J}_{\omega,x^n }(tx^n+(1-t)x^{n+1},w^n,\epsilon_n) \\
&\le t\mathcal{J}_{\omega,x^n }(x^n,w^n,\epsilon_n) +(1-t)\mathcal{J}(x^{n+1},w^n,\epsilon_n)  -  {C} t(1-t) \|x^n-x^{n-1}\|_{\ell_2}^2+(1-t)^2 \omega\|x^n-x^{n+1}\|^2_{\ell_2}.
\end{align*}
We add and subtract the term $(1-t)\omega\|x^n-x^{n+1}\|^2_{\ell_2}$ and gain as a result
\begin{align*}
&\mathcal{J}_{\omega,x^n }(tx^n+(1-t)x^{n+1},w^n,\epsilon_n) \\
&\le t\mathcal{J}_{\omega,x^n }(x^n,w^n,\epsilon_n) +(1-t)\mathcal{J}_{\omega,x^n }(x^{n+1},w^n,\epsilon_n)  -  (C+\omega) t(1-t) \|x^n-x^{n-1}\|_{\ell_2}^2.
\end{align*}
The last inequality { establishes actually} the strong convexity condition for $\mathcal{J}_{\omega,x^n }(\cdot,w^n,\epsilon_n)$ at $x^{n+1}$ in the point $x^n$. If we carry out calculations analogous to the {ones in the } proof of \thref{fc} we obtain 
\begin{align*}
\mathcal{J}_{\omega,x^n}(x^n,w^n,\epsilon_n)- \mathcal{J}_{\omega,x^n}(x^{n+1},w^n,\epsilon_n) \ge {\tilde{C}} \|x^{n+1}-x^n\|_{\ell_2}^2 ,
\end{align*}
where $\tilde{C}=C+\omega$ which is positive if $\omega$ { is } chosen large enough.
$\blacksquare$


Now we would like to prove again that the iterates are getting arbitrarily close :
\begin{lemma}\thlabel{chiave}
Let $A:\mathbb{R}^k\rightarrow \mathbb{R}^m$ be a nonlinear map with $A \in \mathcal{C}^0$ {  and $(x^n)_{n \in \mathbb{N}}$
 and $(w^n)_{n \in \mathbb{N}}$ be the sequences generated by Algorithm 2, so that condition \eqref{Lipcond} holds.}  
Then, for $\omega>0$ large enough
$$\left\|x^{n}-x^{n+1}\right\|^2_{\ell_2} \to 0  \text{  as  }  n \to \infty.$$
\end{lemma}

Proof:
With the monotonicity property above we obtain 
%

$$
\|\mathcal{J}(x^n,w^n,\epsilon_n)-\mathcal{J}(x^{n+1},w^{n+1},\epsilon_{n+1})\|^2_{\ell_2}\geq\| \mathcal{J}_{\omega,x^n}(x^n,w^n,\epsilon_n)-\mathcal{J}_{\omega,x^n}(x^{n+1},w^{n},\epsilon_{n}) \|^2_{\ell_2}.
$$
By \thref{USCC12} we get
\begin{align*}\label{29381}
\| \mathcal{J}_{\omega,x^n}(x^n,w^n,\epsilon_n)-\mathcal{J}_{\omega,x^n}(x^{n+1},w^{n},\epsilon_{n}) \|^2_{\ell_2}\geq \tilde{C}\|x^n-x^{n+1}\|^2_{\ell_2}
\end{align*}
%
 As $\|\mathcal{J}(x^n,w^n,\epsilon_n)-\mathcal{J}(x^{n+1},w^{n+1},\epsilon_{n+1})\|^2_{\ell_2} \rightarrow 0$ as $n\rightarrow \infty$ we also get 
 $$\|x^n-x^{n+1}\|^2_{\ell_2} \rightarrow 0 \text{  as  }n\rightarrow \infty.\blacksquare$$

We are now ready to present the result of convergence which we concisely summarize as follows: As $x^{n}$ is a bounded sequence,  we obtain either the exact minimizer of $\|A(x)-y\|^p_{\ell_p}$ or every cluster point is a critical point of the $\epsilon$-perturbed $\ell_p$-norm residual $f_\epsilon(x)$ at least. 

\begin{theorem}\label{mainth}
Fix $y\in \mathbb{R}^m$, $x^0 \in \mathbb{R}^k$. Let $A:\mathbb{R}^k\rightarrow\mathbb{R}^m$ be a nonlinear map with $A\in \mathcal{C}^1$ for which the boundedness and coercivity condition (BCC) holds, i.e., there exist $\alpha,\beta>0$ such that, for all $z\in\mathcal{B}(0,R^*)$:
\begin{equation*}
\alpha \|x^*-z\|_{\ell_2}\leq \|A(x^*)-A(z)\|_{\ell_p}\leq \beta \|x^*-z\|_{\ell_2}.
\end{equation*}
{ Additionally we require that, for $(x^n)_{n \in \mathbb{N}}$ and $(w^n)_{n \in \mathbb{N}}$ sequences generated by Algorithm 2,
\begin{align} \label{Lipcond2}
& \left |  t   [\|A(tx^n+(1-t)x^{n+1})-y \|_{\ell_2(w^n)}^2  - \|A(x^n)-y \|_{\ell_2(w^n)}^2 ] \right .  \\  \nonumber 
& \left .  +(1-t)[  \|A(tx^n+(1-t)x^{n+1})-y \|_{\ell_2(w^n)}^2 - \|A(x^{n+1})-y \|_{\ell_2(w^n)}^2]  \right | \\ \nonumber
&\leq L t (t-1) \|x^{n} -x^{n+1}\|_{\ell_2}^2.
\end{align}
for some $L>0$ indepedent of $n \in \mathbb N$ and for all $t \in [0,1]$. For $\omega>0$ large enough (and determined according to Lemma \ref{USCC12}) we have the following properties of  Algorithm 2:} 
\begin{enumerate}
\item[(i)] If $\epsilon=\lim\limits_{n\rightarrow\infty} \epsilon_n=0$, then the sequence $(x^n)_{n \in \mathbb{N}}$ converges to a vector $\bar{x}$, which is the solution to the $\ell_p$-minimization problem \eqref{ellp}. Moreover if $y\in Ran(A)$ and $y=A(x^*)$ then $x^*$ is the unique minimizer, thus $\bar{x}$ coincides with $x^*$.
\item[(ii)]if $\epsilon=\lim\limits_{n\rightarrow\infty} \epsilon_n>0$, then all accumulation points of $(x^n)_{n \in\mathbb{N}}$ are critical points of the $\epsilon$-perturbed $\ell_p$-norm residual $f_\epsilon$ as defined in \eqref{xeps} and they all belong to $\mathcal{B}(0,R^*)$.
\end{enumerate}
\end{theorem}

Proof: 
\begin{itemize}
\item[(i)] { the proof of this statement is analogous to the one of \thref{conv}.}
\item[(ii)] 
We already mentioned that $(x^n)_{n\in \mathbb{N}}$ is a bounded sequence in $\mathcal{B}(0,R^*)$ and hence it has accumulation points. Let $(x^{n_\ell})_{\ell\in \mathbb{N}}$  be any convergent subsequence of $(x^n)_{n\in \mathbb{N}_0}$  and let $\bar{x}$ be its limit. We want to show that $\bar{x}$ is a critical point of \eqref{xeps}. \\
Since $w^n_i=[(A_i(x^n)-y_i)^2+\epsilon_n^2]^{(p-2)/2} \leq \epsilon_n^{p-2} \leq \epsilon^{p-2}$, it follows (up to extracting an additional subsequence) that $\lim\limits_{\ell\rightarrow \infty} w^{n_\ell}_i= [(A_i(\bar{x})-y_i)^2+\epsilon^2]^{(p-2)/2}={w}(\bar{x},\epsilon)_i:=\bar{w}_i, i=1,\dots, m$. \\
On the other hand, by invoking \thref{chiave} , we obtain also that $x^{n_\ell+1}\rightarrow \bar{x}, \ell\rightarrow\infty$. { (Notice that here $\epsilon_n \geq \epsilon >0$ and the considerations in Remark \ref{remrem} can be applied in order to justify the assumption \eqref{Lipcond2}.)}
Analogously $w^{n_\ell+1}\rightarrow \bar{w}$ for $\ell\rightarrow\infty$. Since we assume $A \in \mathcal C^1$ then $x \to \mathcal{J}(\cdot ,w^n,\epsilon_{n})$ is actually differentiable. We observe that by \eqref{inner},
$$0  =\nabla_x\mathcal{J}_{\omega,x^{n_\ell}}(x^{n_\ell+1},w^{n_\ell},\epsilon_{n_\ell})=\nabla_x \mathcal{J}(x^{n_\ell+1},w^{n_\ell},\epsilon_{n_\ell})+2\omega(x^{n_\ell+1}-x^{n_\ell})$$ or
$$ -2\omega(x^{n_\ell+1}-x^{n_\ell}) = \nabla_x\mathcal{J}(x^{n_\ell+1},w^{n_\ell},\epsilon_{n_\ell}).$$
Using \thref{chiave} we can conclude that taking the limit $\ell \to \infty$ gives 
$$ 0 \in\nabla_x \mathcal{J}(\bar{x},\bar{w},\epsilon)= \nabla f_\epsilon(\bar{x}). \blacksquare$$ 
(Let us emphasize that, due to our assumption $A \in \mathcal C^1$, in all the steps above, the functions we consider are differentiable. We could consider to lower the smoothness by requiring $A \in \mathcal C^0$ and ask additional properties of subdifferentials, such as outer semicontinuity etc. However, this generalization towards nonsmooth analysis brings little truly additional insights and we keep content with the current formulation.)
\end{itemize}

\begin{remark}
The error decay result (iii) in \thref{conv} stays valid also for Algorithm 2 if condition (c) in \thref{conv} is fulfilled. 
\end{remark}

\section{Numerical Experiments}
We want to illustrate our theoretical results by several numerical experiments. In this section we shall first test the developed algorithms in a simple case to get some intuition of their behavior before trying to apply them on more involved higher dimensional $\ell_p$-minimization problems, whose optimal solution as we know is often not so easy to determine. In a first exemplary case that is computationally easy to handle and visually presentable we study the development of the iterates of  NR-IRLS step by step and compare the algorithm output with standard Matlab optimization routines. Furthermore we would like to test the validity of our theoretical results also in more complex situations and we check them via the correct reconstruction of sparse vectors in the context of nonlinear compressed sensing problems as studied in \cite{Quasi14}. We involve the NR-IRLS in both its original and convexified versions in intermediate steps of a greedy sparse recovery algorithm. If the overall sparse recovery algorithm gives correct results the intermediate steps must have been performed correctly as well. 
Finally we consider the context of data corrupted by impulsive noise, where the sparsification of the residual is desired, and observe the recovery success rates for different levels of severity of the impulsive noise perturbations.

Before we present the results of the numerical tests in more detail, there are some important numerical issues that should be mentioned in the context of IRLS-type methods:
\begin{itemize}
\item Assume that we are in the case where $\epsilon_n\rightarrow 0$, and thus $|A(x^n)_i-y_i|\rightarrow 0$ for all $i\in\left\{1,\dots,m\right\}$ for $n\rightarrow \infty$ due to the definition of $\epsilon_n$.
Thus in this case it will be unavoidable that $(w_n)_i$ grows until reaching the limits of  machine representation. On the one hand practically there is the need of fixing a lower bound $\hat{\epsilon}>0$ to go around this issue. 
On the other hand by introducing $\hat{\epsilon}$ we limit the capability of our IRLS algorithm: it will only be able to approximate the correct solution but not to find it exactly, which does not necessarily match the theoretical analysis of course. In order to obtain sufficient recovery accuracy it is necessary to choose the value $\hat{\epsilon}$ appropriately small. Unfortunately this issue leads to numerical complications again: choosing e.g. $\hat{\epsilon}=1e^{-8}$ and considering approximation residuals $|A(x^n)_i-y_i|\ll1$, then $(w_n)_i$ is of the order of $1e^{+8}$ and thus multiplication or addition can lead to numerical errors that are not negligible and will definitely affect the calculation results. As a consequence we have to face limitations on the recovery accuracy that will most probably not allow results with errors in the range of machine precision. 

\item Additionally to the errors resulting from IRLS itself we are using iterative methods for solving the internal locally convex optimization problem in each step. From that we will have an approximation error determined by the specific termination tolerance of the chosen method lowering further the expected accuracy. 
\end{itemize}

All numerical experiments part of this paper were performed on a MacBook Pro 9.1.
with a 2.6 GHz Intel Core i7 quad-core-processor and 8GB memory.
Computations were run in MATLAB R2012b version 8.0.0.

\subsection{Simple Example for $A:\mathbb{R}\rightarrow \mathbb{R}^2$}

As an introductory test example we consider $$A:\mathbb{R}\rightarrow \mathbb{R}^2, x\mapsto \begin{pmatrix}x\\x^2\end{pmatrix}$$ and a measurement vector $ y\in [0,1]^2$, hence the $\ell_p$-minimizer $x^*:=\argmin\limits_{x}\|A(x)-y\|^p_{\ell_p}$ will be in $[0,1]$ as well. \\
As a first step we check that in this situation the BCC as in \thref{BCC} is fulfilled for $1<p<2$.
We verify that a choice for the lower BCC-bound $\alpha$ is just $1$:
\begin{align*}
\|A(x)-A(x^*)\|_{\ell_p}=(|x-x^*|^p+|x^2-(x^*)^2|^p)^{1/p}\geq  |x-x^*|\geq\alpha \|x-x^*\|_{\ell_2}
\end{align*}
For the upper bound $\beta$ we obtain $(1+2^p)^{1/p}$:
\begin{align*}
\|A(x)-A(x^*)\|_{\ell_p} &=(|x-x^*|^p+|x^2-(x^*)^2|^p)^{1/p}= ( |x-x^*|^p+|x-x^*|^p\cdot|x+x^*|^p)^{1/p}\\
&\leq ( |x-x^*|^p+|x-x^*|^p\cdot2^p)^{1/p}=(1+2^p)^{1/p}|x-x^*|= \beta \|x-x^*\|_{\ell_2}
\end{align*}
In general we have a nonconvex problem in our situation with multiple local minimizers and want to analyze the behaviour of our original version of the NR-IRLS, Algorithm 1. We are interested in the influence of the parameter $p$ . Of course, changing $p$ means also the change of the problem and a different minimizer as well as the appearing of possible local minimizers.
 In this context it is interesting to examine the different minimization results for different $p$ and different starting points. In the following we observe the behavior of the NR-IRLS applied to $\ell_p$-minimization problem for this measurement setting and compare it to the solution obtained with Matlab's $lsqnonlin$-function, which  basically implements a trust-region-reflective or Levenberg-Marquardt strategy, applied directly to the $\ell_p$-minimization problem. \\

In the following, we start with a more detailed description of the concrete test setting.
We consider for measurements $y=(0,0.9)^T$ and observe the recovery results for different values of $p$ ranging between $1$ and $2$, more precisely $p \in \{1.1,1.3,1.7,1.9\}$.  \\
For the specific setting of the algorithm parameters we choose the maximum number of iterations of NR-IRLS itself to be 50.  For the execution of the locally convex minimization in each inner step we chose the MATLAB built-in function $fminunc$ with default settings and the last iterate as the starting point. \\
Moreover  we use MATLAB's default settings as well for running the $lsqnonlin$-function applied directly to the $\ell_p$-minimization problem.\\
For both procedures we try different starting points $x^0$ in the interval $[0,1]$, namely $x^0 \in \{0,0.25,0.5,0.75,1\}$ and observe the convergence to different local minimizers.\\

By means of the graphical analysis, we are able to extract the following observations: Figure 1 and 2 demonstrate the property of NR-IRLS to converge to the local minimizer of the objective function that is closest to the $\ell_2$-critical point obtained in the first step of the algorithm for all values of $p$ depending on the starting point $x^0$ in the first step. \\
Moreover the influence of the starting point for solving the first nonlinear least squares problem becomes obvious in Figures 3 and 4. The MATLAB method converges to the minimizer closest to the starting point while the NR-IRLS converges to the minimizer closest to the results of the $\ell_2$-minimization in the first step, which can be different. Hence NR-IRLS is able to find different local minimizers than standard gradient based methods for the same starting point of the minimization.

\begin{figure}[h!]
  \caption{}
  \centering
    \includegraphics[width=1\textwidth]{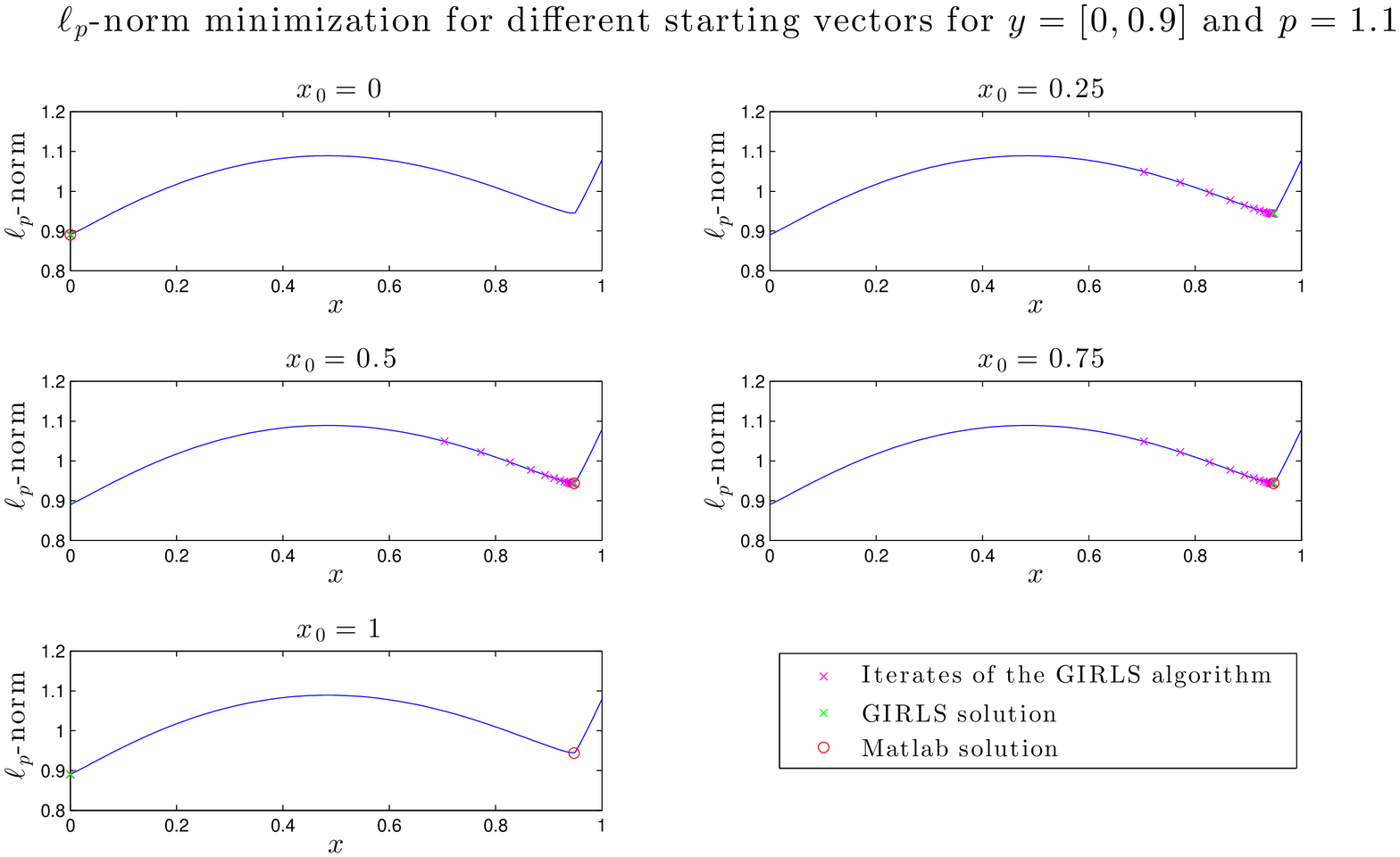}
\end{figure}
\begin{figure}[h!]
  \caption{}
  \centering
    \includegraphics[width=1\textwidth]{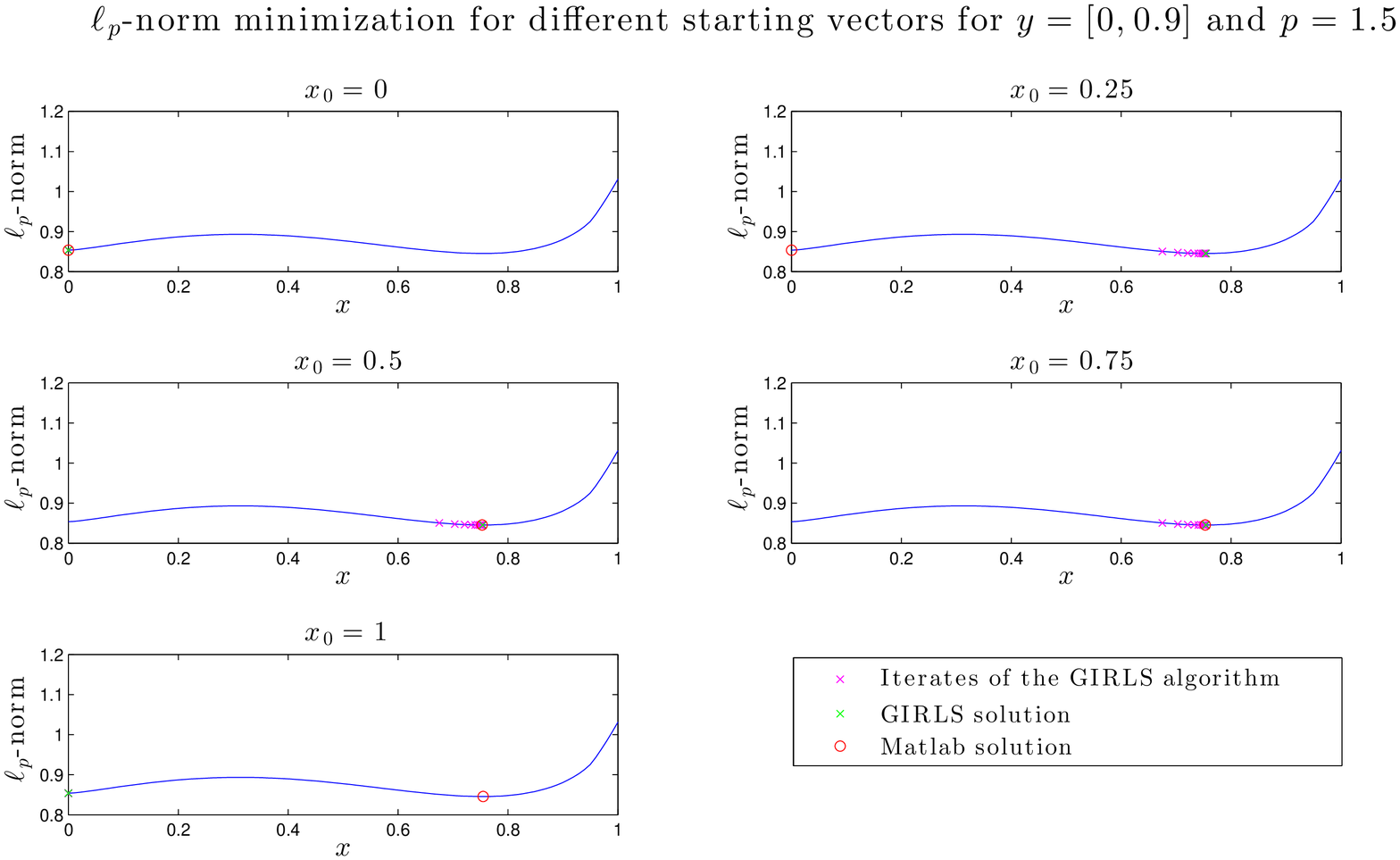}
\end{figure}
\begin{figure}[h!]
  \caption{}
  \centering
    \includegraphics[width=1\textwidth]{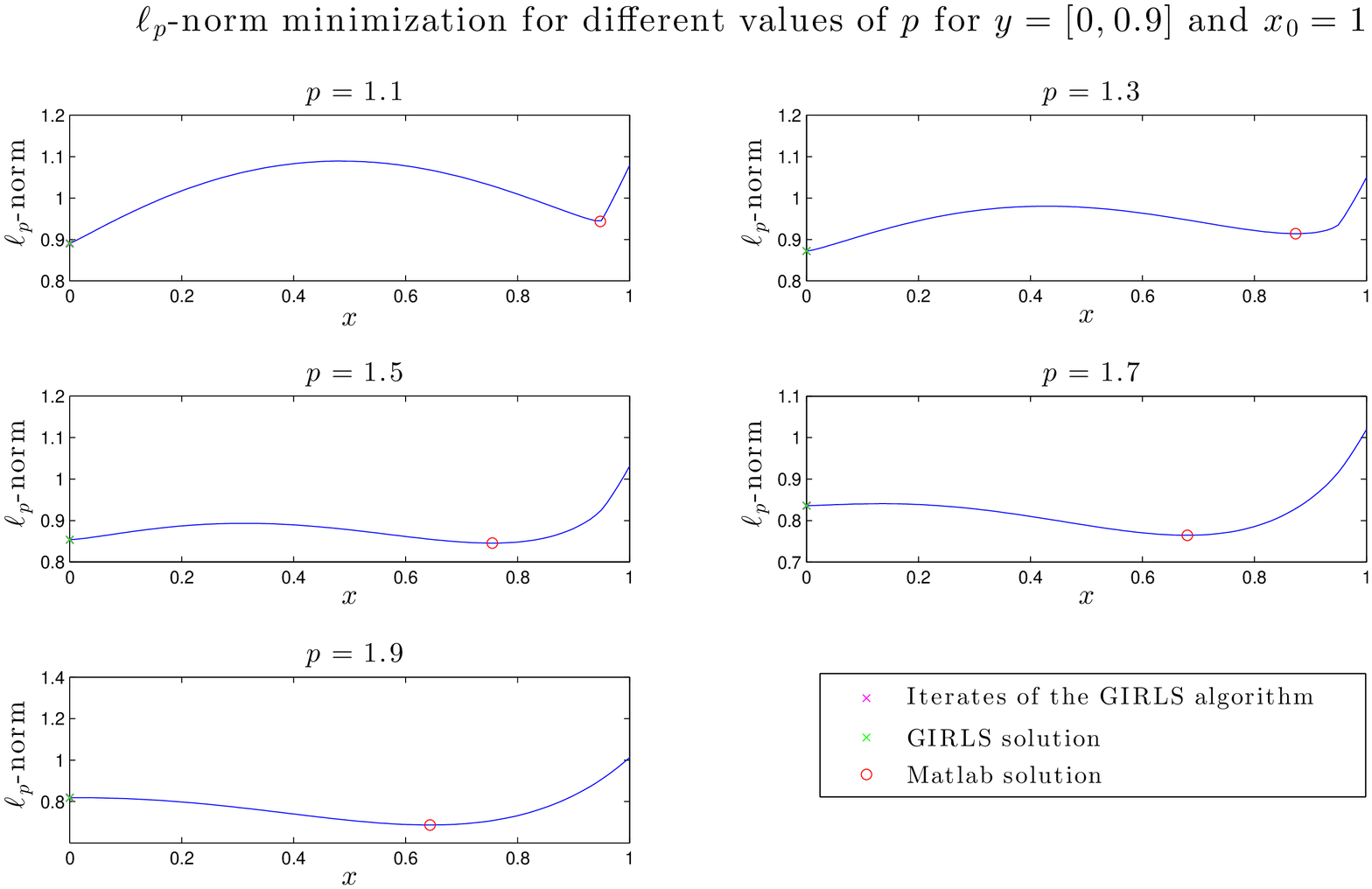}
\end{figure}
\begin{figure}[h!]
  \caption{}
  \centering
    \includegraphics[width=1\textwidth]{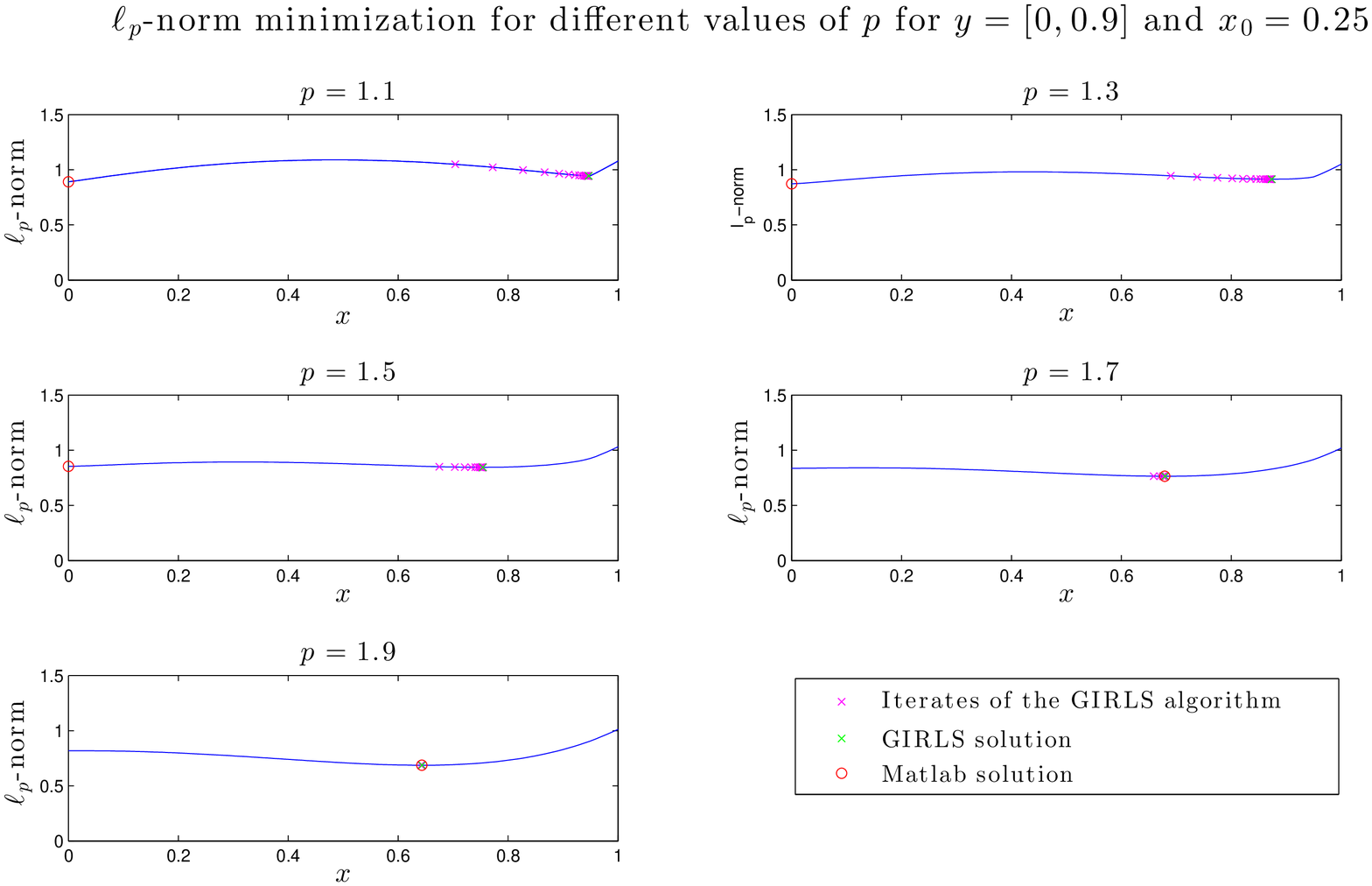}
\end{figure}

\newpage
$\left. \right.$
\newpage
\subsection{High Dimensional Examples in a Nonlinear Compressed Sensing Application Context}
In \cite{Quasi14} we introduced a greedy algorithm (Algorithm 1) for recovering sparse vectors from a  small number of nonlinear measurements (we may like to name such class of problems, expecially if the measurements are generated randomly, {\it nonlinear compressed sensing}). One of the crucial operations of this iterative algorithm is the $k$-variate nonlinear $\ell_p$-minimization problem \eqref{lp}. Roughly speaking, the algorithm searches at its $k$-th step for the vector with at most $k$ nonzero entries which best fits the data in terms of a minimal norm nonlinear residual problem of the type \eqref{lp}. In the problem class considered in \cite{Quasi14}, it was also interesting to consider $p \in [1,2]$ as a norm parameter in a  BCC-like conditions (there named Restricted Isometry Property (RIP)), see, more specifically, formula (3.1) in \cite{Quasi14}. 
In what follows we consider nonlinear functions $A:\mathbb R^k \to \mathbb R^m$ taken as restrictions to $k$-dimensional index subspaces of two types of maps considered in \cite{Quasi14}. The first type are maps which are Lipschitz perturbations of matrices fulfilling
the RIP. 
The second setting refers to quadratic maps $A$ relative to phase retrieval models \cite{Drenth07,Fienup82,Gerchberg72}. 
For both these setting we present the corresponding numerical recovery results using NR-IRLS Algorithm 1 and its locally convexified version Algorithm 2, and we compare them with standard tools of Matlab.\\
An implementation of the greedy algorithms for nonlinear compressed sensing is available at \url{http://www-m15.ma.tum.de/Allgemeines/SoftwareSite}.

\subsubsection{Locally Convex Case: Nonlinear Perturbation of Linear RIP-Matrices}
In \cite[Section 3.2.1]{Quasi14} the following result was obtained.

\begin{proposition}\thlabel{pert}
Assume $k\leq m \leq N$ and ${A_1}\in\mathbb{R}^{m\times N}$ satisfies the $\delta$-RIP of order $2 k$, i.e.,
$$
(1-\delta) \|z\|_{\ell_2^N} \leq \| A_1 z \|_{\ell_2^m} \leq (1+\delta)\|z\|_{\ell_2^N}, 
$$
for all $z \in \mathbb R^N$ with at most $2 k$ nonzero entries. If ${A}_\rho: \mathbb{R}^N\rightarrow\mathbb{R}^m$ is chosen as
\begin{equation}
{A}_\rho({z}):={A_1}{z}+\rho {f}(\|{z}-z^\circ\|_{\ell_2}^2){A_2}{z}, \label{firstqn}
\end{equation}
where $z^\circ\in\mathbb{R}^N$ is some reference vector in $\mathbb R^N$, $f:[0,\infty)\rightarrow\mathbb{R}$ is a  bounded Lipschitz continuous function  with $f(0)=0$, $\rho>0$ is a sufficiently small scaling factor, and ${A_2}\in\mathbb{R}^{m\times N}$ arbitrarily fixed, then there are constants $\alpha,\beta>0$, such that for $p=2$ 
$$  \alpha \|{z}-z^*\|_{\ell_2^N}\leq \|{A}_\rho({z})-{A}_\rho(z^*)\|_{\ell_p^m}\leq \beta \|{z}-z^*\|_{\ell_2^N}$$
 for all $z$ with at most $k$ nonzero entries and $z^*$ is another fixed vector of at most $k$ nonzero entries.
For other $p \in [1,2)$ these inequalities hold again with different constants $\alpha,\beta$, derived, for instance, by equivalence of norms:  for $0<r<q$ we  have $\|z\|_{\ell_q} \leq \|z\|_{\ell_r} \leq N^{1/r-1/q}\|z\|_{\ell_q}$.
\end{proposition}

From this proposition we deduce that any restriction of $A_\rho$ to vectors supported on a certain fixed set of indexes $\Lambda \subset \{1,\dots, N\}$ with $\# \Lambda =k$ satisfies the BCC. For the purposes of this paper, without loss of generality
we can simply assume $\Lambda =\{1,\dots, k\}$ and we can define
 $$A:\mathbb{R}^k\times \mathbb{R}_+ \rightarrow \mathbb{R}^m, (x,\rho)\mapsto A(x,\rho)={A}_\rho({x}^\Lambda),$$ 
where $z=x^\Lambda$ is the zero padding extension of $x$ to a vector in $\mathbb R^N$.

By \thref{full} in the Appendix we have that the first USCC required in \thref{conv} is fulfilled for the linear case of $A(\cdot,0)$, i.e., for $\rho=0$ and $A(x,0) = (A_1)_{|\Lambda}$ reduces to a matrix in $\mathbb R^{m \times k}$.
As $\rho>0$ is small $A(\cdot,\rho)$ is actually just a small nonlinear perturbation of $A(\cdot,0)$. If we additionally assume now that $f \in \mathcal C^2(\mathbb R_+)$ as appearing in the definition of $A_\rho$, then by a rather straightforward continuity argument we can easily extend the first USCC to $A(\cdot,\rho)$ on a small ball around $x^*$. We omit the explicit, tedious, and perhaps rather clear elaboration of this argument.
\\

The next numerical examples are addressed to the problem of recovering a sparse vector $z^* \in \mathbb R^N$ (with at most $k$ nonzero entries, for $k \in [1,10] \cap \mathbb N$ from the given measurements $y=A_\rho(z^*)$. 
We employ \cite[Algorithm 1]{Quasi14}, which requires at each step the solution of a minimal norm nonlinear residual. To perform such optimization we utilize Algorithm 1 of the present paper.
Let us describe the precise setting of the experiments. We fixed the dimension $N=80$, the number of measurements $m=30$ and we draw at random RIP matrices $A_1$ having i.i.d. Gaussian entries and we fix $A_2$ as the matrix with one as an entry everywhere. The function $f$ is chosen to be the squared Euclidean distance from the solution $ f(\|z-z^*\|^2_{\ell_2}) = \|z-z^*\|^2_{\ell_2^N}$.  As pointed out above the fulfilling of the BCC and the USCC depends on a small parameter $\rho>0$, which is indicating the severity of the nonlinearity. We perform experiments with  $\rho \in \{0,0.5,1,3,5,10,20\}$ to observe the influence of this parameter and how the success rate may depend on the nonlinearity.
For each of these parameter combinations we randomly generate a set of 100 synthetic problems.\\
Since we work with synthetic problems, we already have the expected sparse minimizer $z^*$, which can be used to determine the success of the recovery. In particular we claim successful reconstruction when the error is within a $1\%$ of the solution's norm. Possible additional noise on the measurements is here not yet considered in all examples. Let us mention that the greedy algorithm itself was allowed to perform $3k$ steps (hence a larger amount of steps with respect to the expected dimension of the solutions), to give it the opportunity to correct indices, which may have been wrongly added to the support. Furthermore we choose the maximum number of iterations of NR-IRLS to be 50. 
As before, for the execution of the locally convex minimization in each inner step we chose the MATLAB built-in function $fminunc$ with default settings and the origin as the starting point.\\

\begin{figure}[H]
\centering
\subfigure{
{\includegraphics[width=.48\textwidth]{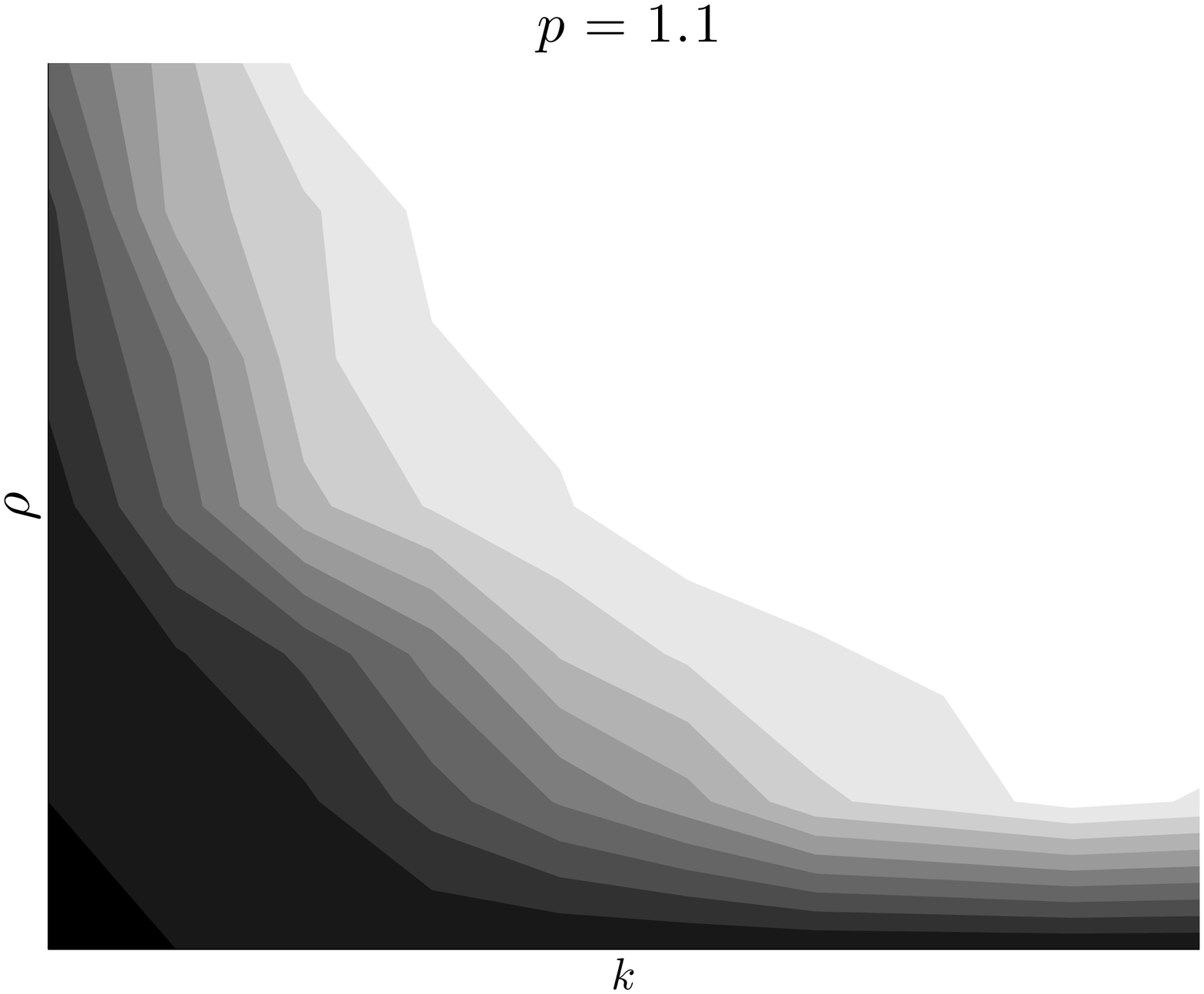}\label{sub:2}}}
\subfigure{
{\includegraphics[width=.48\textwidth]{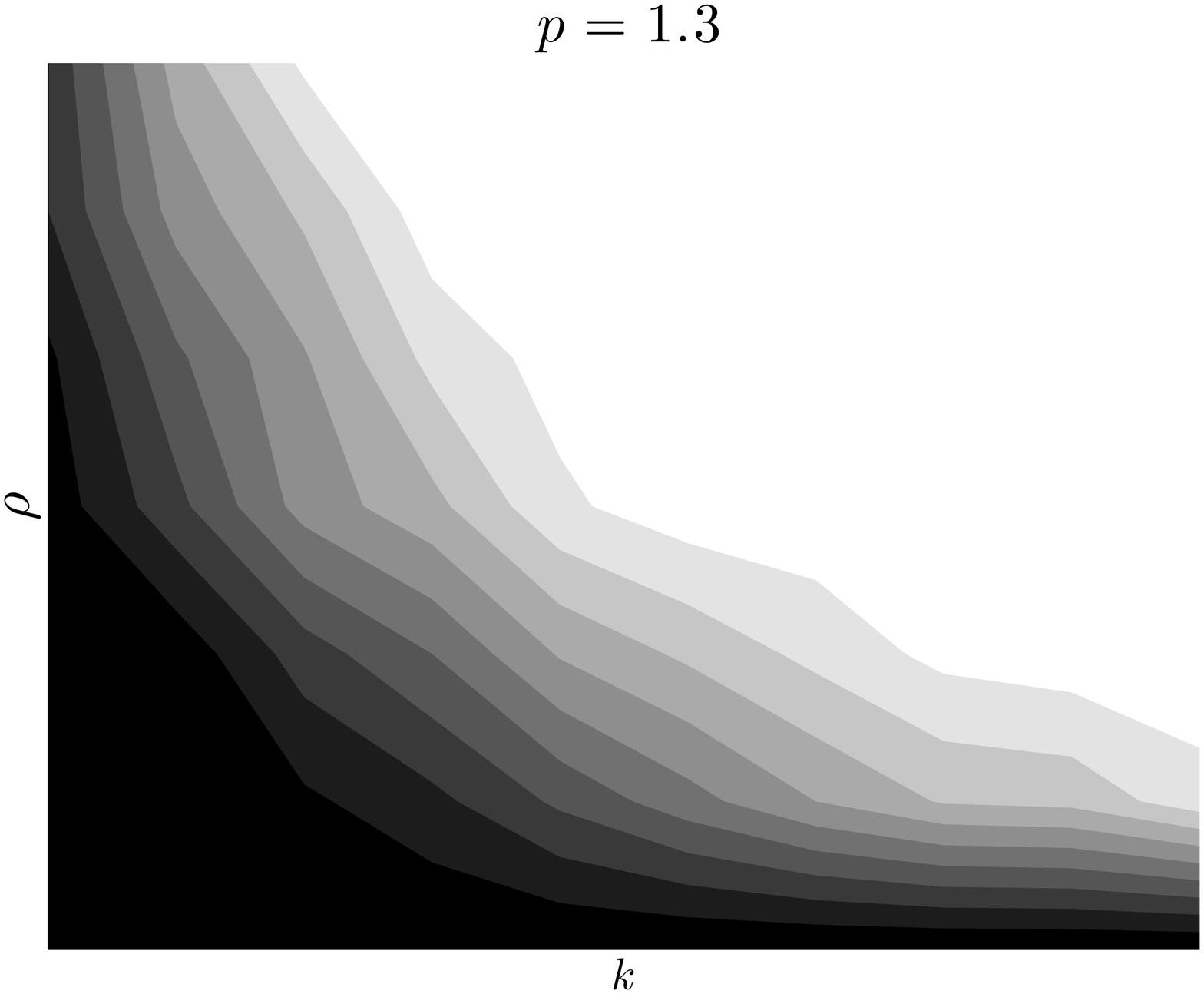}}\label{sub:1}}
\subfigure{
{\includegraphics[width=.48\textwidth]{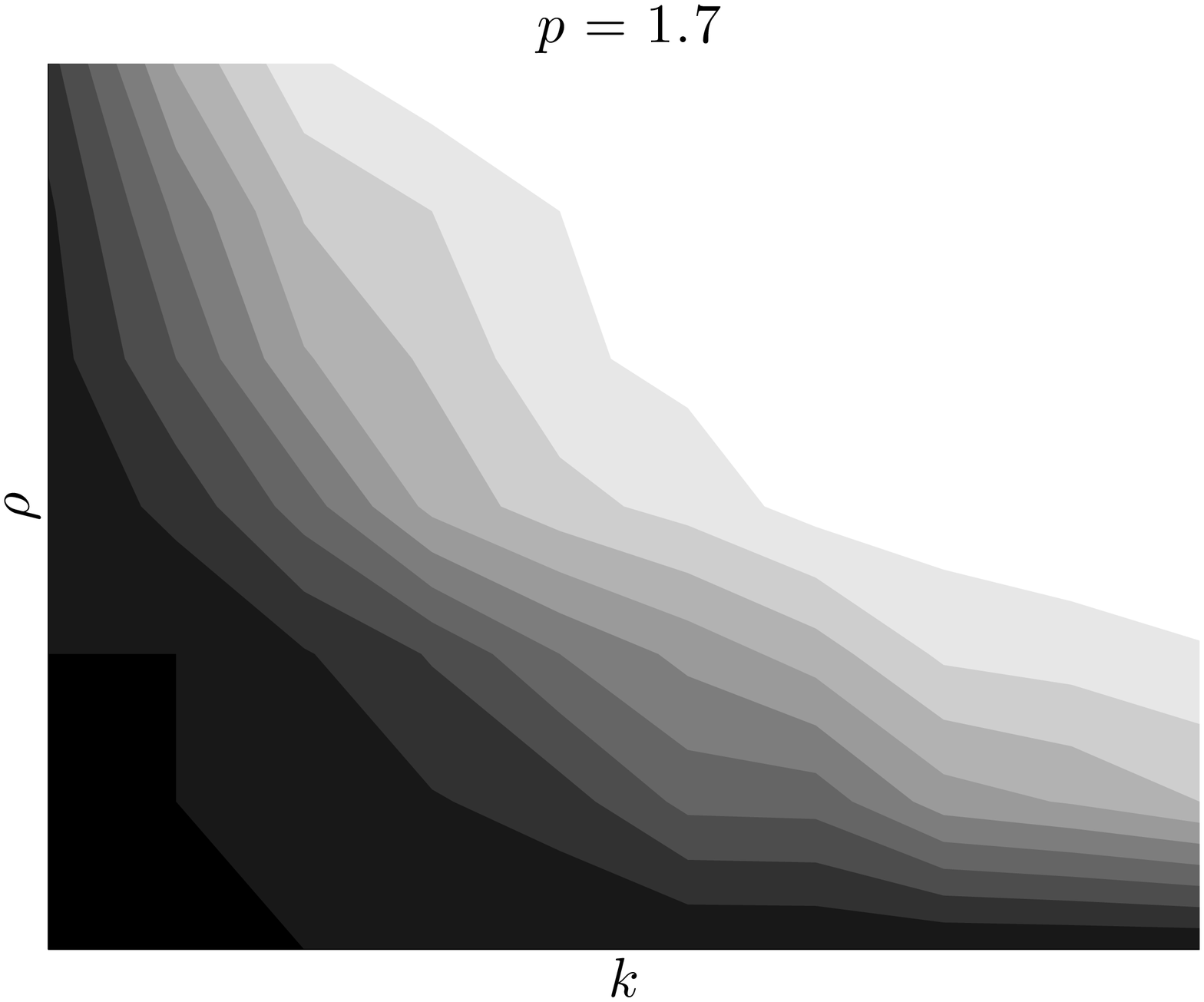}\label{sub:2}}}
\subfigure{
{\includegraphics[width=.48\textwidth]{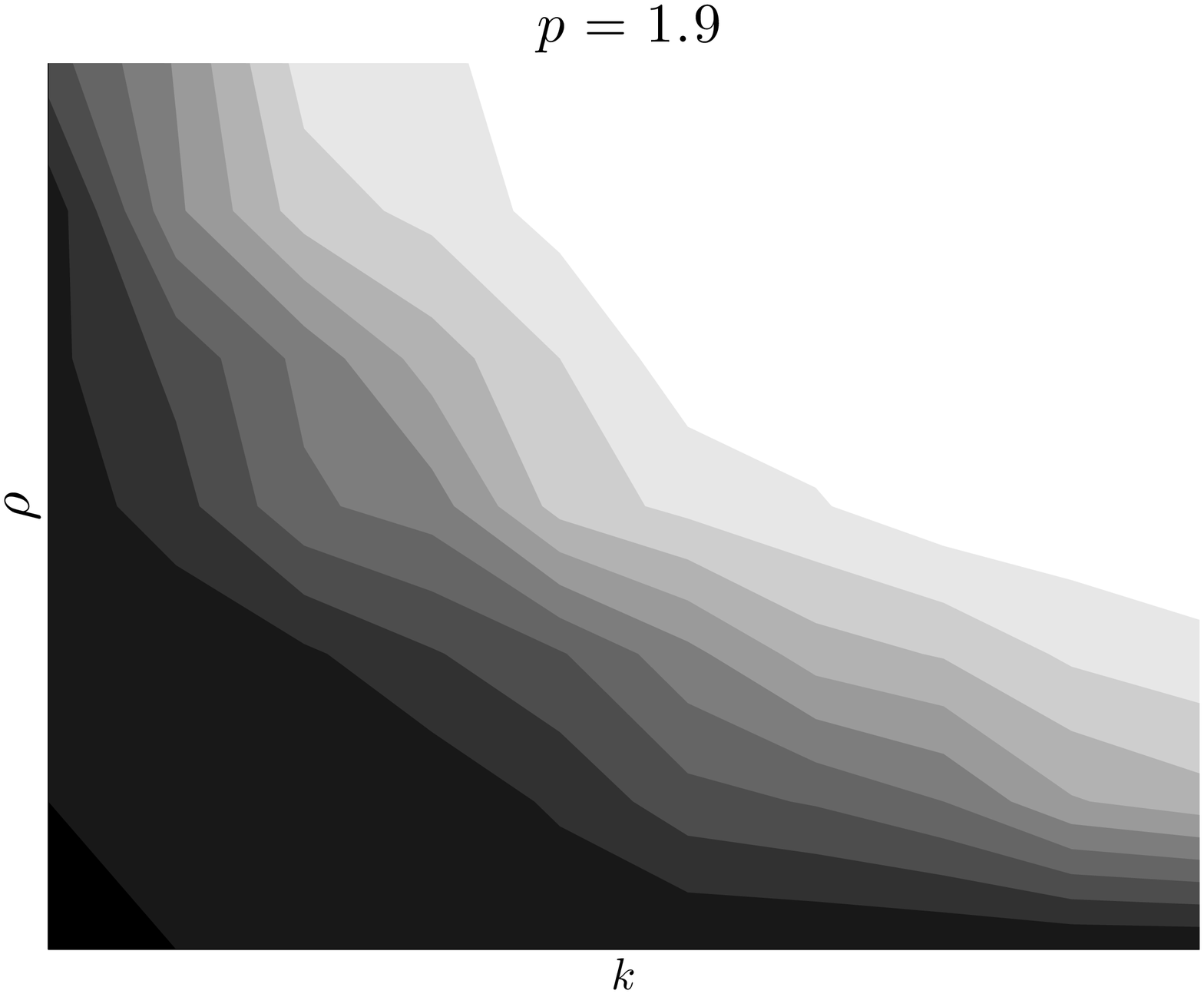}}\label{sub:1}}

\caption{Recovery rates for the greedy strategy developed in \cite{Quasi14} used with the RIP matrix perturbation measurements as above with $N = 80, m = 30$, $A_1$ having i.i.d. Gaussian entries, $A_2$ being the matrix with ones everywhere, $ f(x) = \|x-x_0\|^2_{\ell_2}$, and we use solutions $x^*$ with $\|x^*\|_{\ell_2}=0.015$ . Reconstruction is repeated 50 times for each signal and $k$ to derive stable recovery rates.}
\label{fig:perRIP}
\end{figure}

The plots in Figure \ref{fig:perRIP} show the recovery rates or the empirical probability of successful sparse vector recovery of \cite[Algorithm 1]{Quasi14} implementing  Algorithm 1 of the present paper for performing the $\ell_p$-minimization for different values of $p$. 
As expected, the recovery rates decrease with growing $k$, hence of the dimensionality of the $\ell_p$-minimization problem, and growing Lipschitz perturbation factor $\rho>0$ of $A_\rho$, representing the degree of nonlinearity. We obtain better recovery results for $p$ closer to $2$  probably because the BCC condition in this case has tighter constants $\alpha$ and $\beta$.

\subsubsection{Phase Retrieval Problem}

We consider here as in \cite{Quasi14} a sequence of real Gaussian random vectors $a_i\in\mathbb{R}^N$, $i=1,\ldots,m$ and the nonlinear measurement map
\begin{equation}\label{eq:model phase}
{A}({x})= (|\langle a_1,{x}\rangle|^2,\ldots,|\langle a_m,{x}\rangle|^2)^\top. 
\end{equation}
By \cite[Theorem 3.12]{Quasi14} and according to \cite[Formula (3.14)]{Quasi14} there are constants $\alpha,\beta>0$, such that another BCC-like property, condition in \cite[formula (3.8)]{Quasi14}, holds for $p=1$ (with slight adaption including the Hilbert-Schmidt norm instead of the $\ell_2$-norm on the left- and right hand side of the inequality which does not influence significantly the results above). However, the first USCC is not fulfilled in general for this kind of problem, thus we need to employ the convexified version of NR-IRLS, Algorithm 2. We again fixed the dimension of the signal $N=80$, the number of measurements $m=30$ and i.i.d. Gaussian random vectors $a_i$, $i=1,\dots,m$. Moreover we created synthetic solutions $z^*$ with $\|z^*\|_{\ell_2}=1$ and respective sparsity level $k\in \{1,3,6,9,12,15,18,21\}$. The vectors were constructed such that the nonincreasing rearrangement of the absolute value of their entries satisfies the decay rates $\kappa \in \{1,0.8,0.6,0.4\}$, see the vector class $\mathcal D_\kappa$ in  \cite{Quasi14} for a precise definition. Such decay is required as a sufficient conditions in the convergence results of \cite{Quasi14} and it was verified to be crucial in the numerical experiments done with standard Matlab optimization routines. For each of these parameter combinations we generate a set of 50 synthetic problems that do not consider the occurrence of noise. 
As above we use the solution $z^*$, the expected sparse solution, to determine the success of the recovery, and claim again successful reconstruction when the error is within a $1\%$ of the solution's norm. The algorithmic settings are similar as above: $3k$ steps are  performed by the greedy algorithm \cite[Algorithm 1]{Quasi14}, the maximum number iterations of NR-IRLS is allowed to reach 100. 
The regularization parameter $\omega>0$ is set to 100 and for the execution of the convex minimization in each inner step we choose the MATLAB built-in function $fminunc$ with default settings and random starting points with a norm smaller or equal to the solutions norm.\\

The plots in Figures \ref{fig:PR}-\ref{fig:PR2} show the recovery rates of  \cite[Algorithm 1]{Quasi14} of sparse vectors from measurements of the type \thref{pert}, implementing Algorithm 2 of the present paper for performing the $\ell_p$-minimization for different values of $p$. 
Surprisingly the decay rate of the nonincreasing rearrangement of the solution vector $z^*$ seems not to have notable influence on the recovery results when using NR-IRLS, although we experienced this phenomenon in our earlier paper \cite{Quasi14}, where we used the built-in MATLAB functions $fminunc$, $fminsearch$ or $lsqnonlin$ for solving the internal $\ell_p$-minimization problem.  For comparison we give the corresponding results obtained from an implementation of the greedy algorithm using $lsqnonlin$ for solving the internal $\ell_p$-minimization problem as well. Its results are obviously outperformed by NR-IRLS when the decay rate of the nonincreasing rearrangement of the solution is not sufficiently pronounced.
\begin{figure}[H]
\centering
\subfigure{
{\includegraphics[width=.48\textwidth]{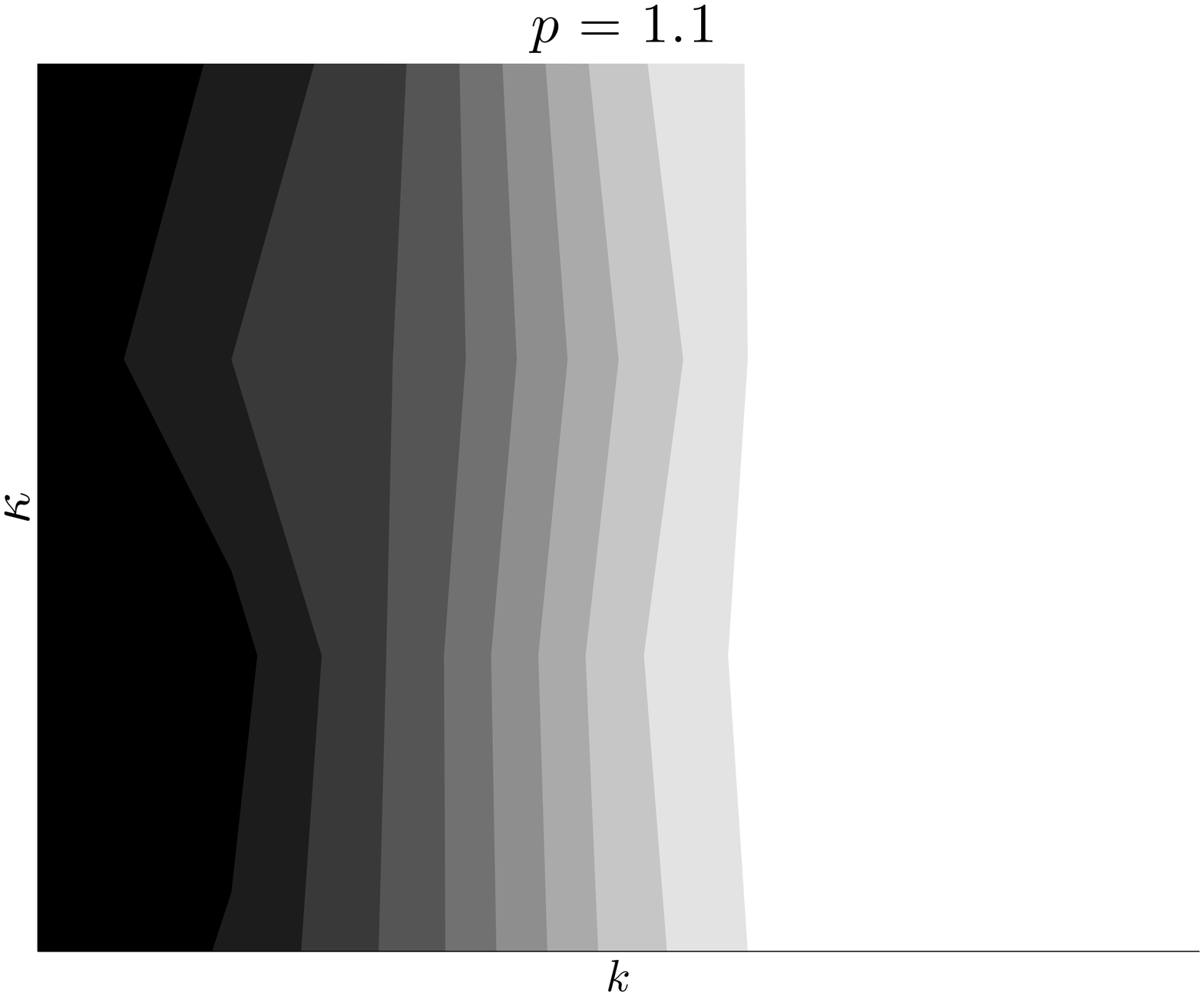}\label{sub:2}}}
\subfigure{
{\includegraphics[width=.48\textwidth]{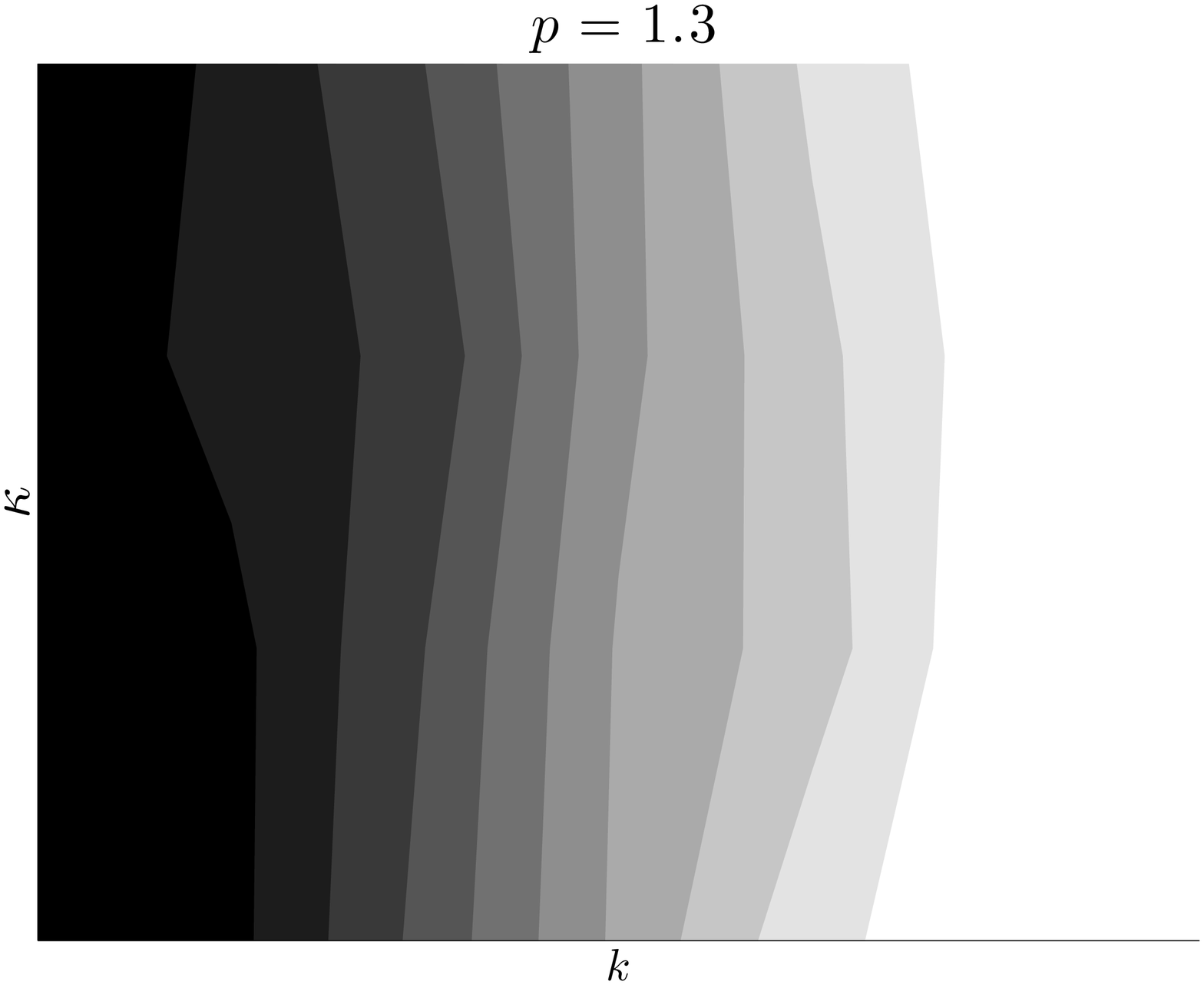}}\label{sub:1}}
\end{figure}
\begin{figure}[H]
\centering
\subfigure{
{\includegraphics[width=.48\textwidth]{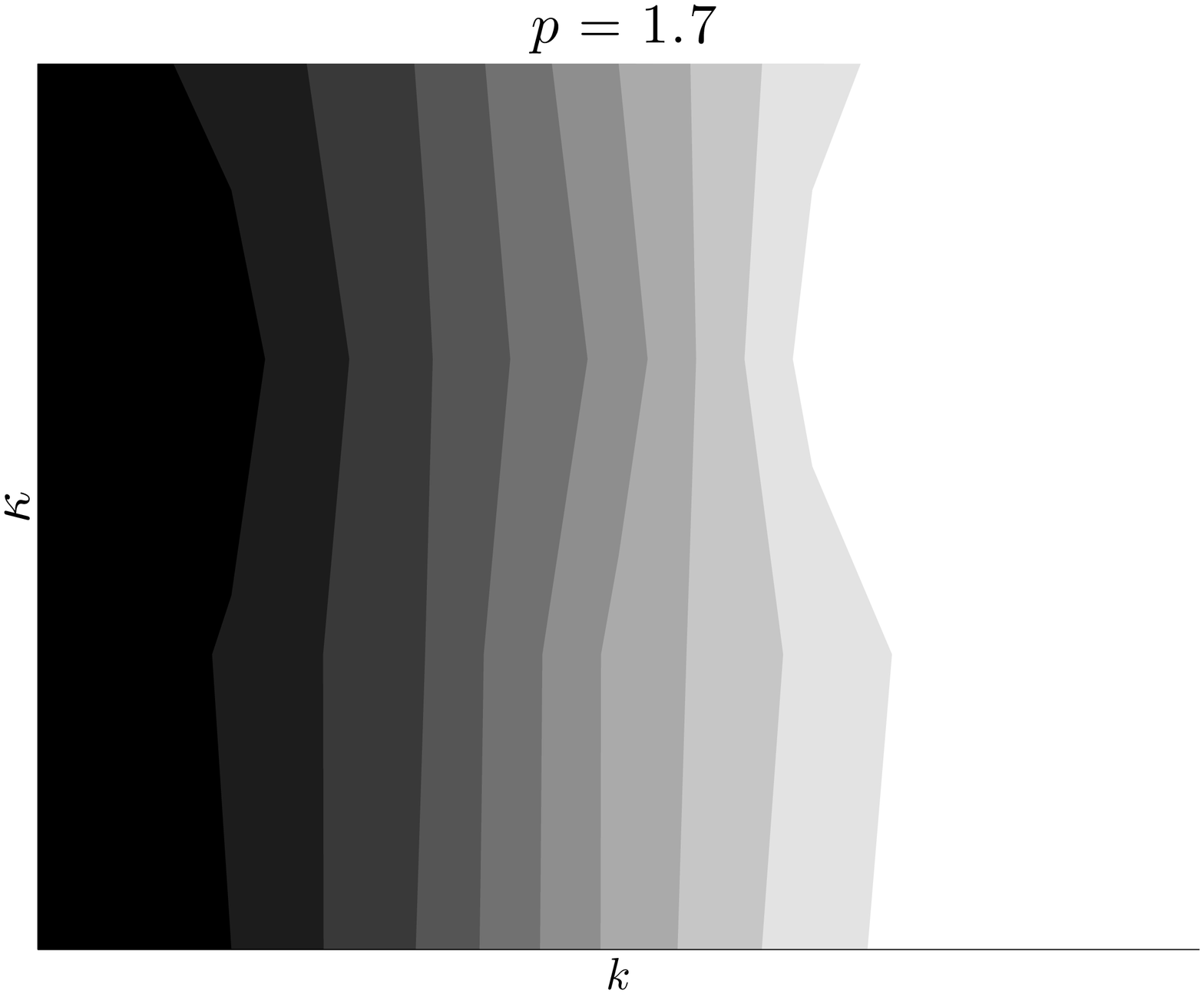}\label{sub:2}}}
\subfigure{
{\includegraphics[width=.48\textwidth]{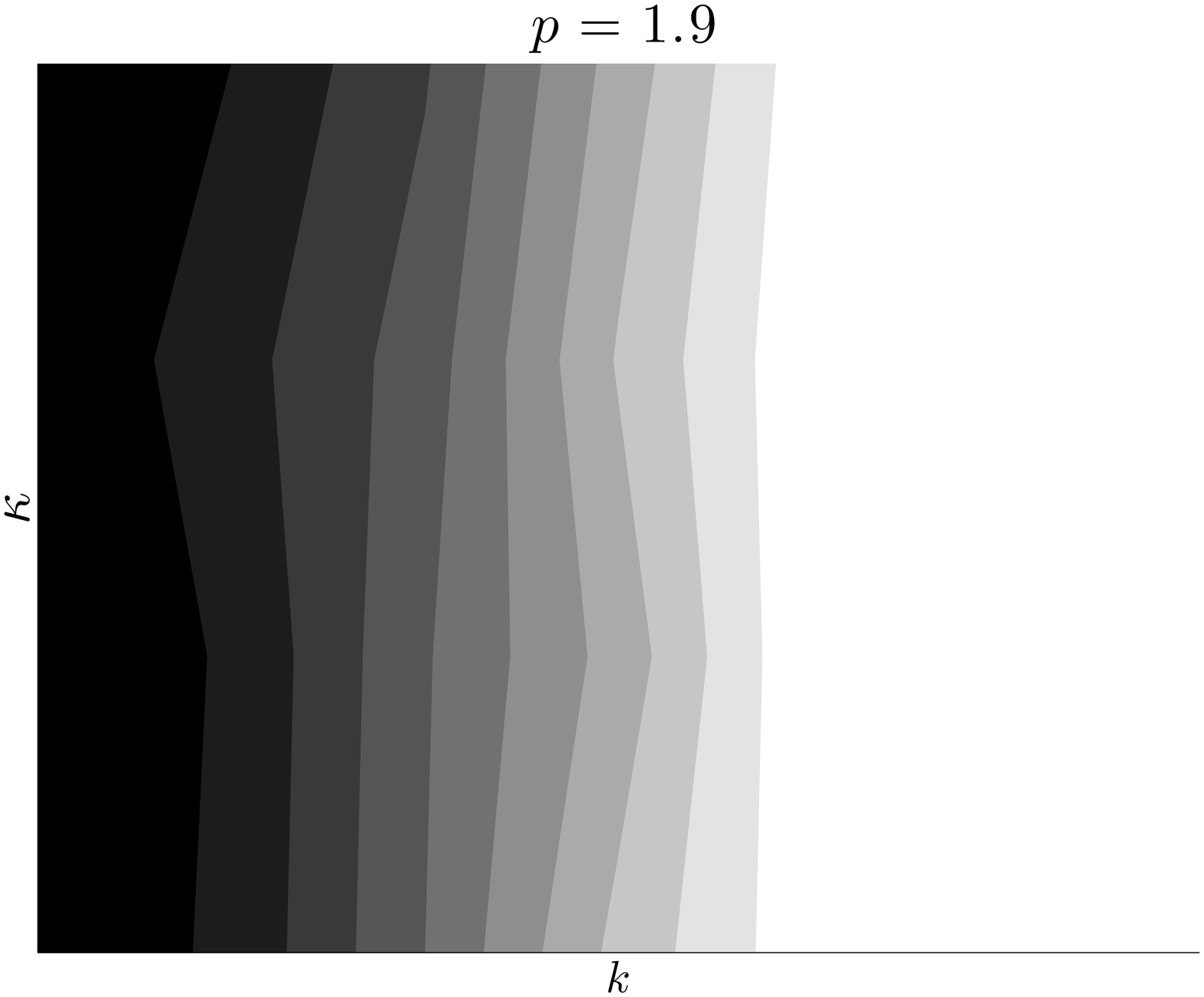}}\label{sub:1}}
\caption{Recovery rates for the greedy strategy developed in \cite{Quasi14} implemented with NR-IRLS used on the phase retrieval problem with Gaussian measurement vectors as above with $N = 80, m = 30$, and we use solutions $x^*$ with $\|x^*\|=1$ . Reconstruction is repeated 50 times for each signal and $k, \kappa$.}
\label{fig:PR}
\end{figure}

\begin{figure}[H]
\centering
\subfigure{
{\includegraphics[width=.48\textwidth]{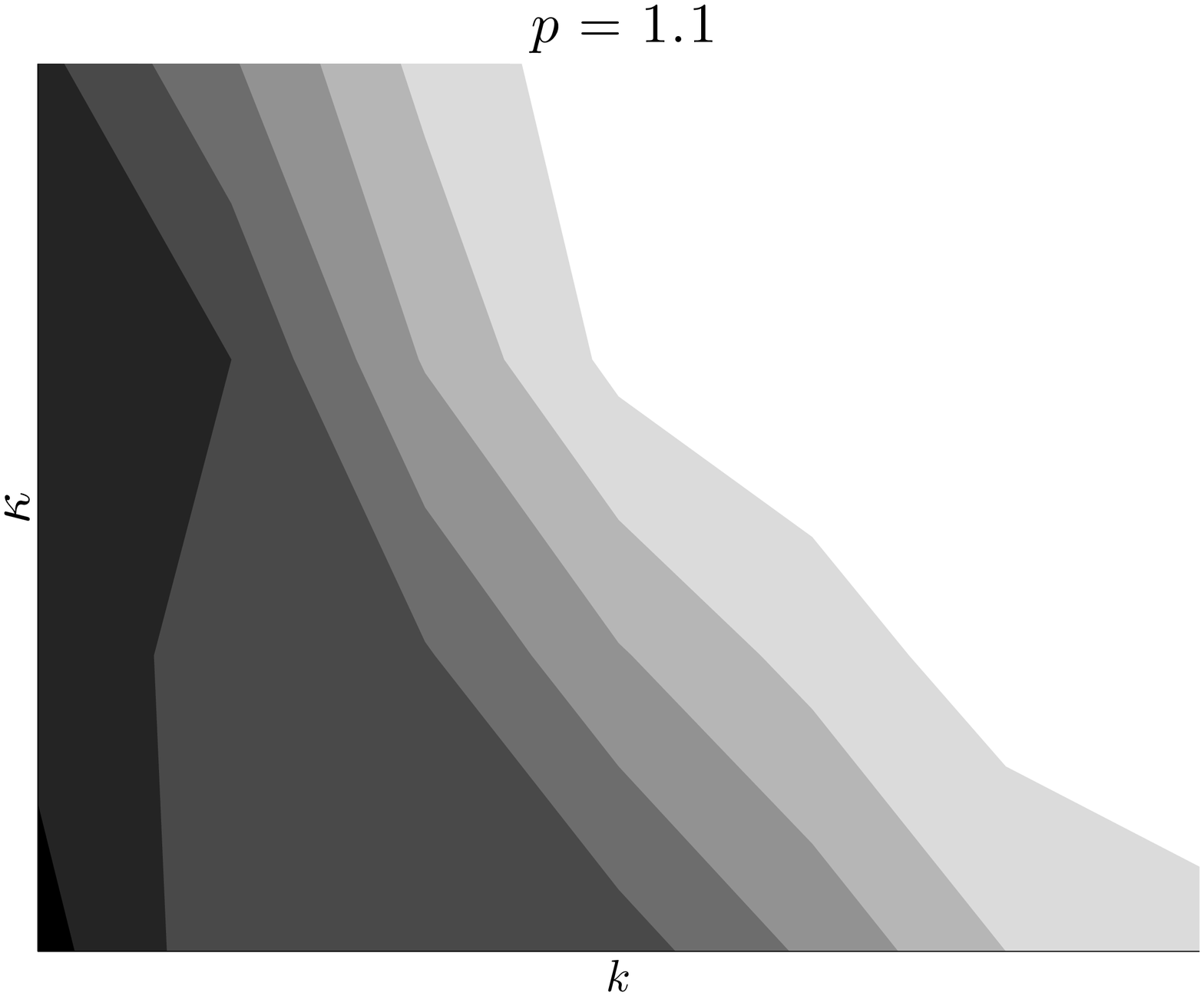}\label{sub:2}}}
\subfigure{
{\includegraphics[width=.48\textwidth]{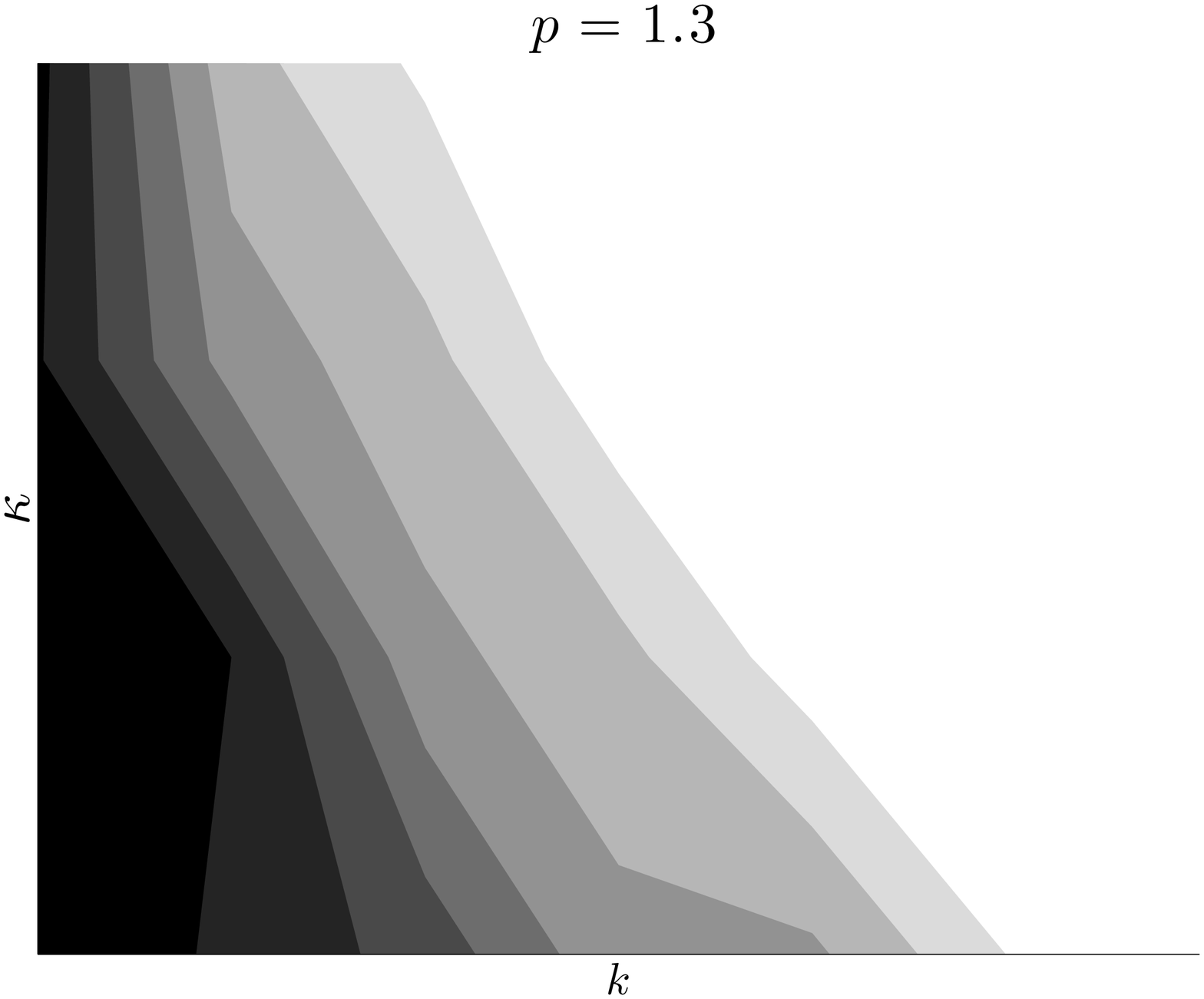}}\label{sub:1}}
\subfigure{
{\includegraphics[width=.48\textwidth]{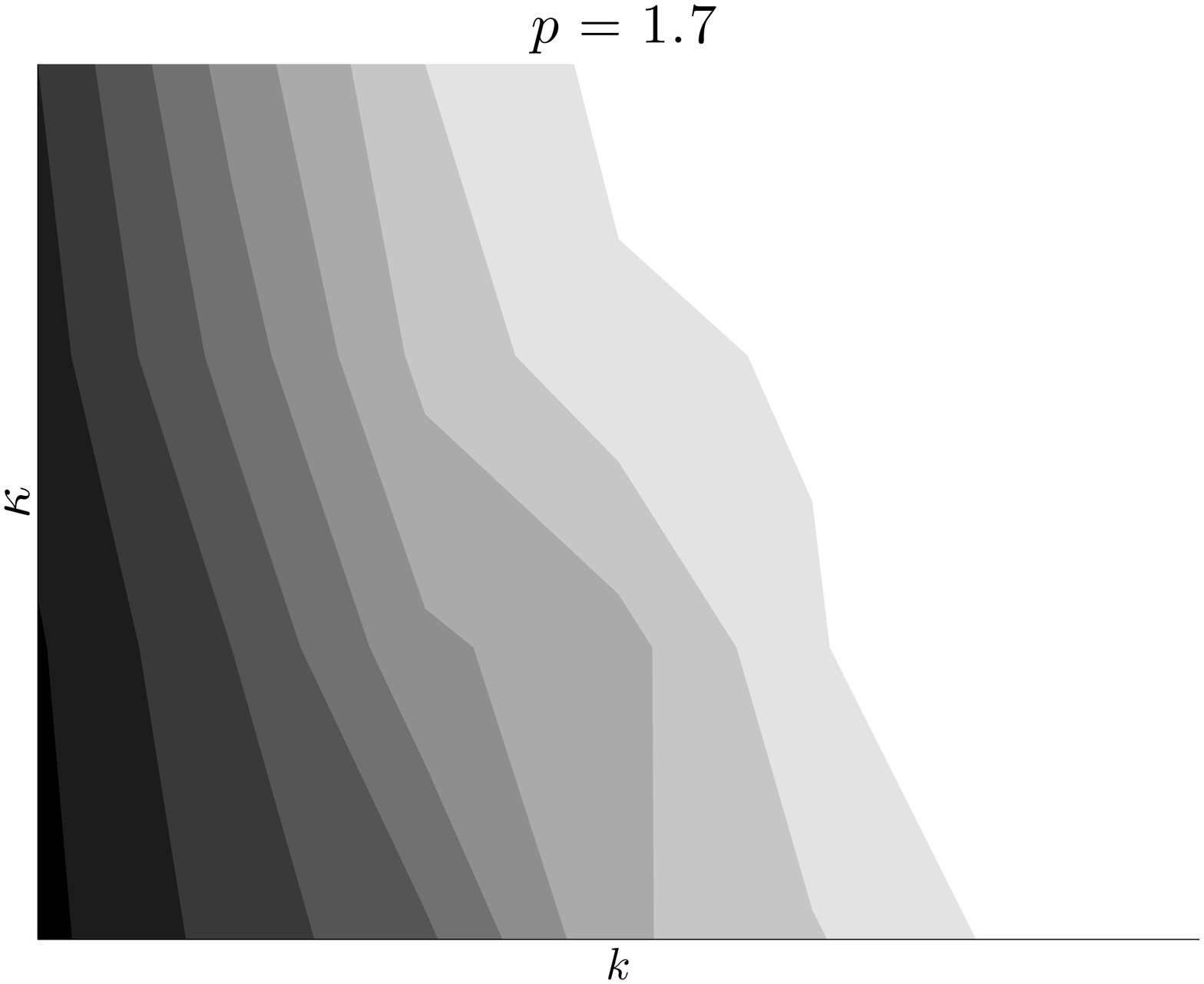}\label{sub:2}}}
\subfigure{ 
{\includegraphics[width=.48\textwidth]{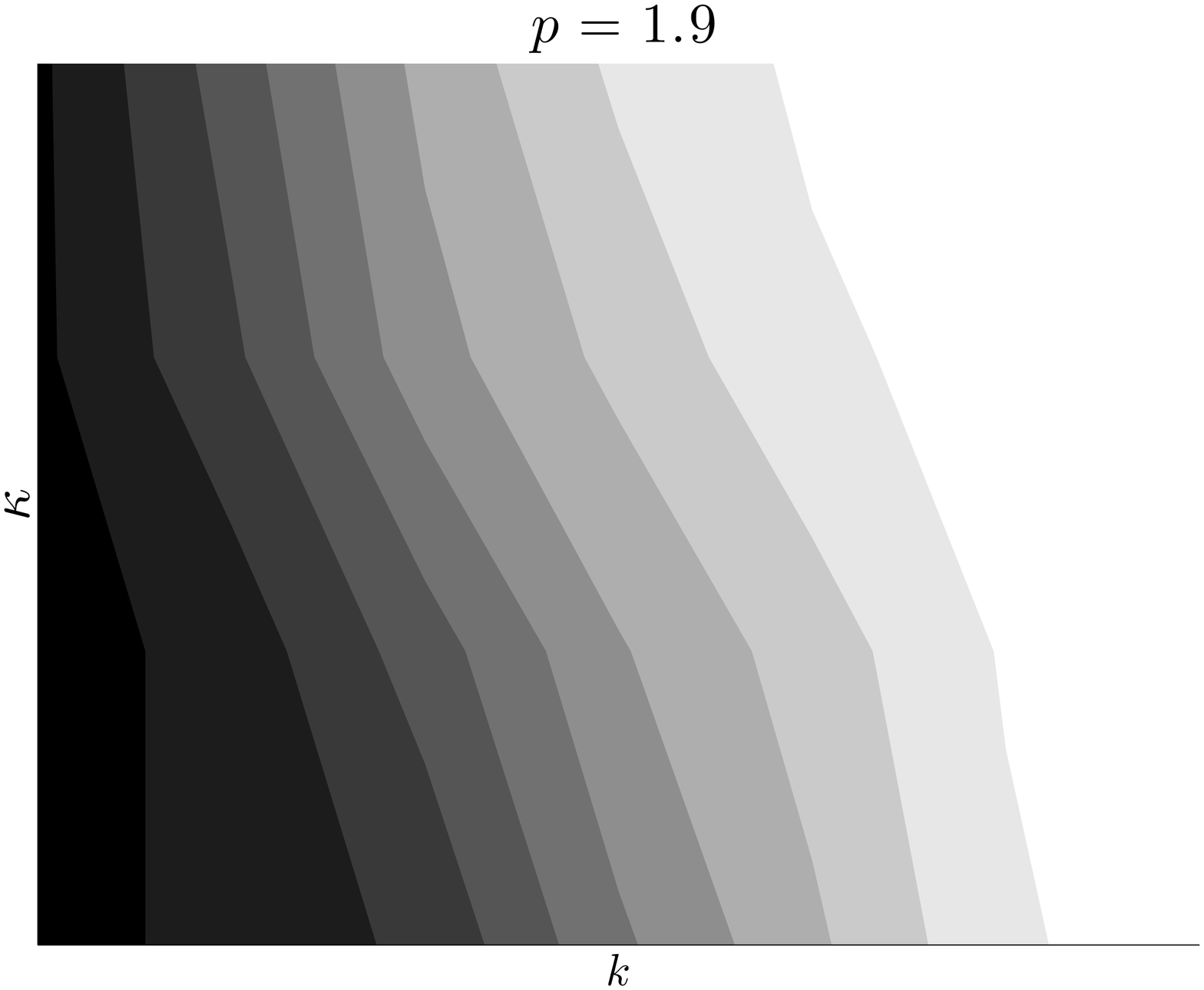}}\label{sub:1}}
\caption{Recovery rates for the greedy strategy developed in \cite{Quasi14} implemented with $lsqnonlin$ used on the phase retrieval problem with Gaussian measurement vectors as above with $N = 80, m = 30$, and we use solutions $x^*$ with $\|x^*\|=1$ . Reconstruction is repeated 50 times for each signal and $k, \kappa$.}
\label{fig:PR2}
\end{figure}

\subsection{Recovery from Data with Impulsive Noise Perturbation}
In the case of measurements $y$ additionally corrupted by noise, the optimal type of loss function for the nonlinear residual minimization has to be chosen depending on the particular kind of noise. The most common type is so-called white noise with Gaussian distribution, a continuous smooth perturbation of the original signal, where usual $\ell_2$-least squares methods are well suited and widely used. In contrast, impulsive noise are random occurrences of instantaneous signal perturbations taking the shape of spikes or pulses having random amplitude. This means the appearance of measurement distortions is sparse. For that  $\ell_1$-minimization of the nonlinear residual is a better choice than $\ell_2$-minimization.\\
In the following we want to consider noisy measurements in the context of the phase retrieval problem as described above and adopt the numerical test setting in great parts with only slight adaption mentioned below. The goal of this section is to examine the influence of the choice of $1\leq p\leq2$ on the recovery success rates for impulsive noise perturbations.\\
First of all we need to fix a possible statistical model for the impulsive noise. To this end we combine a binary-valued random sequence model of the
time of occurrence of impulsive noise with a continuous-valued random process model of impulse amplitude.\\
An important statistical process for modeling impulsive noise as an amplitude modulated
binary sequence is the Bernoulli-Gaussian process \cite{Vaseghi01}.
In a Bernoulli-Gaussian model of an impulsive noise process, the random
time of occurrence of the impulses is modeled by a binary Bernoulli process $B_{\alpha_p}$ with success probability $\alpha_p$ and the amplitude of the impulses is modeled by a Gaussian process $\mathcal{N}_{(0,1)}$ with mean $0$ and standard deviation $1$.
Having introduced a proper model for impulsive noise, we  apply the NR-IRLS Algorithm 2 on the noisy phase retrieval problem. We give a complete description of the measurement setting as follows.
As above we again fixed the dimension of the signal $N=80$, the number of measurements $m=30$ and i.i.d. Gaussian random vectors $a_i$, $i=1,\dots,m$. Moreover we created synthetic solutions $z^*$ with $\|z^*\|_{\ell_2}=1$ and respective sparsity $k \in \{1,2,3,5,7,9\}$. The vectors were constructed such that the nonincreasing rearrangement of the absolute value of their entries satisfies the decay rate $\kappa=0.5$, see the vector class $\mathcal D_\kappa$ in  \cite{Quasi14} for a precise definition. We generated impulsive Bernoulli-Gaussian noise with parameters $\alpha_p \in \{0.5,0.4,0.3,0.2,0.1,0.0\}$ respectively that was scaled to the norm of the measurements and added it to the originally generated measurement data itself.
For each of these parameter combinations we generate a set of 100 synthetic problems. As above we use the solution $z^*$, the expected sparse solution, to determine the success of the recovery and claim again successful reconstruction when the error is within a $5\%$ of the solution's norm. \\
The algorithmic settings are similar to above: $3k$ steps are performed by the greedy algorithm \cite[Algorithm 1]{Quasi14}, the maximum number iterations of NR-IRLS itself is allowed to reach 50. 
The regularization parameter $\omega$ is set to $100$ and for the execution of the convex minimization in each inner step we choose the MATLAB built-in function $fminunc$ with default settings and random starting points with a norm smaller or equal to the solutions norm.

\begin{figure}
\centering
\subfigure{
\includegraphics[width=.48\textwidth]{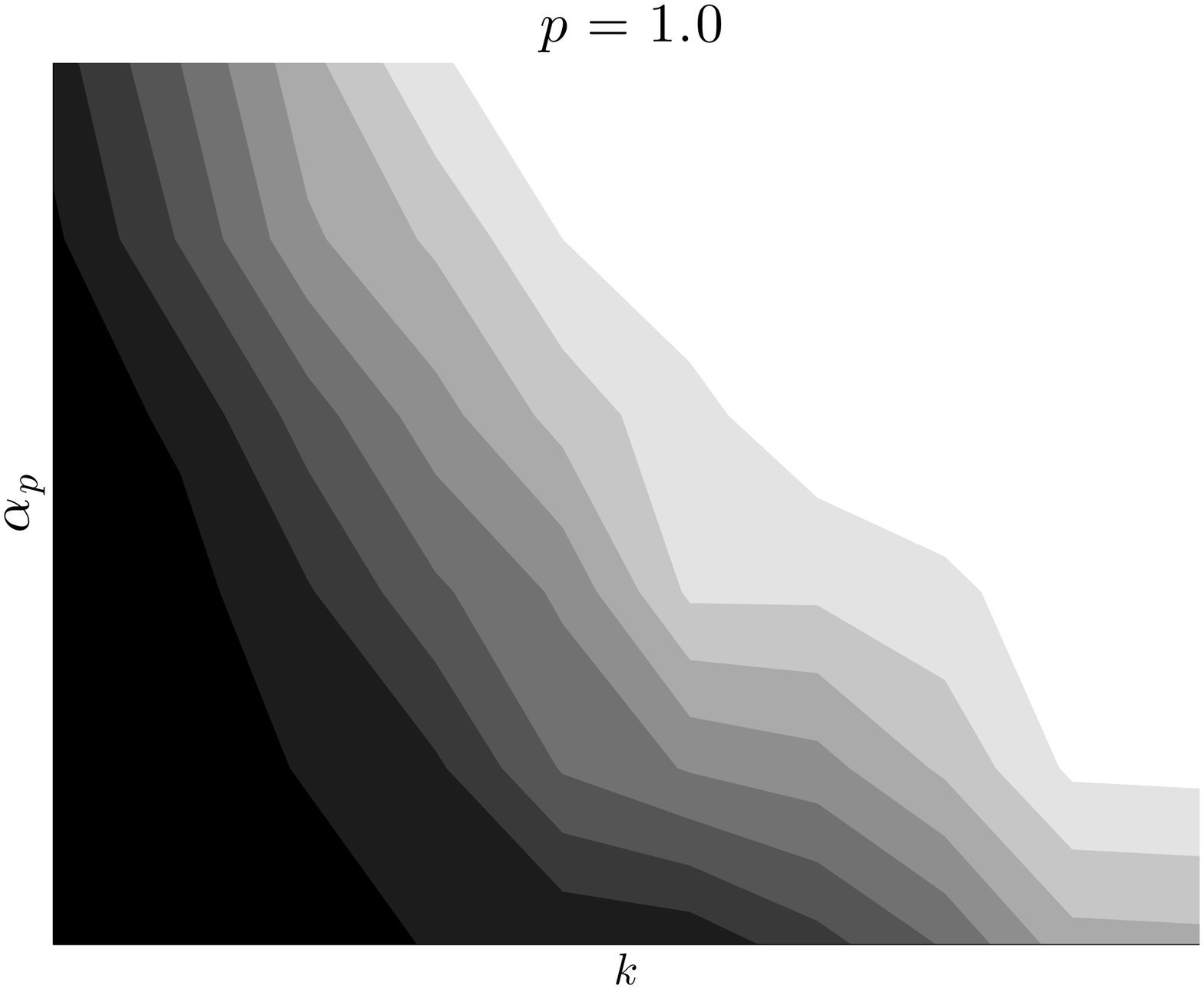}\label{sub:2}}
\subfigure{
\includegraphics[width=.48\textwidth]{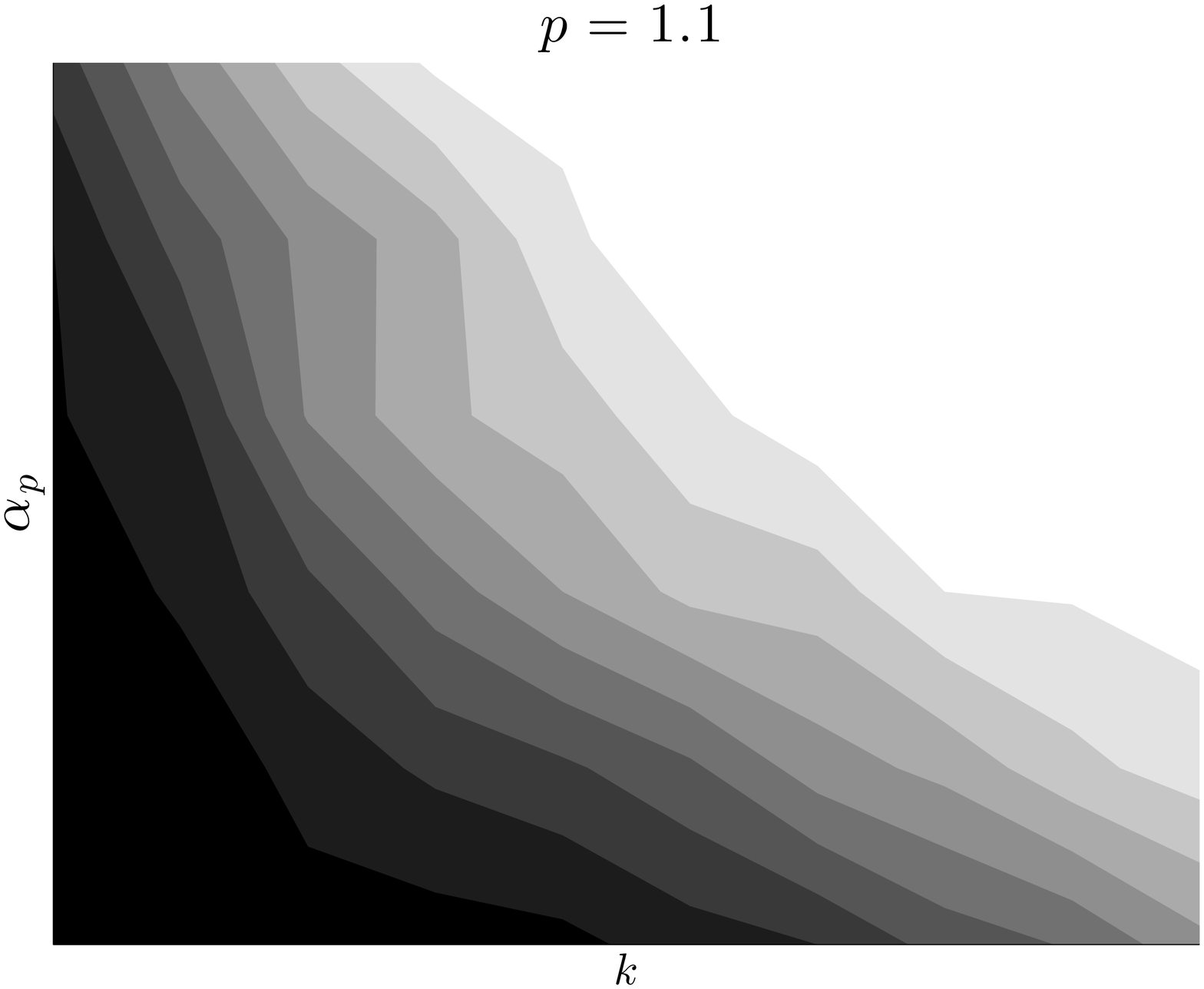}}\label{sub:1}
\subfigure{
\includegraphics[width=.48\textwidth]{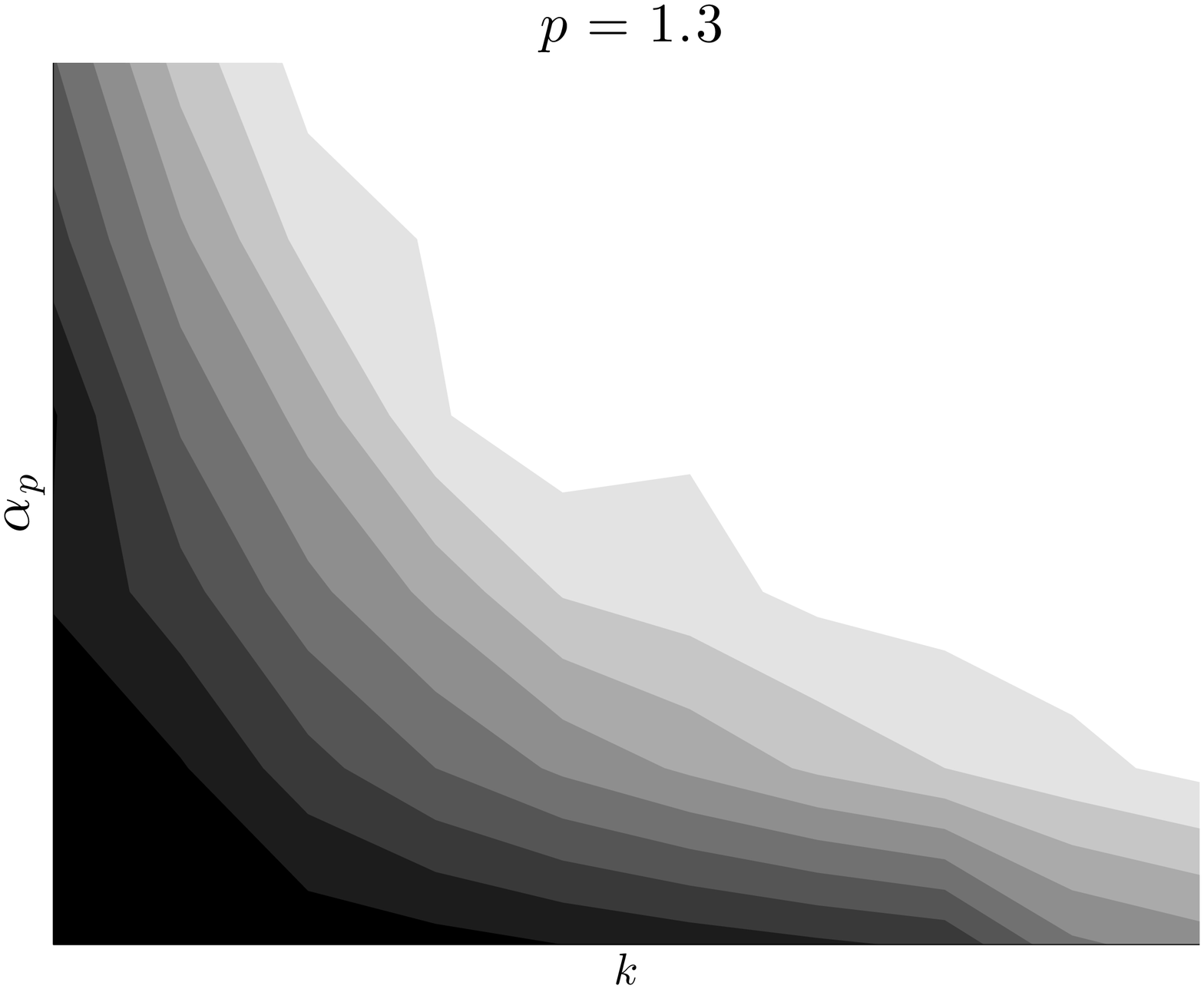}\label{sub:2}}
\subfigure{ 
\includegraphics[width=.48\textwidth]{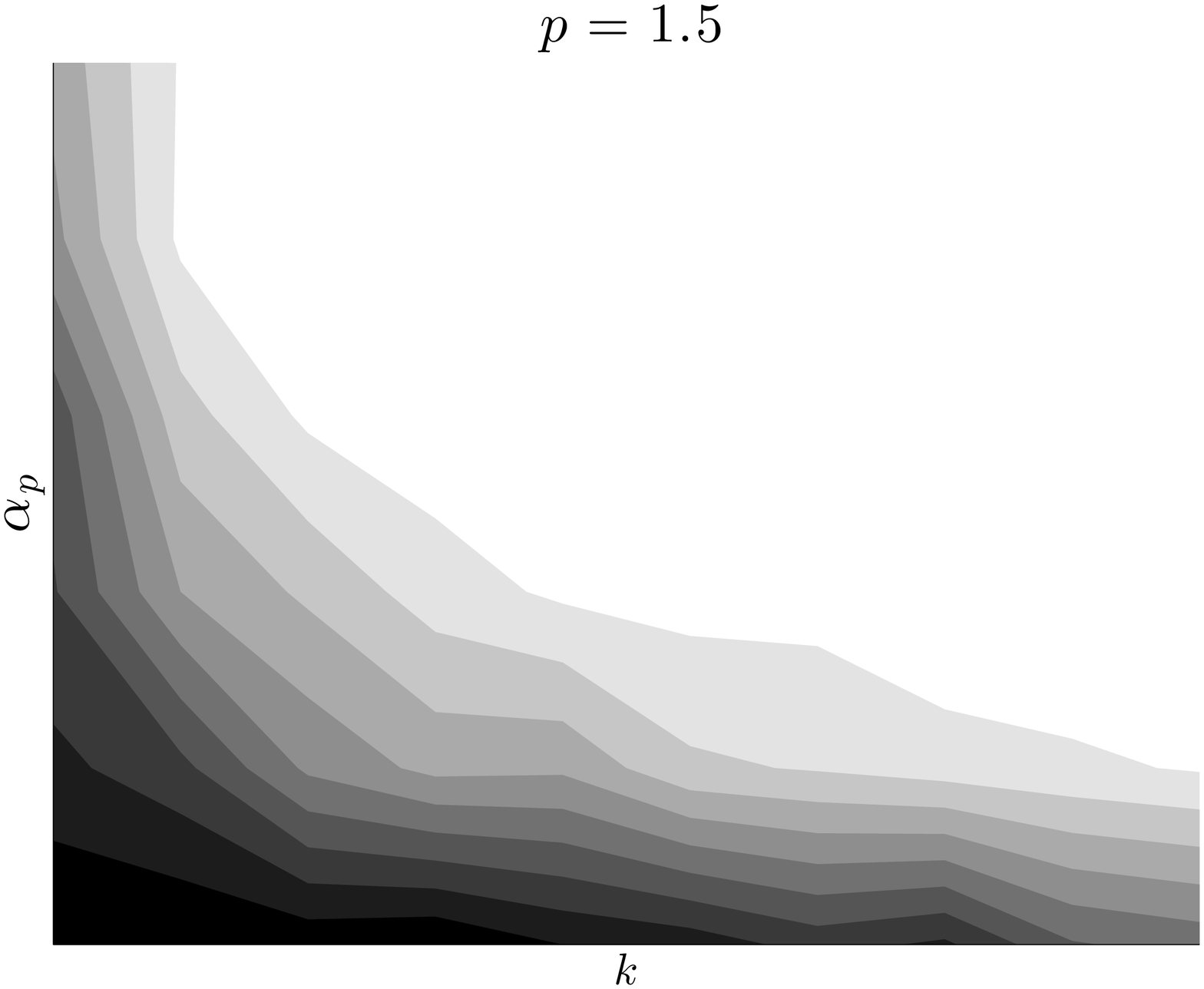}}\label{sub:1}
\subfigure{ 
\includegraphics[width=.48\textwidth]{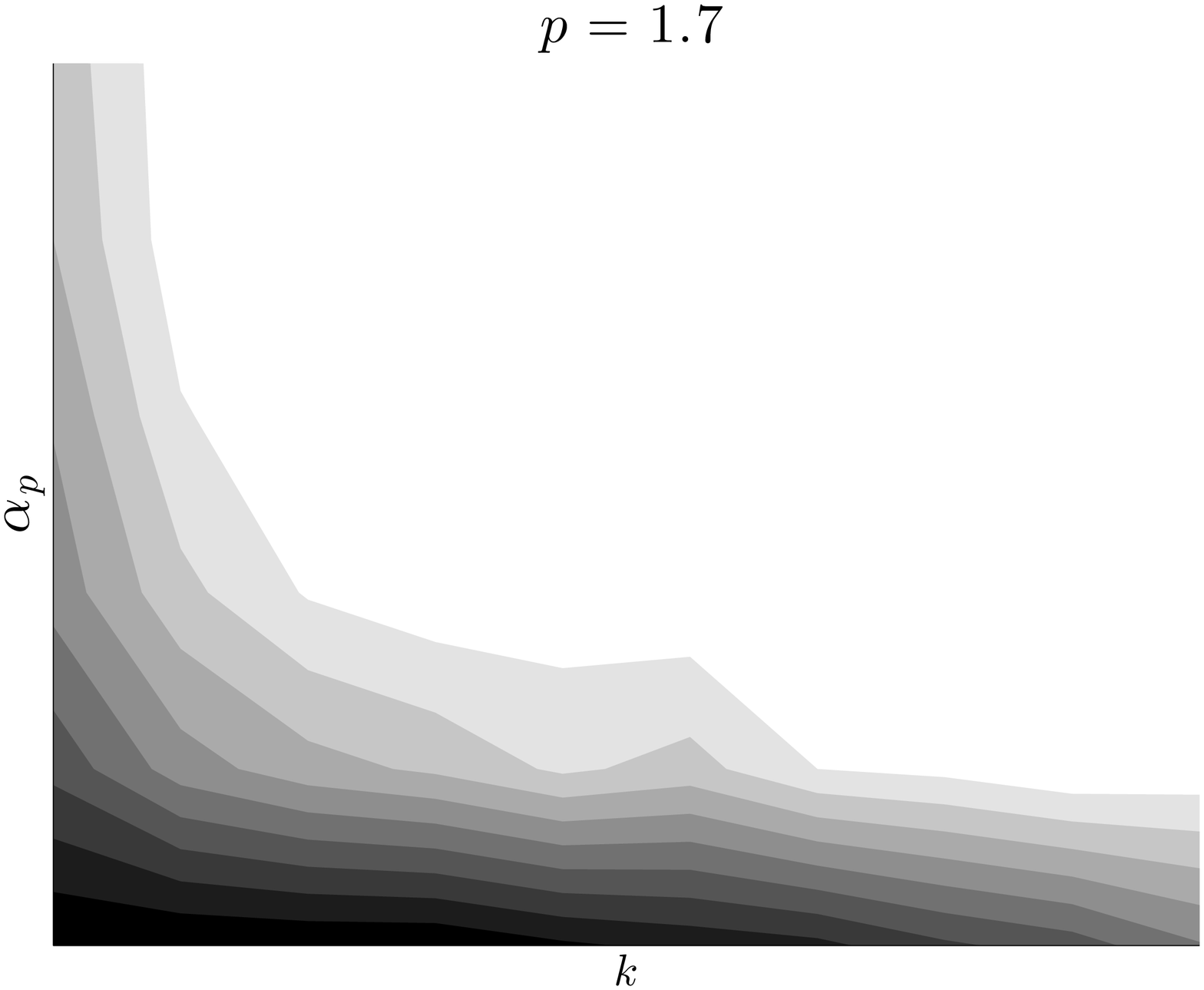}}\label{sub:1}
\subfigure{ 
\includegraphics[width=.48\textwidth]{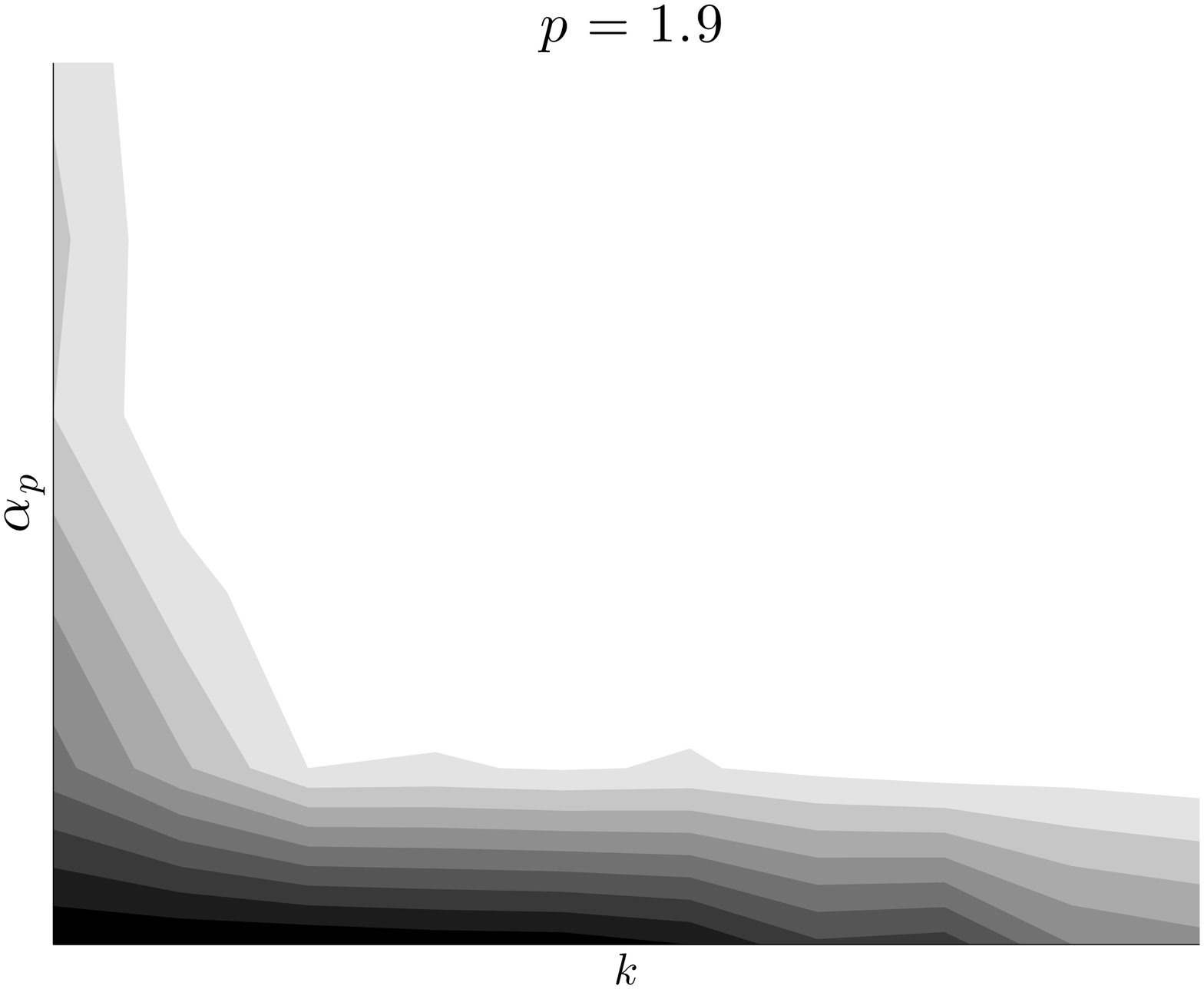}}\label{sub:1}
\subfigure{ 
\includegraphics[width=.48\textwidth]{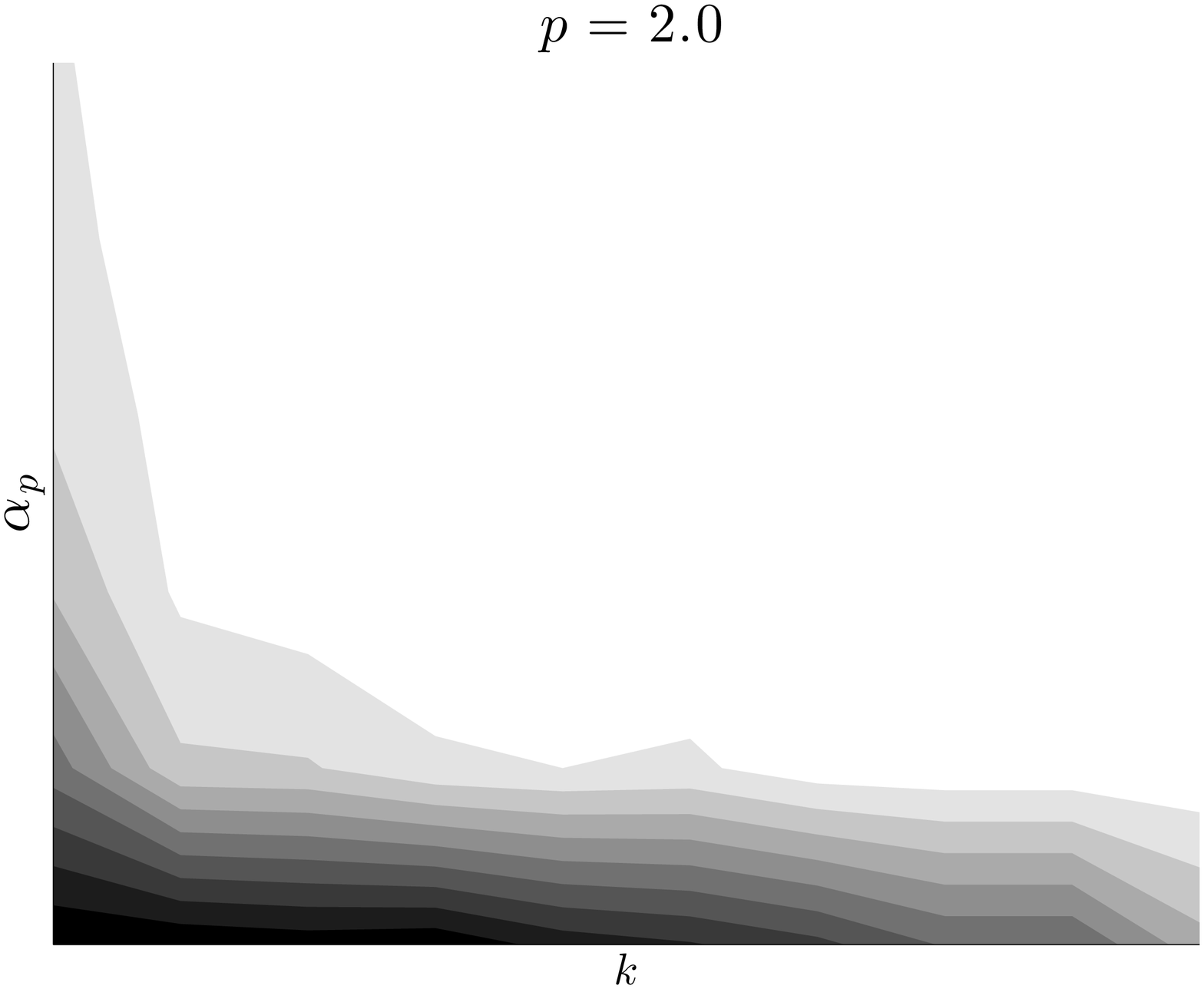}}\label{sub:1}
\caption{Recovery rates for the greedy strategy developed in \cite{Quasi14} implemented with NR-IRLS used on the phase retrieval problem with Gaussian measurement vectors as above with $N = 80, m = 30$, and we use solutions $x^*$ with $\|x^*\|=1$ . Reconstruction is repeated 50 times for each signal with sparsity $k$ and the particular noise perturbation as given above.}
\label{fig:PR4}
\end{figure}

The plots in Figure \ref{fig:PR4} show the recovery rates of \cite[Algorithm 1]{Quasi14} for sparse vectors from measurements of the type \thref{pert} affected by impulsive noise, implementing Algorithm 2 for performing the $\ell_p$-minimization for different values of $p$. 
As a short synthesis from the visual analysis of the phase transition diagrams, we can state that smaller values of $p$ and thus less smooth loss functions for residual minimization should clearly be preferred to standard $\ell_2$-least squares, as expected from the specific noise model. Moreover we note that for $p$ close to $1$ and a small number of vector entries recovery is still very robust also for strong perturbations with impulsive noise.

\section{Acknowledgement}
Juliane Sigl acknowledges the support of the  DFG$/$\"OAW- project FO 767/5-1.
She thanks Massimo Fornasier and Francesco Solombrino for the inspiring discussions on several issues appearing during the preparation of this work. 
\section{Appendix}

\subsection{Minimization of the functional $\mathcal{J}$ for a linear map $A$}

Let $\epsilon>0,w,y \in \mathbb{R}^m$ with $w>0$ be fixed. Moreover we first consider a general  map $A \in { \mathcal{C}^{2}}(\mathbb{R}^k,\mathbb{R}^m)$ and  bounded not explicitly requiring the linearity of $A$. (For the sake of an intuitive presentation we do not consider weaker assumptions on $A$ here.)\\
We examine the minimization of the functional $ \mathcal{J}(\cdot,\epsilon,w)$ in general first and then comment on the difference between the linear and nonlinear case.
 As we assume that the problem is smooth enough the minimizer $x^s$ of $ \mathcal{J}(\cdot,\epsilon,w)$ has to fulfill the necessary and sufficient conditions of a zero gradient
$$\nabla \mathcal{J}(x^s,\epsilon,w)=p\sum_{i=1}^m w_i(A(x^s)_i-y_i)\nabla A(x^s)_i=0,$$
and a positive definite Hessian
$$\nabla^2 \mathcal{J}(x^s,\epsilon,w)=p \sum_{i=1}^m w_i \left[ \nabla A(x^s)_i \nabla A(x^s)_i^T+(A(x^s)_i-y_i)\nabla^2 A(x^s)_i \right] >0.$$
%
In the linear case, where $A(x)=A$, $\nabla^2 \mathcal{J}(x,\epsilon,w)>0$ for every $x$ and $\epsilon, w$ fixed as soon as $\nabla A(x)=A^T$ has full rank $k$. This is due to the fact that in this special case $\nabla^2 A(x)_i=0$ for $i=1,\dots,m$ and $\sum_{i=1}^m w_i  \nabla A(x)_i \nabla A(x)_i^T>0$ for $k$ linear independent vectors $\nabla A(x)_i$. \\
For the sake of completeness we  prove the latter property in more generality here.
\begin{lemma} \thlabel{full}
If $a_i \in \mathbb{R}^k$, $i=1,\dots,m$ and $A=\left[a_1|\dots,|a_i|\dots|a_m\right]$ has full rank $k\leq m$, then for all $w=(w_1,\dots,w_i,\dots,w_m)$ with $w_i >0$ for $i=1,\dots,m$ it holds that 
$$\sum\limits_{i=1}^m w_ia_ia_i^T >0$$
is a positive definite matrix. 
\end{lemma}
Proof: First observe that $\sum \limits_{i=1}^mw_ia_ia_i^T \geq0$ is obvious because
$$0\leq h^T\left( \sum \limits_{i=1}^mw_ia_ia_i^T\right)h= \sum \limits_{i=1}^mw_i (a_i^Th)^2  \text{ for all  }h \in \mathbb{R}^k.$$
Such a quantity is vanishing $0=\sum \limits_{i=1}^mw_i (a_i^Th)^2$ iff $a_i^Th=0$ for all $i=1,\dots,m$. \\
But as $A$ has full rank there exist $a_{i_1},\dots, a_{i_k}$ linear independent vectors in $\mathbb{R}^k$ for which $a_{i_l}^Th=0$ for all $l=1,\dots,k$ and this implies $h\equiv0$.
$\blacksquare$\\

Hence in the linear case the functional of \thref{J} is strictly convex in the variable $x$. This implies that every stationary point is a global minimizer and this optimization problem can be solved efficiently.

\bibliography{Literature}
\end{document}